\documentclass[accepted]{uai2022} 

\usepackage[american]{babel}

\usepackage[round]{natbib}
\usepackage{booktabs} 
\usepackage{tikz} 

\title{$\ell_{\infty}$-Bounds of the MLE in the BTL Model under General
Comparison Graphs}

\author[1]{\href{mailto:<wanshanl@andrew.cmu.edu>?Subject=Your UAI 2022
paper}{Wanshan~Li}{}}
\author[1]{\href{mailto:<sshrotri@andrew.cmu.edu>?Subject=Your UAI 2022
paper}{Shamindra~Shrotriya}{}}
\author[1]{\href{mailto:<arinaldo@andrew.cmu.edu>?Subject=Your UAI 2022
paper}{Alessandro~Rinaldo}{}}
\affil[1]{%
    Department of Statistics \& Data Science \\
    Carnegie Mellon University \\
    Pittsburgh, Pennsylvania, USA }

\usepackage{etoolbox}
\newtoggle{uai2022version} 
\newtoggle{arxiv2022version} 

\toggletrue{arxiv2022version} 

\iftoggle{arxiv2022version}{
  \pagestyle{numbered}
  }{
}

\usepackage{tcolorbox}
\usepackage{float}

\usepackage{csquotes}
\newenvironment{itquote}
{\begin{quote}\itshape} {\end{quote}}

\newcommand{\btldima}{n}
\newcommand{\btldimb}{N_{\text{comp}}}
\newcommand{\btldimc}{\theta}
\newcommand{\btl}{BTL}

\usepackage{inconsolata}

\newcommand{\specialcell}[2][c]{%
  \begin{tabular}[#1]{@{}c@{}}#2\end{tabular}}

\usepackage{graphicx} 
\graphicspath{{figures/}} 

\usepackage{000_custom_macros_comb}

\iftoggle{arxiv2022version}{
  \usepackage{hyperref}       
\hypersetup{colorlinks = true, linkcolor = blue, anchorcolor = blue, citecolor =
            blue, filecolor = blue, urlcolor = blue } }{
  \usepackage{xr}
  \usepackage{xr-hyper}
  \usepackage[breaklinks]{hyperref}[2019/11/10]
}

\usepackage[labelfont=bf]{caption}
\usepackage{url}
\usepackage[capitalize,sort,compress]{cleveref}

\iftoggle{arxiv2022version}{
}{
\externaldocument[]{li_649-supp}

\makeatletter
\newcommand*{\addFileDependency}[1]{
  \typeout{(#1)}
  \@addtofilelist{#1} \IfFileExists{#1}{}{\typeout{No file #1.}} }
\makeatother

\newcommand*{\myexternaldocument}[1]{%
    \externaldocument{#1}%
    \addFileDependency{#1.tex}%
    \addFileDependency{#1.aux}%
}

\myexternaldocument{li_649-supp}
}

\begin{document}
\maketitle

\begin{abstract}
  The Bradley-Terry-Luce (\btl{}) model is a popular statistical approach for
  estimating the global ranking of a collection of items using pairwise
  comparisons. To ensure accurate ranking, it is essential to obtain precise
  estimates of the model parameters in the $\ell_{\infty}$-loss. The difficulty
  of this task depends crucially on the topology of the pairwise comparison
  graph over the given items. However, beyond very few well-studied cases, such
  as the complete and Erd\"os-R\'enyi comparison graphs, little is known about
  the performance of the maximum likelihood estimator (MLE) of the \btl{} model
  parameters in the $\ell_{\infty}$-loss under more general graph topologies.
  In this paper, we derive novel, general upper bounds on the $\ell_{\infty}$
  estimation error of the \btl{} MLE that depend explicitly on the algebraic
  connectivity of the comparison graph, the maximal performance gap across items
  and the sample complexity. We demonstrate that the derived bounds perform well
  and in some cases are sharper compared to known results obtained using
  different loss functions and more restricted assumptions and graph topologies.
  We carefully compare our results to \citet{yan2012sparsecompbtl}, which is
  closest in spirit to our work. We further provide minimax lower bounds under
  $\ell_{\infty}$-error that nearly match the upper bounds over a class of
  sufficiently regular graph topologies. Finally, we study the implications of
  our $\ell_{\infty}$-bounds for efficient (offline) tournament design. We
  illustrate and discuss our findings through various examples and simulations.
\end{abstract}

\section{Introduction}\label{sec:introduction}

Simultaneous or `global' ranking of a set of items is a practical problem that
arises naturally in a variety of domains. For example, one may wish to ascertain
a `best player' or `best team' in a given sports league. Designing a principled
statistical approach to global ranking of items is challenging due to data
limitations and complex domain-specific relationships between the underlying
items to be ranked.

A popular and practicable solution  to estimating global ranking  is to utilize
pairwise comparison information across the items to be ranked, which is easily
accessible across many application domains. The \btl{} model
\citep{bradley1952rank,Luce59} is a popular statistical model for pairwise
comparison data. A similar model was also originally studied in
\cite{Zermelo1929}. The continued practical and theoretical interest in the
\btl{}  model stems from its relatively simple parametric form which provides a
good balance between interpretability and tractability for theoretical analysis.
The \btl{} model is domain-agnostic, making it an ideal benchmarking tool across
a variety of ranking applications \eg sports analytics
\citep{FaT1994,MaV2012,CMV2012}, and bibliometrics \citep{St1994, Va2016}.

Formally, we can describe the \btl{}  model as follows. Suppose that we have $n$
distinct items, each with a (fixed but unobserved) positive strength or
preference score $w^{*}_{i}$, $i \in [n]$, quantifying item $i$'s propensity to
beat other items in pairwise comparisons. The \btl{} model assumes that the
comparisons between different pairs are independent and the outcomes of
comparisons between any given pair, say item $i$ and item $j$, are \iid
Bernoulli random variables, with \textit{winning probability} $p_{ij}$, defined
as
\begin{equation}\label{neqn:bradley_terry_prob_succ}
  p_{ij}
  \defined \Prb{i \text{ beats } j}
  \defined \frac{w^{*}_{i}}{w^{*}_{i} + w^{*}_j}, \: \forall \; i,j \in [n].
\end{equation}
A common reparametrization is to set, for each $i$,  $w^{*}_{i} =
  \exp(\theta^*_i)$, where $\boldsymbol{\theta}^* \defined (\theta^*_{1},
  \ldots, \theta^*_{n})^{\top} \in \reals^{n}$. By convention, we assume that
$\sum_{i \in [n]} \theta^*_i = 0$ for parameter identifiability.

From a theoretical perspective, much attention in the \btl{} literature has been
paid to two popular estimators, namely the maximum likelihood estimator (MLE)
and the spectral method \citep{jain2020spectralmethodscarcedata}. Recently,
\cite{chen2020partialtopkranking} show that the MLE attains a sharper minimax
rate of the Hamming top-$k$ loss compared to the spectral method. In this paper,
we thus focus on the MLE, which we formally define later in
\Cref{sec:upper-bounds}.

\noindent{\bf General pairwise comparison
  graphs}\label{subsec:optimality-mle-gen-topology}

Given $n$ items to be compared, the pairwise comparison scheme among them can be
expressed through an undirected simple graph $\mclG(V, E)$, where the vertex set
$V \defined [n]$ and the edge set $E \defined \{(i,j): i \text{ and } j \text{
    are compared }\}$ is determined by the comparison scheme.  Correspondingly, if
we define the directed edge set as $E_d\defined \{(i,j,k): (i \text{ beats } j)
  \text{ $k$ times}\}$, then the induced directed simple graph $\mclG(V,E_d)$ is
called a \textit{directed} comparison graph. It is a classical result
\citep{ford1957,simons1999,hunter2004mm} that the \btl{} model is identifiable
if and only if $\mclG(V,E)$ is connected, and the MLE of the model parameters
exists and is consistent if and only if $\mclG(V,E_d)$ is strongly connected.
Henceforth, \textit{comparison graph} refers to the undirected pairwise
comparison graph.

Typically one is interested in getting sharp bounds for the estimation risk,
which could be based on a norm-induced metric $ \|\hat{\bbrtheta} -
  \bbrtheta^*\|_{p}$ or a ranking metric, \eg, Kendall’s tau distance
\citep{kendall1938}. What makes risk analysis of \btl{} model estimators
particularly challenging is a combination of the type of estimation risk loss
considered, and the assumptions on the topology of $\mclG(V,E)$.

\noindent{\bf Core questions of
  interest}\label{subsec:core-question-of-interest}

Among all the metrics measuring uncertainty of estimators of \btl{} parameters,
the $\ell_\infty$-loss directly connects with ranking metrics, \eg binary and
Hamming top-$k$ (partial) ranking loss \citep[see,
  \eg][]{chen2019spectralregmletopk,chen2020partialtopkranking}.

It is thus natural to study the MLE for the \btl{} parameters in the
$\ell_{\infty}$-loss, to better understand the risk optimality of the MLE and
further justify its use for practical global and partial ranking problems. In
this spirit, \cite{yan2012sparsecompbtl} focus specifically on proving
$\ell_{\infty}$-error bounds for the \btl{} MLE for general comparison graphs.
However, a notable limitation in their setting is that they impose a strictly
dense comparison graph assumption, which may be impractical in many real world
applications. This leaves a gap in the literature, summarized in the following
questions:

\begin{tcolorbox}
  \begin{itquote}
    \textbf{Core questions:} For the \btl{} model, how does the MLE perform with
    respect to the $\ell_{\infty}$ loss, under much weaker assumptions on the
    pairwise comparison graph compared to \citet{yan2012sparsecompbtl}? That is,
    assuming only that the comparison graph is connected. Moreover, what are the
    implications of such bounds in applications?
  \end{itquote}
\end{tcolorbox}

Providing a sharp analysis to these questions with a detailed comparison to
recent theoretical results in the \btl{} literature motivates our work in this
paper.

\noindent{\bf Relevant and related literature}

We give a brief overview of the work that addresses the challenge of comparison
graph topology in ranking. When the comparison graph is a complete graph,
\cite{simons1999} give a high-probability upper bound for the $\ell_{\infty}$
loss, \ie, $\|\hat{\bbrtheta} - \bbrtheta^*\|_{\infty}$ and obtain the
asymptotic distribution of the MLE. In the setting where the comparison graph
follows the Erd\"os-R\'enyi graph model,
\cite{chen2015spectralmletopkpairwisecomparison},
\cite{chen2019spectralregmletopk}, \cite{chen2020partialtopkranking} and
\cite{han2020asymptoticsparsebradleyterry} derive high-probability upper bounds
for the $\ell_{\infty}$ loss. Moreover,
\cite{chen2019spectralregmletopk}
show that both MLE and spectral method are minimax optimal in terms of the
binary top-$k$ ranking loss, \ie, whether the items with the highest $k$ out of
$n$ preference scores are perfectly identified;
\cite{chen2020partialtopkranking} consider a Hamming Loss for top-$k$ items and
show that the MLE is minimax optimal compared to the spectral method with
differences arising in constant factors.

For a broader class of comparison graphs beyond complete and Erd\"os-R\'enyi
graph, researchers have studied the explicit dependence of the estimation risk
on graph topology. In particular, \cite{yan2012sparsecompbtl} give a
high-probability upper bound for the $\ell_{\infty}$-loss for relatively dense
graphs.
\cite{hajek2014minimaxinferencepartialrank,shah2015estimationfrompairwisecomps}
give a high probability upper bound for the $\ell_2$ or Euclidean loss
$\|\hat{\bbrtheta} - \bbrtheta^*\|_{2}$, establish upper and lower bounds of
$\mbbE{\|\hat{\bbrtheta} - \bbrtheta^*\|_{2}}$ and show the minimax optimality
of the constraint MLE across a wide range of graph topologies. Recently,
\cite{agarwal2018acceleratedspectralranking} give sharp upper bounds for a novel
spectral method in the $\ell_1$-loss $\|\hat{\bbrpi} - \bbrpi^*\|_{1}$ for
$\bbrpi^* = \bfw^*/\|\bfw^*\|_1$ instead of $\bbrtheta^*$.
\cite{hendrickx2019graphresistance, hendrickx2020minimaxpairwisebtl} propose a
weighted least square method to estimate $\bfw^*$ and prove a sharp upper bound
for their estimator in $\mathbb{E}[\sin^2(\hat{\bfw},\bfw^*)]$ or equivalently
in $\mathbb{E}\|\hat{\bfw}/\|\hat{\bfw}\|_2 - \bfw^*/\|\bfw^*\|_2\|^2_2$, in the
sense that this upper bound matches a instance-wise lower bound up to constant
factors. \looseness=-1

\noindent{\bf Contributions}

Our contributions in this paper are fourfold and are summarized as follows:
\bitems
\item \textbf{Upper bounds:} We derive a novel upper bound for the
$\ell_{\infty}$-error of the regularized MLE in \btl{} model allowing for
general graph topology. Our upper bounds hold under minimal assumptions on graph
topologies, \ie, assuming only that the comparison graph is connected. Given
such generality, we show our $\ell_\infty$ bound is tighter than existing
results under a broad range of graph topologies, and works well in general. In
particular, we carefully compare our work analytically and in simulation to
\citet{yan2012sparsecompbtl}, which is closest in spirit to our work.

A minor corollary of our techniques results in the state of the art
$\ell_{2}$-loss bounds for the Erd\"os-R\'enyi graph.
\item \textbf{Lower bounds:} We derive minimax lower bounds for \btl{} parameter
estimation in $\ell_{\infty}$-loss. We analyze specific graph topologies
satisfying certain regularity connectivity conditions under which the \btl{} MLE
is nearly minimax optimal.
\item \textbf{Implications for tournament design:} We show that the \btl{} MLE
in $\ell_{\infty}$-loss satisfies a unique subadditivity property, and how our
$\ell_{\infty}$ bounds can exploit this property for efficient (offline)
tournament design.
\item \textbf{Extension to the unregularized \btl{} model:} We also extend our
upper bounds under $\ell_{\infty}$-loss to the unregularized (`vanilla') \btl{}
MLE, which is also frequently used in practice.
\eitems
Due to the more complicated form of the vanilla \btl{} MLE upper bounds and
space limitations, we present these analagous results and their proofs
separately in \Cref{sec:vanilla-mle}. Henceforth, MLE refers to the
regularized \btl{} MLE unless stated otherwise. In addition to our theoretical
contributions a core aspect throughout our paper is to emphasize the
interpretability of our results, the associated assumptions, and implications
for practical ranking tasks.

\noindent {\bf Organization of the paper}

The rest of the paper is organized as follows. In \Cref{sec:upper-bounds}, we
present our main results for the upper bound in \Cref{nthm:thm1} and an
interpretation of the key components of the bound. In \Cref{sec:lower-bounds},
we discuss minimax lower bounds using the $\ell_{\infty}$ risk loss in
\Cref{nthm:thm2-lb}. In \Cref{sec:implications-of-work}, we show some practical
implications of our results in efficient tournament design from a ranking
perspective. In \Cref{sec:simulations}, we conduct extensive numerical
simulations to validate the optimality of our bounds compared to related results
in the literature.

\noindent {\bf Notation}

We typically use lowercase for scalars, \eg, $(x, y, z, \ldots)$, boldface
lowercase for vectors, \eg, $(\bfx, \bfy, \bfz, \ldots)$, and boldface uppercase
for matrices, \eg $(\bfX, \bfY, \bfZ, \ldots)$. We denote the finite set
$\theseta{1, \ldots, n}$ by $[n]$. For asymptotics, we denote $x_n\lesssim y_n$
or $x_n = O(y_n)$ and $u_n\gtrsim v_n$ or $u_n = \Omega(v_n)$ if $\forall n$,
$x_n \leq c_1 y_n$ and $u_n\geq c_2 v_n$ for some constants $c_1,c_2>0$. We
denote $\bfe_i$ as a vector whose entries are all $0$ except that the $i$-th
entry is $1$. $a_n = o(b_n)$ means $a_n/b_n\rightarrow 0$ as $n\rightarrow
  \infty$ and conversely, $a_n = \omega(b_n)$ means $b_n/a_n\rightarrow 0$ as
$n\rightarrow \infty$. We denote $\textbf{1}_n \in \reals^{n}$ to be a vector of
ones.

\section{Upper bounds}\label{sec:upper-bounds}

Recall that given $n$ items to be compared, the comparison scheme among them
defines the comparison graph $\mclG(V, E)$, where $V = [n]$ and $E = \{(i,j): i
  \text{ and } j \text{ are compared }\}$. We denote the corresponding adjacency
matrix as $A\in \mathbb{R}^{n\times n}$, and its $(i,j)^{\text{th}}$ entry is
$A_{ij} \defined 1\{(i,j)\in E\}$. The associated (unnormalized) graph Laplacian
is the symmetric, positive-semidefinite matrix   $\mclL_{\bfA}\defined\bfD -
  \bfA$, where $\bfD = \mathrm{diag}(n_1,\ldots, n_n)$, with
$n_{i}\defined\sum_{j=1}^n A_{ij}$ the degree of node $i$. It is well known that
the smallest eigenvalue of $\mclL_{A}$ is $0$ with an eigenvector
$\textbf{1}_n$. Let $\lambda_2(\mclL_\bfA)$ be the second smallest eigenvalue of
$\mclL_\bfA$, known as the algebraic connectivity of $\mathcal{G}$
\citep{laplacian2004}, then $\mclG$ is connected if and only if
$\lambda_2(\mclL_\bfA) >0$. Following the standard in the \btl{} literature we
assume a that for each edge $(i,j)$ of the comparison graph, the corresponding
items $i$ and $j$ are compared $L$ times, each leading to an independent outcome
$y_{ij}^{(l)}\in \{0,1\}$, where $l \in [L]$. If pairs are compared different
number of times, we take $L$ to be the smallest number of pairwise comparisons
over the edge set, as a worst-case scenario. The corresponding sample averages
are denoted with $\bar{y}_{ij} = \frac{1}{L}\sum_{l = 1}^L y_{ij}^{(l)}$ and are
sufficient statistics for the model parameters. The $\ell_{2}$-regularized MLE
is defined as
\begin{equation}\label{eq:reg.mle}
  \hat{\bbrtheta}_{\rho} = \argmin_{{\bf 1}_n^\top \bbrtheta = 0} \ell_{\rho}(\bbrtheta;\bfy),\  \ell_{\rho}(\bbrtheta;\bfy) =\ell(\bbrtheta;\bfy) + \frac{\rho}{2} \|\bbrtheta\|_2^2,
\end{equation}
where $\ell(\bbrtheta;\bfy)$ is the negative log-likelihood, given by
\begin{equation}
  \begin{split}
    \ell(\bbrtheta;\bfy) \defined - &\sum_{1\leq i<j\leq n} A_{ij}\lbrace\bar{y}_{ij}\log {\psi(\theta_i - \theta_j)} \\
    &+ (1 - \bar{y}_{ij})\log [{1 - \psi(\theta_i - \theta_j)}]\rbrace,
  \end{split}
  \label{eq:reg_loglikelihood}
\end{equation}
and $t \in \mathbb{R} \mapsto \psi(t) = {1}/{[1 + e^{-t}}]$ the sigmoid
function.

Under this notational setup, we are ready to state the $\ell_{\infty}$ upper
bound of the \btl{} MLE in \Cref{nthm:thm1}.

\bnthm \label{nthm:thm1} Assume the \btl{} model with parameter $\bbrtheta^* =
  (\theta^*_1,\ldots,\theta^*_n)^\top$ such that ${\bf 1}_n^\top \bbrtheta^* =0$
and a comparison graph $\mathcal{G} = \mclG([n],E)$ with adjacency matrix
${\bf A}$, algebraic connectivity $\lambda_2(\mclL_\bfA)$ and maximum and
minimum degrees $n_{\max}$ and $n_{\min}$. Suppose that each pair of items
$(i,j)\in E$ are compared $L$ times. Let $\kappa = \max_{i,j}|\theta^*_i -
  \theta^*_j|$ and $\kappa_E = \max_{(i,j)\in E}|\theta^*_i - \theta^*_j|$ and
set $\rho \geq c_{\rho}\kappa^{-2}e^{-2.5\kappa_E}n^{-4}n_{\max}^{1/2}$.
Assume that $\mclG$ is connected or $\lambda_2(\mclL_A) >0$. Then with
probability at least $1 - O(n^{-4})$, the regularized MLE
$\hat{\bbrtheta}_{\rho}$ from \eqref{eq:reg.mle} satisfies
\begin{align}
  \|\hat{\bbrtheta}_{\rho} - \bbrtheta^*\|_{\infty} \lesssim & \frac{e^{2\kappa_E}}{\lambda_2} \frac{n_{\max}}{n_{\min}}\parens{\sqrt{\frac{n+r}{L}} + \rho \kappa\sqrt{\frac{n}{n_{\max}}}} \nonumber \\
                                                             & + \frac{e^{\kappa_E}}{\lambda_2} \sqrt{\frac{n_{\max} (\log n + r)}{L}},
  \label{eq:thm1}                                                                                                                                                                                      \\
  \|\hat{\bbrtheta}_{\rho} - \bbrtheta^*\|_2 \lesssim        & \frac{e^{\kappa_E}}{\lambda_2} \parens{ \sqrt{\frac{n_{\max}(n+r)}{L}} + \rho\kappa \sqrt{n}}
  \label{eq:l2.rate}
\end{align}
where $\lambda_2 = \lambda_2(\mclL_\bfA)$, $r \defined \kappa_E + \log \kappa$
provided that $L\leq n^8e^{5\kappa_E}\max\{1,\kappa\}$, and $L$ is large enough
so that the right hand side of \Cref{eq:thm1} is smaller than a sufficiently
small constant $C>0$. In particular, if we set $\rho = c_\rho
  /{\kappa}\sqrt{{n_{\max}}/{L}}$ for some $c_\rho>0$, then
\begin{align}
  \|\hat{\bbrtheta}_{\rho} - \bbrtheta^*\|_{\infty} \lesssim & \frac{e^{2\kappa_E}}{\lambda_2} \frac{n_{\max}}{n_{\min}}\sqrt{\frac{n + r}{L}} + \frac{e^{\kappa_E}}{\lambda_2} \sqrt{\frac{n_{\max} (\log n + r)}{L}}, \nonumber \\
  \|\hat{\bbrtheta}_{\rho} - \bbrtheta^*\|_2 \lesssim        & \frac{e^{\kappa_E}}{\lambda_2} \sqrt{\frac{n_{\max}(n+r)}{L}}.
  \label{eq:thm1_optimal}
\end{align}
\enthm
As a brief sketch, the proof is based on a gradient descent procedure
initialized at $\bbrtheta^{(0)} = \bbrtheta^*$ and the idea is to control
$\|\bbrtheta^{(T)} - \hat{\bbrtheta}_{\rho}\|_{\infty}$ using the linear
convergence property and $\|\bbrtheta^{(T)} - {\bbrtheta^*}\|_{\infty}$ using
the leave-one-out technique in \cite{chen2019spectralregmletopk} and
\cite{chen2020partialtopkranking}. In fact, our work confirms that such a line
of argument extends to more general graph topologies beyond the Erd\"os-R\'enyi
graph, which is non-trivial. The proof details can be found in
\Cref{sec:prf_upper_bounds}.

\noindent {\bf Interpretation of key terms}

The upper bound in \Cref{eq:thm1} contains several distinct terms, which
interact with each other in non-trivial ways and express different aspects of
the intrinsic difficulty of the estimation task.
\begin{itemize}
  \item
        The factor  $\frac{e^{\kappa_E}}{\lambda_2(\mclL_{\bfA})}$ combines two
        sources of statistical hardness: the \textit{maximal gap} in performance
        $\kappa_E$ among the ranked items over the edge set $E$, and the
        \textit{algebraic connectivity} $\lambda_2(\mclL_{\bfA})$ of the
        comparison graph. It is intuitively clear that the larger the
        performance gap among the compared items, the more difficult it is to
        accurately estimate the model parameters. Furthermore, the smaller the
        algebraic connectivity, the less connected the comparison graph is, due
        to the presence of bottlenecks\footnote{Here, bottlenecks can be
          formally described as small connected subgraphs with very few edges
          separating dense portions of the graph.}. This in turn will increase the
        chance of obtaining a highly erroneous ranking or of gathering data from
        which a global ranking cannot be elicited at all. The minimal and
        maximal degrees $n_{\min}$ and $n_{\max}$ further quantify the impact of
        the connectivity of the comparison graph.
  \item We note that the factor $\frac{1}{\lambda_2(\mclL_{\bfA})}$ can be
        equivalently replaced with $\frac{1}{\lambda_2(\mclI)}$ (see
        \Cref{lm:lem8} in \Cref{sec:prf_upper_bounds}). Here, $\mclI \defined
          \nabla^2 \ell_0(\bbrtheta^*;\bfy)$ is the Fisher information matrix at
        $\bbrtheta^*$ and $\lambda_2(\mclI)$ its smallest  non-zero eigenvalue.
        The fact that the bound depends on the Fisher information is not too
        surprising. This is so, since this quantity in exponential families
        quantifies the curvature of the likelihood and the intrinsic difficulty
        of estimating $\bbrtheta^*$.
  \item Our bounds depend on both $\kappa$ and $\kappa_E$, which is non-standard
        in the literature. By definition, $\kappa_E \leq \kappa$ and in many
        cases, $\kappa_E$ can be much smaller than $\kappa$. We discuss this
        further in \Cref{sec:simulations}.
  \item The term $r \defined \kappa_E  +\log\kappa$ shows the impact of large
        $\kappa$ and $\kappa_E$. When $\kappa\lesssim n$ and $\kappa_E\lesssim
          \log n$, $r$ is negligible. We will consider this parameter range
        throughout the paper unless stated otherwise.
  \item The term $\sqrt{\frac{n}{L}}$ describes explicitly the impact of a
        high-dimensional parameter space on the estimation problem in relation
        to $L$, the number of samples for each comparison, which can be thought
        of as a measure of the sample size required for each of the $n$
        parameters. The inverse root dependence on $L$ is to be expected and, we
        conjecture, not improvable.
\end{itemize}

\bnrmk\label{nrmk:thm1-conditions} In the case of dense graphs, \eg, complete
graphs, $\lambda_2(\mclL_{\bfA})$ is large enough so that even $L = 1$ will
ensure a consistent estimator as $n\rightarrow \infty$. But for sparse graphs,
$L$ needs to be larger to compensate for weaker connectivity. The assumption
that $L\leq n^8e^{5\kappa_E}\max\{1,\kappa\}$ is a technical condition. There is
nothing special in the exponent for $n$. Any fixed number larger than $8$ can be
used which will only affect the constants in the bounds. The condition $L \leq
  n^8e^{5\kappa_E}\max\{1,\kappa\}$ may seem counter-intuitive, since it places an
upper bound on the sample size. But a control over $L$ is needed because as $L$
gets larger, the optimal choice of the regularization parameter $\rho = c_\rho
  \frac{1}{\kappa}\sqrt{\frac{n_{\max}}{L}}$ gets smaller and, accordingly, the
convergence rate of the gradient descent procedure upon which our proof is based
degrades. The optimal choice $\rho =
  c_{\rho}\frac{1}{\kappa}\sqrt{\frac{n_{\max}}{L}}$ depends on $\kappa$, which is
unknown before an estimator is produced, however, one can set $\rho =
  c_{\rho}\sqrt{\frac{n_{\max}}{L}}$ and the upper bound will only change by a
factor $\max\{1,\kappa\}$ in the first term of \Cref{eq:thm1_optimal}.
\enrmk

\subsection{Comparison to other work}\label{sec:comparison-to-other-work}

To the best of our knowledge, \cite{yan2012sparsecompbtl,
  hajek2014minimaxinferencepartialrank,shah2015estimationfrompairwisecomps,negahban2017rankcentralitypairwisecomparisons,
  agarwal2018acceleratedspectralranking,
  hendrickx2019graphresistance,hendrickx2020minimaxpairwisebtl} are the only
existing papers that study estimation error for the \btl{} model on a
comparison graph with general topology. Since
\cite{negahban2017rankcentralitypairwisecomparisons,agarwal2018acceleratedspectralranking,hendrickx2019graphresistance,hendrickx2020minimaxpairwisebtl}
estimate the the preference scores $\bfw^{*}$ rather than $\bbrtheta^*$, we
cannot directly compare our results with theirs because there is no tight
two-sided relationship between their metrics of error and ours. Therefore, here
we only compare our results to those in \cite{yan2012sparsecompbtl,
  hajek2014minimaxinferencepartialrank, shah2015estimationfrompairwisecomps}, as
is summarized in \Cref{tab:compare}.
We include the comparison to the other four papers in
\Cref{sec:comparison_detail}.

\textbf{$\ell_{\infty}$ loss: }\cite{yan2012sparsecompbtl} establish an
$\ell_{\infty}$-bound depending on $n_{ij}$, the number of common neighbors of
item $i$ and item $j$ in the comparison graph, under a strong assumption that
$n_{ij}\geq cn$ for some constant $c\in (0,1)$. This constraint on graph
topology is stronger than ours since it requires the graph to be dense. In
particular, when the comparison graph comes from an Erd\"os-R\'enyi model
$ER(n,p)$, $\min_{i,j}n_{ij}\asymp np^2$. Then the conditions in
\cite{yan2012sparsecompbtl} requires $p$ to be bounded away from $0$ and their
bound becomes $\frac{e^\kappa}{p}\sqrt{\frac{\log n}{npL}}$, while our bound is
$\frac{e^{2\kappa_E}}{\sqrt{p}}\sqrt{\frac{\log n}{npL}}$. Our bound is tighter
for moderate or small $\kappa_E$, and importantly, allows $p$ to vanish.
Furthermore, in \Cref{sec:simulations}, we show by some specific examples that
$\min_{i,j}n_{ij}$ could be 0 even for many fairly dense graphs, to illustrate
that the \cite{yan2012sparsecompbtl} upper bound cannot apply to many realistic
settings.

\textbf{$\ell_{2}$ loss: }\cite{hajek2014minimaxinferencepartialrank,
  shah2015estimationfrompairwisecomps} consider constrained MLE
$\hat{\bbrtheta}:=\min_{\|\bbrtheta\|_{\infty}\leq B}\ell_0(\bbrtheta)$ for a
known parameter $B$ such that $\| \bbrtheta^*\|_\infty \leq B$. Setting aside
the fact that their results require stricter conditions than ours, our
$\ell_{2}$ bound is tighter than theirs for general parameter settings with
moderate $B,\kappa$ and for a broad range of graphs with moderate
$\lambda_{2}(\mclL_{\bfA})$, i.e., not too sparse or irregular.

\begin{table}[htb!]
  \centering
  \begin{tabular}{c|c|c}
    \hline
    \textbf{Norm}                               & \textbf{Reference}          &
    \textbf{Upper bound}
    \\ \hline\hline
    $\|\cdot\|_\infty$                          & \cite{yan2012sparsecompbtl} &
    $\frac{e^{\kappa}}{\min_{i,j}n_{ij}}\sqrt{\frac{n_{\max}\log n}{L}}$
    \\
    \cline{2-3}                                 & \textbf{Our work}           &
    See \Cref{nthm:thm1}
    \\
    \hline
                                                &
    \cite{hajek2014minimaxinferencepartialrank} & $e^{8B}\frac{|E|\log
        n}{\lambda_2(\mathcal{L}_A)^2L}$
    \\ \cline{2-3} $\|\cdot\|^2_2$     &
    \cite{shah2015estimationfrompairwisecomps}  & $e^{8B}\frac{n\log
        n}{\lambda_2({\mclL}_A)L}$
    \\ \cline{2-3} &  \textbf{Our work}     &
    $\frac{e^{2\kappa_E}}{\lambda_2(\mclL_\bfA)^2} \frac{n_{\max}n}{L}$
    \\
    \hline
  \end{tabular}
  \caption{Comparison of results in literature.}
  \label{tab:compare}
\end{table}

We re-emphasize that \cite{hendrickx2020minimaxpairwisebtl} also provide upper
bounds for a general fixed comparison graph that matches an instance-wise lower
bound, for their parameter of interest $\mathbf{w}^* \defined
  (e^{\theta_1^*},\ldots,e^{\theta_n^*})^{\top}$, instead of $\bbrtheta^*$.
However, their error metric, \ie, $\sin(\hat{\bfw}, \bfw^*)$, is quite different
from other similar papers in the \btl{} literature, including our work. As such,
it is not clear how to compare to their results. Furthermore, as noted in
\Cref{sec:introduction}, from the perspective of ranking, an entry-wise metric
like $\|\cdot\|_{\infty}$ is more informative than vector-level metrics like
$\|\cdot\|_2$ and $\sin(\cdot, \cdot)$.

\subsection{Special cases of  graph topologies}\label{sec:special-case} We can
check some common types of comparison graph topologies and see in what order the
necessary sample complexity $N_{\textnormal{comp}} =|E|L$ needs to be to achieve
consistency, i.e., $\|\hat{\bbrtheta} - \bbrtheta^*\|_{\infty} = o(1)$. The
results are summarized in \Cref{tab:special_cases}. For path and star graphs, we
used the specialized bounds in \Cref{prop:path,prop:tree}.
\begin{table}[!h]
  \centering
  \begin{tabular}{c|c|c}
    \hline
    \textbf{Graph}                &
    \specialcell{$\mathbf{N_{\textnormal{comp}}}$
    \\
    \citep{yan2012sparsecompbtl}} &
    \specialcell{$\mathbf{N_{\textnormal{comp}}}$
    \\
      \textbf{(Our work)}}
    \\
    \hline \hline
    \textbf{Complete}             & $\Omega(n^2)$ & $\Omega(n^2)$
    \\
    \textbf{Bipartite}            & \texttt{N/A}  & $\Omega(n^2)$
    \\
    \textbf{Path}                 & \texttt{N/A}  & $\omega(e^{2\kappa_E}n^2\log
      n)$
    \\
    \textbf{Star}                 & \texttt{N/A}  & $\omega(e^{2\kappa_E}n \log
      n)$
    \\
    \textbf{Barbell}              & \texttt{N/A}  & $\omega(e^{2\kappa_E}n^5\log
      n)$
    \\
    \hline
  \end{tabular}
  \caption{Magnitude of $N_{\textnormal{comp}}$ to ensure $\|\hat{\bbrtheta} -
    \bbrtheta^*\|_{\infty} = o(1)$.}
  \label{tab:special_cases}
\end{table}
As shown in \Cref{tab:special_cases}, our bound now applies to a much broader
class of graph topologies under the $\ell_{\infty}$-norm compared to
\citet{yan2012sparsecompbtl}.

\bnrmk \label{rmk:sample-complexity} For the path graph, star graph, and barbell
graph, the necessary sample complexity induced by directly applying our
$\ell_{\infty}$ bound is larger than the sample complexity induced by the
$\ell_2$ bound in \cite{shah2015estimationfrompairwisecomps}, though they
require more stringent conditions than ours. Thus we provide specialized sharp
upper bounds in the case of path and star graph in \Cref{prop:path} and
\ref{prop:tree}. Additionally, in \Cref{sec:implications-of-work}, we illustrate
that by applying a \textit{unique} sub-additivity property of
$\ell_{\infty}$-loss, we can achieve a much smaller sample complexity in graphs
with bottlenecks like the barbell graph.
\enrmk

\textbf{Erd\"os-R\'enyi graph: }By applying a union bound on
$\lambda_2(\mclL_{\bfA})$, $n_{\max}$, and $n_{\min}$ to the sample-wise bounds
in \Cref{nthm:thm1}, we obtain a corollary in the setting where the comparison
graph follows the Erd\"os-R\'enyi model $ER(n,p)$.

\bncor[Erd\"os-R\'enyi graph]\label{cor:cor_ER} As a corollary to
\Cref{nthm:thm1}, suppose that the comparison graph comes from an
Erd\"os-R\'enyi graph $ER(n,p)$, then under the same conditions, with
probability at least $1 - O(n^{-4})$, it holds that
\begin{equation*}
  \|\hat{\bbrtheta}_{\rho} - \bbrtheta^*\|_{\infty} \lesssim e^{2\kappa_E} \sqrt{\frac{\log n}{np^2L}},
  \|\hat{\bbrtheta}_{\rho} - \bbrtheta^*\|_2 \lesssim  {e^{\kappa_E}} \sqrt{\frac{1}{pL}}.
\end{equation*}
\encor
The full form of \Cref{cor:cor_ER} with a proof can be found at the end of
\Cref{sec:prf_upper_bounds}. For the Erd\"os-R\'enyi comparison graph $ER(n,p)$,
the tightest $\ell_{\infty}$-norm error bound $e^{2\kappa}\sqrt{\frac{\log
      n}{npL}}$ is proved in \cite{chen2019spectralregmletopk} and
\cite{chen2020partialtopkranking}. \cite{han2020asymptoticsparsebradleyterry}
establish an $\ell_{\infty}$-norm upper bound of $e^{2\kappa}\sqrt{\frac{\log
      n}{np}}\cdot \frac{\log n}{\log (np)}$.
\cite{negahban2017rankcentralitypairwisecomparisons} obtain an $\ell_{2}$-norm
upper bound of $e^{4\kappa}\frac{\log n}{pL}$ and a lower bound of $e^{-\kappa}
  \frac{1}{pL}$. Thus the derived $\ell_{2}$-bound in \Cref{cor:cor_ER} in
Erd\"os-R\'enyi case is minimax optimal.

In this case our derived $\ell_{\infty}$-bound cannot achieve the rate
established in \cite{chen2019spectralregmletopk},
\cite{chen2020partialtopkranking}, though our $\ell_{2}$-bound exhibits the
optimal rate proved in \cite{negahban2017rankcentralitypairwisecomparisons}. The
reason why our bound does not imply the optimal  $\ell_{\infty}$-rate under a
Erd\"os-R\'enyi comparison graph is that our bound is a sample-wise bound and
thus cannot leverage some regular property of Erd\"os-R\'enyi graph beyond
algebraic connectivity and degree homogeneity that is exhibited with high
probability.

\textbf{Tree graphs:} For extremely sparse graphs like tree graphs, the general
upper bound in \Cref{nthm:thm1} is loose compared to the lower bound in
\Cref{nthm:thm2-lb}. Therefore, we separately prove some sharp upper bounds for
path and star graphs as a complement to our general theory, in these frequently
studied cases. For example, single-elimination sports tournaments are commonly
designed as a binary tree graph. By the spectral property of path and star
graphs (see \Cref{sec:appendix-special-case}), one can verify that the upper
bounds in both norms match the $\ell_{\infty}$ lower bound in
\Cref{nthm:thm2-lb} and the $\ell_{2}$ lower bound in
\cite{shah2015estimationfrompairwisecomps}, up to $\sqrt{\log n}$ and
$e^{2\kappa_E}$ factors.

\bnprop[Path graph]
\label{prop:path}
Suppose the comparison graph is a path graph $([n],E)$ with $E = \{(i,
  i+1)\}_{i\in [n - 1]}$ and $L>c e^{2\kappa_E}n\log n$ for some universal
constant $c$, then with probability at least $1 - n^{-4}$, the vanilla MLE
$\hat{\bbrtheta}_0$ satisfies
\begin{equation*}
  \begin{split}
    \|\hat{\bbrtheta}_0 - \bbrtheta^*\|_{\infty}&\lesssim e^{\kappa_E}\sqrt{\frac{n\log n}{L}}, \\
    \|\hat{\bbrtheta}_0 - \bbrtheta^*\|_{2} &\lesssim e^{\kappa_E}n\sqrt{\frac{\log n}{L}}.
  \end{split}
\end{equation*}
\enprop
\bnprop[General tree graph]
\label{prop:tree}
Suppose the graph is a tree graph $([n],E)$ where each item $i$ and $j$ are
compared $L$ times such that $L>c e^{2\kappa_E}n\log n$ for some universal
constant $c$. Then with probability at least $1 - n^{-4}$, the vanilla MLE
$\hat{\bbrtheta}_0$ satisfies
\begin{equation*}
  \begin{split}
    \|\hat{\bbrtheta}_0 - \bbrtheta^*\|_{\infty}&\lesssim e^{\kappa_E}\sqrt{\frac{D\log n}{L}},\\
    \|\hat{\bbrtheta}_0 - \bbrtheta^*\|_{2}&\lesssim e^{\kappa_E}\sqrt{\frac{Dn\log n}{L}},
  \end{split}
\end{equation*}
where $D:=\max_{i,j}|{\rm path}(i,j)|$ is the diameter. In particular, for star
graph, the upper bound is given by $D = 1$.
\enprop
The full form of \Cref{prop:path} and \Cref{prop:tree} with proofs are found in
\Cref{sec:prf_upper_bounds}. Briefly, the proofs leverage the closed-form
solution of vanilla MLE under the tree graph.

\section{Lower bounds}\label{sec:lower-bounds} In this section, we derive a
minimax lower bound for the $\ell_{\infty}$ loss. Towards that end, we first
introduce some new notation. Let $N_{\textnormal{comp}}$ be the total number of
comparisons that have been observed, so in our setting, $N_{\textnormal{comp}} =
  |E|L$ where $|E|$ is number of edges in the comparison graph $\mclG$. Denote the
two items involved in the $i$-th comparison as $(i_1,i_2)$ such that $i_1<i_2$.
Let $\tilde{\mclL}_A =
  \frac{1}{N_{\textnormal{comp}}}\sum_{i=1}^{N_{\textnormal{comp}}}(\bfe_{i_1} -
  \bfe_{i_2})(\bfe_{i_1} - \bfe_{i_2})^\top$ be the normalized graph Laplacian
with pseudo inverse $\tilde{\mclL}_A^{\dagger}$ and eigenvalues
$0=\lambda_1(\tilde{\mclL}_A) \leq \lambda_2(\tilde{\mclL}_A)\leq\cdots\leq
  \lambda_n(\tilde{\mclL}_A)$. With the main notation in place, our minimax lower
bound is summarized in the following result.

\bnthm\label{nthm:thm2-lb} Assume that the comparison graph $\mclG$ is connected
and the sample size $N_{\textnormal{comp}} \geq \frac{c_{2}
    \operatorname{tr}\left(\tilde{\mclL}_A^{\dagger}\right)}{e^{2\kappa}
    \kappa^{2}}$, any estimator $\widetilde{\bbrtheta}$ based on
$N_{\textnormal{comp}}$ comparisons with outcomes from the \btl{} model
satisfies
\begin{align*}
  \sup_{\bbrtheta^*\in \Theta_\kappa} \mathbb{E} & \left[\|\widetilde{\bbrtheta} - \bbrtheta^{*}\|^2_{\infty}\right]
  \gtrsim \frac{e^{-2\kappa}}{n N_{\textnormal{comp}} } \enspace \times                                              \\
                                                 & \max \Big\{{n}^2,
  \max_{n^{\prime} \in\{2, \ldots, n\}}  \sum_{i=\left\lceil 0.99 n^{\prime}\right\rceil}^{n^{\prime}} [\lambda_{i}(\tilde{\mclL}_A)]^{-1} \Big\}
\end{align*}
where $\Theta_{\kappa} = \{\theta\in \mathbb{R}^n:\mathbf{1}_n^\top \bbrtheta =
  0,\ \|\bbrtheta\|_{\infty}\leq \kappa\}$.
\enthm

The proof of \Cref{nthm:thm2-lb} largely leverages the lower bound construction
from Theorem 2 in \cite{shah2015estimationfrompairwisecomps}. The main
modification in adapting it to our setting is to construct an
$\ell_{\infty}$-packing set. This is done by utilizing the \textit{tight}
topological equivalence of $\ell_{\infty}$ and $\ell_{2}$ norms in finite
dimensions.

We can compare this lower bound with the upper bound in \Cref{nthm:thm1}. In our
setting, the comparisons distribute evenly over all pairs, so
$N_{\textnormal{comp}} = |E|L$, and $\lambda_i(\tilde{\mclL}_A) = \frac{1}{|E|}
  \lambda_i(\mathcal{L}_A)$. Thus, given a comparison graph with
$\lambda_2(\tilde{\mclL}_A) \asymp \frac{1}{n}$, the lower bound becomes
\[
  \sup _{\bbrtheta^*\in \Theta_{\kappa}} \mathbb{E}\left[\|\widetilde{\bbrtheta}-\bbrtheta^{*}\|_{\infty}\right]
  \gtrsim e^{-\kappa} \sqrt{ \frac{n}{N_{\textnormal{comp}}} }
\]
In $ER(n,p)$ case, this lower bound becomes $e^{-\kappa}\sqrt{\frac{1}{npL}}$
which matches the upper bound in \cite{chen2019spectralregmletopk}. For some
``regular'' graph topology with $\lambda_2(\tilde{\mclL}_A) \asymp \frac{1}{n}$
like complete graph, expander graph with $\phi=\Omega(n)$ and complete bipartite
graph with two partition sets of size $\Omega(n)$, the upper bound becomes
\begin{equation*}
  \|\hat{\bbrtheta}_{\rho} - \bbrtheta^*\|_{\infty} \lesssim e^{2\kappa} \sqrt{\frac{n\log n}{N_{\textnormal{comp}}}}.
\end{equation*}
Therefore, when the comparison graph topology is sufficiently regular, our upper
bound matches the lower bound up to a $\log n$ factor and a factor of
$e^{3\kappa}$. As a final remark,
\cite{negahban2017rankcentralitypairwisecomparisons} show that the minimax lower
bound for $\ell_{2}$-loss and Erd\"os-R\'enyi comparison graph $ER(n,p)$  is
$e^{-\kappa} \frac{1}{pL}$, which matches  our $\ell_{2}$ upper bound up to a
factor of $e^{2\kappa}$.

\section{Implications for tournament design}\label{sec:implications-of-work}

In this section, we discuss how our results can be leveraged to construct more
efficient tournament design from a ranking perspective in sports leagues.

As discussed in \Cref{sec:special-case}, for some comparison graphs with small
$\lambda_2(\mclL_{\bfA})$, the requirement on $L$ and $N_{\textnormal{comp}}$
for consistency is stringent. However, as we show next, we can significantly
relax the requirement on the sample complexity $N_{\textnormal{comp}}$ by
adaptively varying the number pairwise comparisons observed over different
subsets of the items in a manner that leverages different degrees of
connectivity of the comparison graphs.

The basic idea is that model parameters corresponding to a subset of items
inducing a highly connected sub-graph require relatively few observations. On
the other hand,  the outcomes of comparisons with items corresponding to nodes
of the comparison graph that are part of a ``graph bottleneck'' are especially
important in yielding accurate global ranking and, therefore, should be more
heavily sampled (in the sense of having a larger number $L$ of observations).
The case of a Barbell graph consisting of two complete sub-graphs connected by
few ``bridge'' edges (as is shown in \Cref{fig:barbell-bridge}) is an extreme
illustration of this situation and will be discussed below. In this case, it is
clear that the parameters corresponding the items adjacent to the bridge edges
ought to be estimated with higher accuracy and therefore, for those items $L$
should be set larger. Furthermore, it is possible to estimate the model
parameters separately over different sub-graphs and combine these estimators in
a way that could lead to an improved rate, compared to a joint or omnibus
estimator. Indeed, the next result shows that the $\ell_{\infty}$-error rate of
the combined estimator is bounded by the sum of the error rates for estimating
the parameters of the individual sub-graphs.

Formally, let $I_1,I_2,I_3$ be three subsets of $[n]$ such that $\cup_{j=1}^3
  I_j = [n]$ and, for each $j \neq k$, $I_j  \not\subseteq I_k$ and for $i=1,2$,
$I_i \cap I_3 \neq \emptyset$. Assume that the sub-graphs induced by $I_j$'s
are connected and the number of comparisons for all pairs can be different
across sub-graphs. Let $\bbrtheta^*$ be the vector of preference scores in the
\btl{} model over $n$ items and $\hat{\bbrtheta}_{(j)}$ be the MLE of
$\bbrtheta^*_{(j)}\in \mathbb{R}^{|I_j|}$ for the \btl{} model involving only
items in $I_j$, $j=1,2,3$. Also define the augmented version
$\tilde{\bbrtheta}_{(j)}\in \mathbb{R}^n$ such that
$\tilde{\bbrtheta}_{(j)}(I_j) = \hat{\bbrtheta}_{(j)}$.

Now take two nodes $t_1\in I_1\cap I_3$, $t_2\in I_2\cap I_3$, and let $\delta_3
  = \tilde{\bbrtheta}_{(1)}(t_1) - \tilde{\bbrtheta}_{(3)}(t_1)$, $\delta_2 =
  \tilde{\bbrtheta}_{(3)}(t_2) - \tilde{\bbrtheta}_{(2)}(t_2)$. An ensemble
estimator \textit{add-MLE} $\hat{\bbrtheta}\in \mathbb{R}^n$ is a vector such
that $\hat{\bbrtheta}(I_1) = \hat{\bbrtheta}_{(1)}$, $\hat{\bbrtheta}(S_2) =
  \tilde{\bbrtheta}_{(2)}(S_2) + \delta_3 + \delta_2$, and $\hat{\bbrtheta}(S_3) =
  \tilde{\bbrtheta}_{(3)}(S_3) + \delta_3$, where $S_2 = I_2\setminus I_1$ and
$S_3 = I_3\setminus (I_1\cup I_2)$. Notice that the value of $\hat{\bbrtheta}$
depends on the choice of $t_1,t_2$, but the estimation error of all ensemble
estimators can be well-bounded, as is shown in
\Cref{nlem:subbadditivity-ellinf-norm}. \bnprop[Subadditivity of
  $\ell_{\infty}$-loss in \btl{}]\label{nlem:subbadditivity-ellinf-norm} Under
the setting above, for any add-MLE $\hat{\bbrtheta}\in\mathbb{R}^n$ based on
$\hat{\bbrtheta}_{(1)}, \hat{\bbrtheta}_{(2)}, \hat{\bbrtheta}_{(3)}$, it
holds that
\begin{equation}
  d_{\infty}(\hat{\bbrtheta},\bbrtheta^*) \leq 4\sum_{i=1}^3 d_{\infty}(\hat{\bbrtheta}_{(i)},\bbrtheta^*_{(i)}),
\end{equation}
where $d_{\infty}(\bfv_1,\bfv_2) \defined \|(\bfv_1 - {\rm
  avg}(\bfv_1)\mathbf{1})-(\bfv_2 - {\rm avg}(\bfv_2)\mathbf{1})\|_{\infty}$ and
${\rm avg}(\bfx):=\frac{1}{n}\mathbf{1}_n^\top \bfx$ for $\bfx\in
  \mathbb{R}^n$.
\enprop
The proof of \Cref{nlem:subbadditivity-ellinf-norm} is found in
\Cref{sec:others}. For some types of graph topologies the above result can be
used to devise a {\it divide-and-conquer strategy} for estimating the model
parameters with better sample complexity than that of an omnibus estimator,
i.e., the joint-MLE in our setting. Indeed, as discussed in
\Cref{sec:special-case}, for a barbell graph containing two size-$n/2$ complete
sub-graphs connected by a single edge, we need $N_{\textnormal{comp}} =
  \Omega(n^5\log n)$ for an $o(1)$ error bound of the joint-MLE. From a practical
perspective, we note that such a divide and conquer strategy gives flexibility
in the number of comparisons in each sub-graph. For example, if we set $L=1$ for
the two complete sub-graphs to get MLEs $\hat{\bbrtheta}_1$,
$\hat{\bbrtheta}_2$, and set $L = n$ for the two items linking the two
sub-graphs to get an MLE $\hat{\bbrtheta}_3$, and combine them by shifting
$\hat{\bbrtheta}_2$ by the difference of two entries of $\hat{\bbrtheta}_3$,
then a total sample complexity \textit{reduction} to $N_{\rm comp}=\Omega(n^2)$
will ensure $\ell_{\infty}$-norm error of order $O(e^{2\kappa_E}\sqrt{\log
    n}/\sqrt{n})=o(1)$, because for a complete graph of size $m$, the
$\ell_{\infty}$-norm error is $O(e^{2\kappa_E}\sqrt{\log m}/\sqrt{mL})$. In
\Cref{exa:n_ij=0} and \Cref{subsec:additivity}, we show some simulation results
illustrating the advantage of using subadditivity in estimation, where we
generalize the add-MLE to Island graph and Barbell graph with multiple bridge
edges that can have more than $3$ dense sub-graphs.

Note that such flexible tournament design is similar to the idea of
\textit{active ranking} \citep{heckel2019activerankingpairwisecomps,
  ren2019activeranking}, but there is still a substantial difference between our
setting and active ranking. Active ranking assumes that one can design the
tournament in an \textit{online} manner, so that the next pair of items to be
compared is determined by the newest outcomes of comparisons. However, in
practice many tournaments can only be designed \textit{offline}, \ie, before any
games are played. Under this common setting, our $\ell_{\infty}$-subadditivity
property provides a useful offline approach to efficient tournament design.

\section{Examples and simulations}\label{sec:simulations} In this section, we
conduct numerical experiments on simulated data with two main goals. First, we
illustrate the utility of the subadditivity property in
\Cref{nlem:subbadditivity-ellinf-norm} in the case of Island graphs (see
\Cref{exa:n_ij=0}). Second, we demonstrate the relative tightness of our
$\ell_\infty$ upper bound compared to \cite{yan2012sparsecompbtl}, since their
work is closest in spirit to ours. Specifically, we compare the two bounds in a
setting where analytical comparison is not directly feasible (see
\Cref{exa:bridge}). All of our reproducible code is openly
accessible\footnote{Repo: \url{https://github.com/MountLee/btl_mle_l_inf}}.

In the \btl{} model, the maximal winning probability is $p_{\max}(\kappa) =
  1/({1 + e^{-\kappa}})$. To get a sense, $p_{\max}(2.20) = 0.900$,
$p_{\max}(4.59) = 0.990$. A winning probability larger than 0.99 is fairly
rare in practice, so it would not be too constraining to set $\kappa = 2.2$ in
our simulation. But analytically our result allows $\kappa$ to diverge with
$n$.

In our experiments, we set $\theta^*_i = \theta_1^* + (i-1)\delta$ for $i>1$
with $\delta = \kappa / (n-1)$. We additionally assign $\theta_1^*$ to ensure
that $\mathbf{1}_n^\top \bbrtheta^* = 0$, for parameter identifiability. Under
this setting, for some special graphs, e.g., the Island graph in
\cref{exa:n_ij=0}, $\kappa_E$ can be much smaller than $\kappa$, showing an
advantage of our upper bound in representing the dependency on the maximal
performance gap $\kappa_E$ along the edge set, rather than $\kappa$ the whole
vertex set. However, there may be some cases where the majority of edges have
small performance gaps and only a few edges have large gaps. Here, the control
in the upper bound purely by $\kappa_E$ can again be loose. An interesting
future direction is to make upper bounds tighter in such cases by including more
structural parameters, like the proportion of small-gap edges. We include some
illustrative examples
in \Cref{sec:additional-experiments}.

\begin{figure}[!h]
  \centering
  \includegraphics[width=0.4\textwidth]{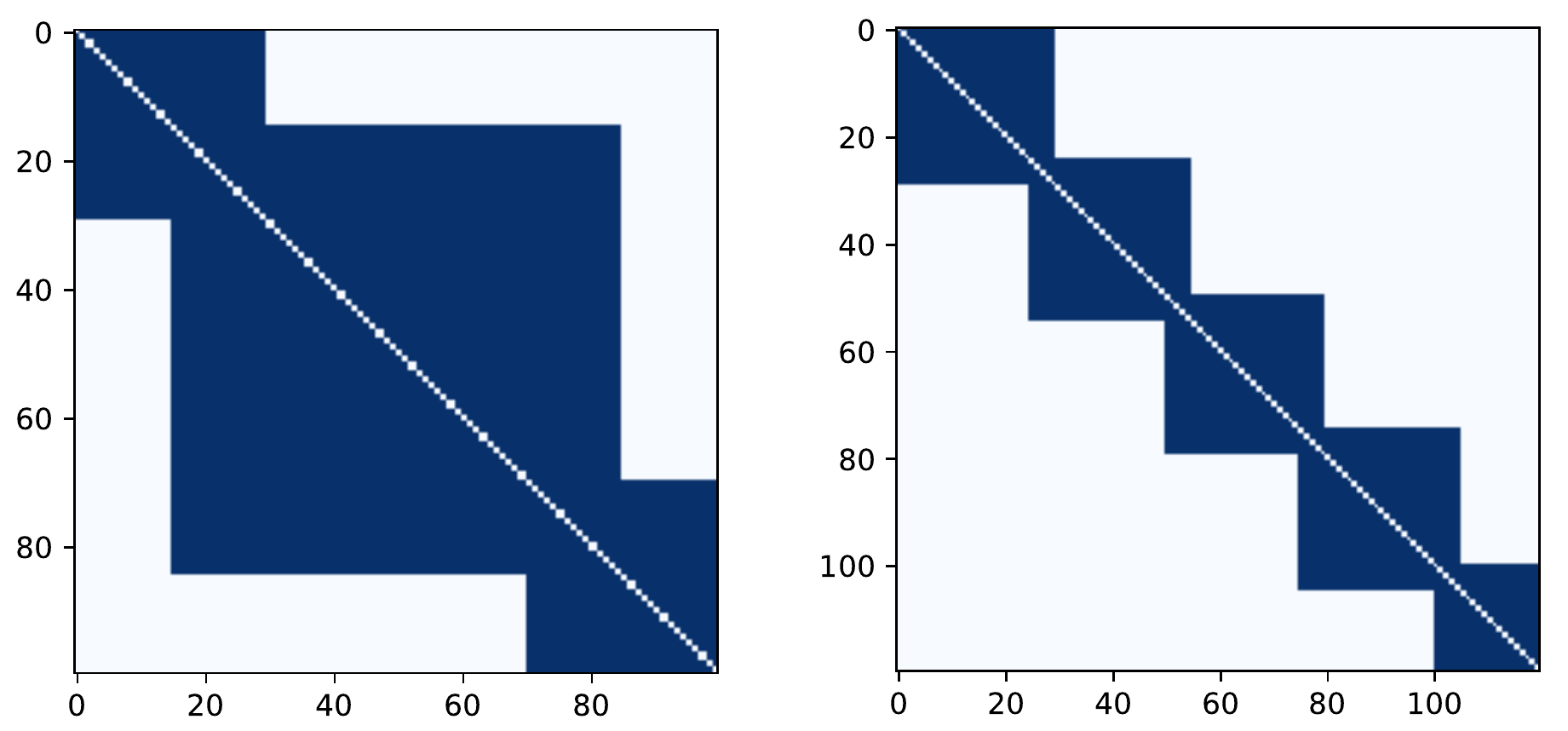}
  \includegraphics[width=0.4\textwidth]{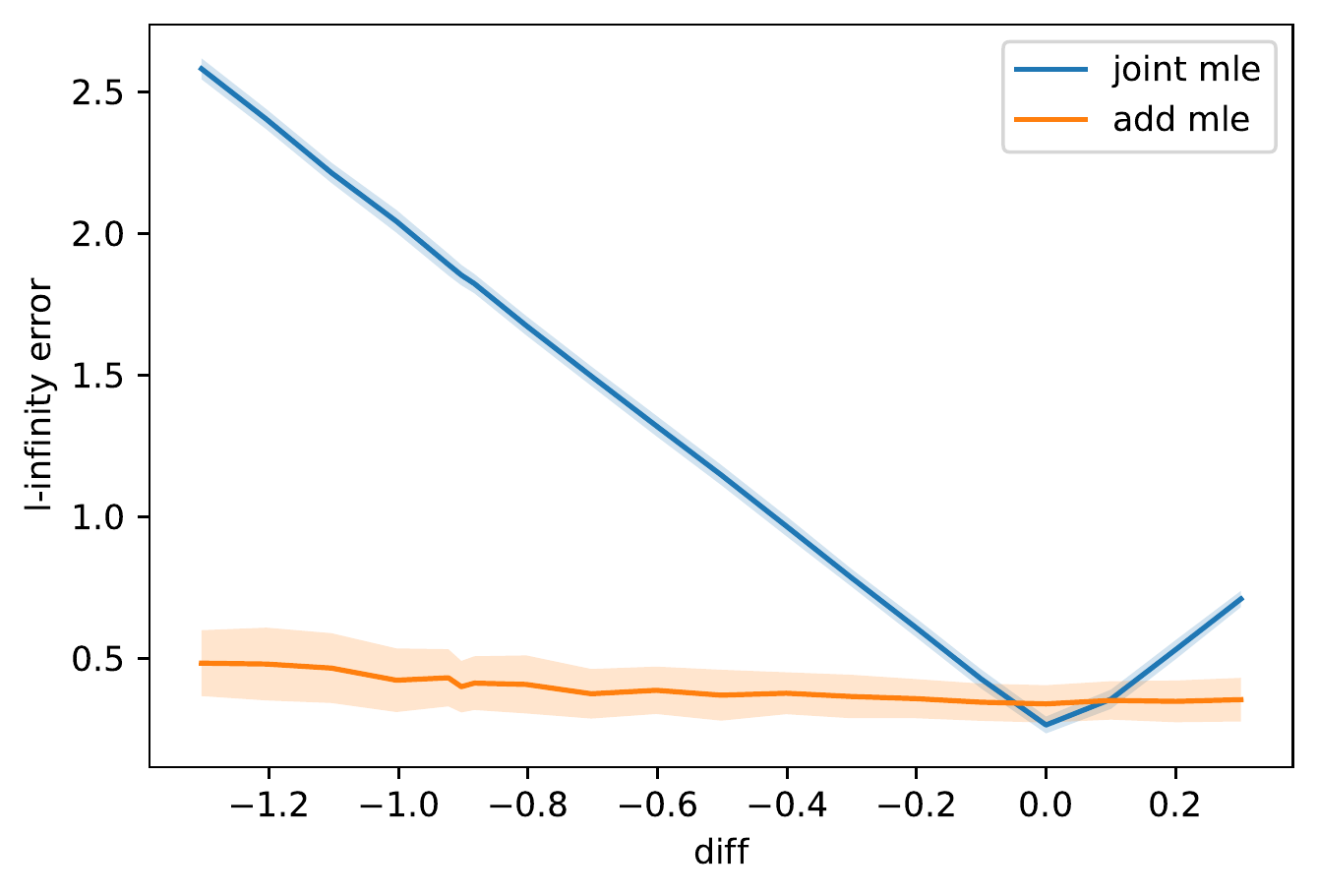}
  \caption{Left: Adjacency matrix of a 3-Island graph, with yellow indicating 1
    and purple indicating 0; $\lambda_2(\mclL_{\bfA}) = 11.92$. Right: Adjacency
    matrix of a general Island graph, with $n_{\text{island}} = 30$,
    $n_{\text{overlap}} = 5$, $n = 120$; $\lambda_2(\mclL_{\bfA}) = 1.19$.
    Bottom: comparison of the error of the joint-MLE and the add-MLE. The curve
    is obtained as the average of 100 trials with one standard deviation shown
    by the colored area.}
  \label{fig:island}
\end{figure}

\bnexa[Graph with $\min_{i,j}n_{ij}=0$] \label{exa:n_ij=0} In this case, we
intend to illustrate  that $\min_{i,j}n_{ij}$ could be $0$ or quite close to $0$
for even fairly dense graphs, making the upper bound in
\cite{yan2012sparsecompbtl} less effective. Consider a \textit{3-Island}
comparison graph $\mclG$ with $n$ nodes. The induced sub-graphs on node sets
$V_1$, $V_2$, $V_3$ with $|V_i| = n_i$ are complete graphs, where $V_1\cap V_3 =
  \emptyset$, $V_1\cup V_2\cup V_3 = [n]$, and $V_i\cap V_2\neq \emptyset$ for $i
  = 1,3$. There is no edge except for those within $V_1,V_2,V_3$. This graph
$\mclG$ is connected, and can be fairly dense if we make $n_2$ large, but
$\min_{i,j}n_{ij}=0$ always holds since $V_1\cap V_3 = \emptyset$ and the two
induced sub-graphs are complete. See \Cref{fig:island} left panel for a
visualization of the adjacency matrix of such a graph.

We can also consider more general \textit{Island} graphs. A general Island graph
is determined by $n$, the size of the graph, $n_{\text{island}}$, the size of
island sub-graphs, and $n_{\text{overlap}}$, the number of overlapped nodes
between islands. Each island sub-graph is a complete graph, and there is no edge
outside islands. For Island graphs, it holds that $\min_{i,j}n_{ij}=0$ and
$\kappa_E\approx \kappa \cdot n_{\text{island}}/n$. \Cref{fig:island} top panel
shows the adjacency matrix of two Island graphs. \Cref{fig:island} bottom panel
shows the comparison of the $\ell_{\infty}$-error of the joint-MLE and the
add-MLE (see the detailed definition in \Cref{sec:additional-experiments}) while
varying the difference in the average of preference scores of each island
sub-graph, where we set $n_{\text{island}} = 50, n_{\text{overlap}} = 5,L = 10$.
Every point on the lines is the average of 100 trials. It can be seen that the
add-MLE by the divide-and-conquer strategy largely dominates the joint-MLE in
$\ell_{\infty}$-error.
\enexa
In \Cref{exa:n_ij=0}, we show a common family of graphs which is fairly dense
while $\min_{i,j}n_{ij}=0$, so that the upper bound in
\cite{yan2012sparsecompbtl} does not hold. Next in \Cref{exa:bridge} we consider
another family of graphs where their upper bound holds but still relatively
looser than our bound.

\bnexa[Barbell graph with random \textit{bridge} edges] \label{exa:bridge}
Consider a generalized Barbell graph $\mclG$ containing $n = n_1 + n_2$ nodes,
where the induced sub-graph on nodes $\{1,\cdots,n_1\}$ and $\{n_1+1,\cdots,n\}$
are complete graphs, and the two sub-graphs are connected by some bridge edges
$(i,j)$ for some $1\leq i\leq n_1$  and $n_1 + 1\leq j\leq n$. Denote the set of
bridge edges as $E_l$, then $|E_l|/(n_1 n_2)$ quantifies the connectivity of
$\mclG$: the larger $|E_l|/(n_1 n_2)$ is, the denser or more regular $\mclG$ is.

\begin{figure}[htb!]
  \centering
  \includegraphics[width=0.46\textwidth]{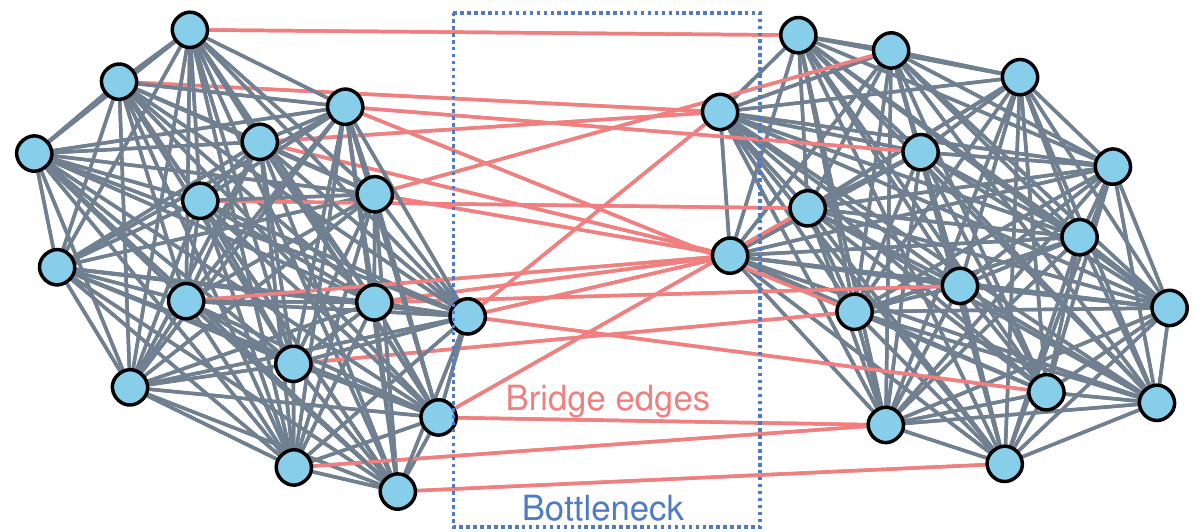}
  \includegraphics[width=0.44\textwidth]{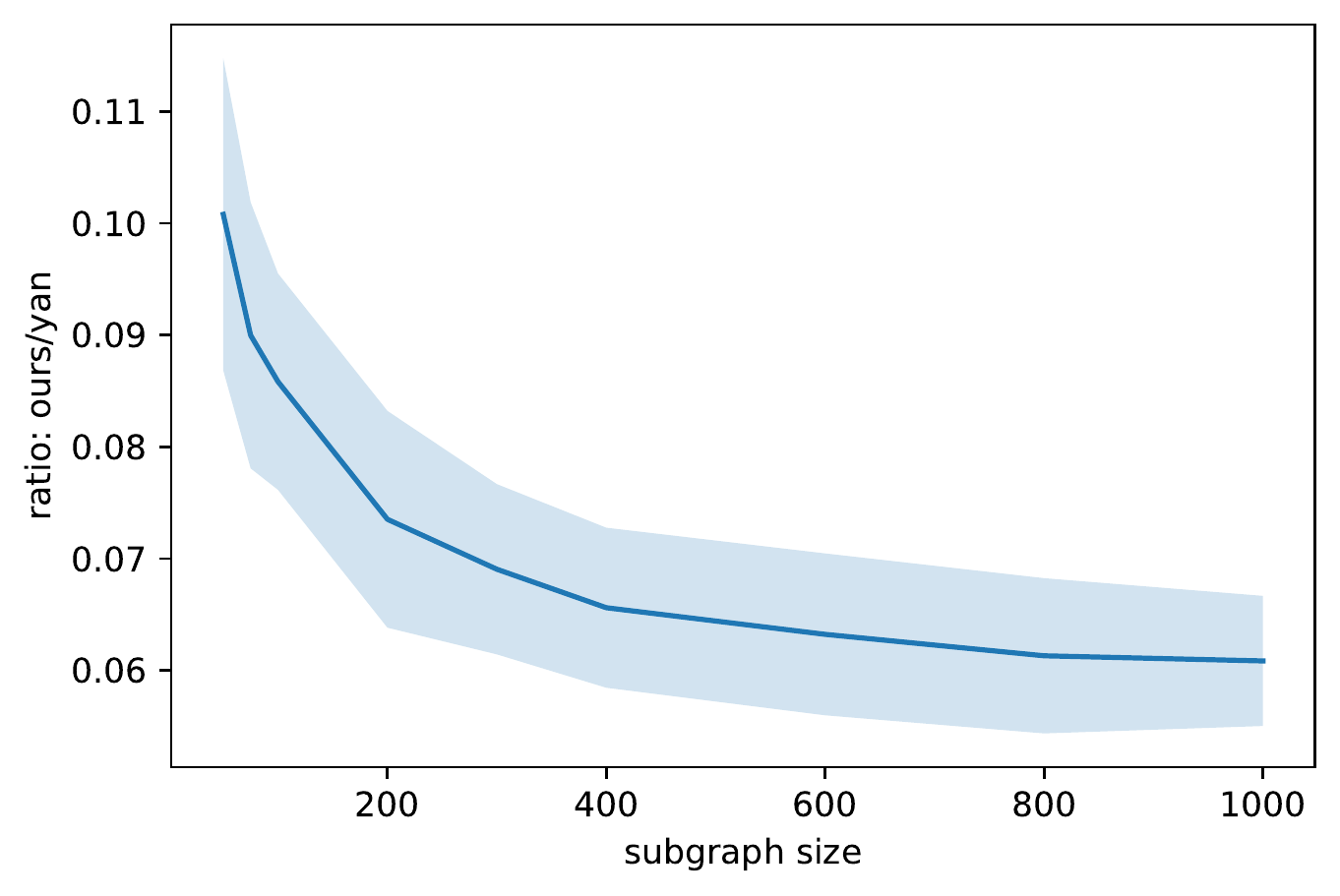}
  \caption{Top: visualization of a Barbell graph with random bridge edges.
    Bottom: The ratio of our bound and the bound in \cite{yan2012sparsecompbtl}
    under the Barbell graph with random bridge edges and sub-graph size $n_1 =
      n_2=n_s$ varying. The curve is obtained as the average of 100 trials with
    one standard deviation shown by the colored area.}
  \label{fig:barbell-bridge}
\end{figure}

In \Cref{fig:barbell-bridge} we show a comparison of the real
$\ell_{\infty}$-loss $\|\hat{\bbrtheta} - \bbrtheta^*\|_\infty$, and the upper
bounds of $\ell_{\infty}$-error in \cite{yan2012sparsecompbtl} and our paper. We
include this relative comparison to numerically demonstrate that our bound is in
general tighter than \citet{yan2012sparsecompbtl}, since there is no
\textit{known} analytical relationship between $\min_{i,j}n_{ij}$ and
$\lambda_2(\mathcal{L}_A)$ for general graphs. In our experiment, we set $n_1 =
  n_2 = n_s$, $L =10$, and randomly link $|E_l|=n_1 n_2 p$ edges between the two
complete sub-graphs, and vary $n_s$ from $50$ to $1000$ with $p =
  3\log(n_s)/n_s$. Every point on the line is the average of 100 trials. It can be
seen that our upper bound has a faster vanishing rate, compared to
\cite{yan2012sparsecompbtl} for this simulated scenario. This is evident as the
plotted ratio of our upper bound relative to the upper bound in
\citet{yan2012sparsecompbtl} has a steady decreasing trend, as $n$ increases. It
should be noted that there are leading constant factors in both upper bounds,
and for convenience we set them to be 1 for both bounds. Thus, one should focus
on the trend of the curve rather than the magnitude of the ratio in
\Cref{fig:barbell-bridge}.
\enexa

\section{Discussion}\label{sec:discussion}

In this work we provide a sharp risk analysis of the MLE for the \btl{} global
ranking model, under a more general graph topology, in the $\ell_{\infty}$-loss.
This addresses a major gap in the \btl{} literature, in extending the comparison
graph to more general and thus more practical settings, compared to dense graph
setting in \citet{yan2012sparsecompbtl}. Specifically we derive a novel upper
bound for the  $\ell_{\infty}$ and $\ell_2$-loss of the \btl{}-MLE, showing
explicit dependence on the algebraic connectivity of the graph, the sample
complexity, and the maximal performance gap between compared items. We also
derive lower bounds for the $\ell_{\infty}$-loss and analyze specific topologies
under which the MLE is nearly minimax optimal. We also show that the
$\ell_{\infty}$-loss satisfies a unique subadditivity property for the \btl{}
MLE and utilize our derived bounds for efficient tournament design. We note that
our upper bound is suboptimal in the cases where the graph topology is extremely
sparse or irregular. Although we provide sharp upper bounds for path and star
graphs as separate propositions, we still miss optimality other graph
topologies. A good future direction would be to optimize the upper and lower
bounds in such comparison graph regimes. Another promising direction is to
extend this analysis to the multi-user ranking models as in
\citet{jin2020rankaggviahetthurstmod}.

\noindent{\bf Acknowledgments}\label{subsec:acknowledgments}

We would like to thank Heejong Bong from the Carnegie Mellon University (CMU)
Department of Statistics \& Data Science, for his valuable feedback and
discussions during this work. We would also like to thank the anonymous
reviewers for their feedback which greatly helped improve our exposition.

\bibliography{refs}

\iftoggle{arxiv2022version}{
  \clearpage
  \onecolumn
  \appendix
  \section{Appendix}
  The Appendix contains following parts:
\begin{itemize}
  \item[\ref{sec:comparison_detail}] Detail on comparison of results.
  \item[\ref{sec:prf_upper_bounds}] Proof in \Cref{sec:upper-bounds}, i.e., upper bounds in
    \Cref{nthm:thm1,cor:cor_ER,prop:path,prop:tree}.
  \item[\ref{subsec:prf_lower_bounds}] Proof in \Cref{sec:lower-bounds}, i.e., lower bounds in
    \Cref{nthm:thm2-lb}.
  \item[\ref{sec:additional-experiments}] Additional experiments.
  \item[\ref{sec:others}] Other supporting results, including the proof for
    \Cref{nlem:subbadditivity-ellinf-norm}.
  \item[\ref{sec:appendix-special-case}] Special cases of comparison graphs,
    details related to \Cref{sec:comparison-to-other-work}.
  \item[\ref{sec:vanilla-mle}] Upper bound for unregularized/vanilla MLE.
\end{itemize}

\subsection{Comparison of results}
\label{sec:comparison_detail}

This section is a complement to \Cref{sec:comparison-to-other-work}. We will
summarize all existing works on the estimation error of Bradley-Terry model in
two tables, and then compare our results with the results in
\cite{negahban2017rankcentralitypairwisecomparisons,agarwal2018acceleratedspectralranking,hendrickx2020minimaxpairwisebtl}
in detail.

For simplicity, in \Cref{tab:ER} and \Cref{tab:general} we use $\kappa$ to
replace $B$ for results in \cite{hajek2014minimaxinferencepartialrank,
shah2015estimationfrompairwisecomps} as $\kappa \asymp B$ when $\bm{1}^\top
\bbrtheta^* = 0$, and in \Cref{tab:general} we omit the lower bound as they are
usually in fairly complex forms.

\begin{table}[htb!]
\centering
\begin{tabular}{|c|c|c|}
\hline
\textbf{Norm} & \textbf{Reference} &  \textbf{Upper Bound}  \\ \hline
     &   \cite{simons1999}    &      $p=1,\quad \lesssim
     e^{\kappa}\sqrt{\frac{\log n}{nL}}$           \\ \cline{2-3} &
     \cite{yan2012sparsecompbtl}  &       $\lesssim
     e^{2\kappa}\frac{1}{p}\sqrt{\frac{\log n}{npL}}$       \\ \cline{2-3}
     $\|\cdot\|_\infty$     &    \cite{han2020asymptoticsparsebradleyterry}   &
     $\lesssim e^{2\kappa}\sqrt{\frac{\log n}{np}}\cdot \frac{\log n}{\log
     (np)}$               \\ \cline{2-3} &   \cite{chen2019spectralregmletopk},
     \cite{chen2020partialtopkranking}    &    $\lesssim
     e^{2\kappa}\sqrt{\frac{\log n}{npL}}$               \\ \cline{2-3} &
     \textbf{Our work}    &       $\lesssim e^{2\kappa_E}\sqrt{\frac{\log
     n}{np^2L}}$            \\ \hline
     &  \cite{hajek2014minimaxinferencepartialrank}     &      $\lesssim
e^{8\kappa}\frac{\log n}{pL}$            \\ \cline{2-3} $\|\cdot\|^2_2$     &
\cite{shah2015estimationfrompairwisecomps}     &      $\lesssim
e^{8\kappa}\frac{\log n}{pL}$,\quad $\gtrsim e^{-2\kappa}\frac{1}{pL}$
\\ \cline{2-3} & \cite{negahban2017rankcentralitypairwisecomparisons}  &
$\lesssim e^{4\kappa}\frac{\log n}{pL}$, \quad    $\gtrsim e^{-\kappa}
\frac{1}{pL}$            \\ \cline{2-3} &  \cite{chen2019spectralregmletopk},
\cite{chen2020partialtopkranking}     &    $\lesssim e^{2\kappa}\frac{1}{pL}$
\\ \hline
$\sin^2(\cdot,\cdot)$ ($\hat{\bfw}$)     &
\cite{hendrickx2020minimaxpairwisebtl}      &    $\lesssim
e^{2\kappa}\frac{1}{pL\|w\|_2^2}$                 \\ \hline
$\|\cdot\|_1$ ($\hat{\bfw}$)     & \cite{agarwal2018acceleratedspectralranking}
&    $\lesssim e^{\kappa}\sqrt{\frac{\log n}{L}}$                 \\ \hline
\end{tabular}
\caption{Comparison of results under $ER(n,p)$ in literature.}
\label{tab:ER}
\end{table}

\begin{table}[htb!]
\centering
\begin{tabular}{|c|c|c|}
\hline
\textbf{Norm} & \textbf{Reference} &  \textbf{Upper Bound}  \\ \hline
$\|\cdot\|_\infty$    &  \cite{yan2012sparsecompbtl}  &
  $\frac{e^{\kappa}}{\min_{i,j}n_{ij}}\sqrt{\frac{n_{\max}\log n}{L}}$
  \\ \cline{2-3} &   \textbf{Our work}    &
  $\frac{e^{2\kappa_E}}{\lambda_2(\mclL_\bfA)}
  \frac{n_{\max}}{n_{\min}}\sqrt{\frac{n}{L}} +
  \frac{e^{\kappa_E}}{\lambda_2(\mclL_\bfA)} \sqrt{\frac{n_{\max} \log n}{L}}$
  \\ \hline
  &  \cite{hajek2014minimaxinferencepartialrank}     &      $\lesssim
e^{8\kappa}\frac{|E|\log n}{\lambda_2(\mathcal{L}_A)^2L}$            \\
\cline{2-3} $\|\cdot\|^2_2$     &  \cite{shah2015estimationfrompairwisecomps}
&      $\lesssim e^{8\kappa}\frac{n\log n}{\lambda_2({\mclL}_A)L}$  \\
\cline{2-3} &  \textbf{Our work}     &
$\frac{e^{2\kappa_E}}{\lambda_2(\mclL_\bfA)^2} \frac{n_{\max}n}{L}$  \\ \hline
$\sin^2(\cdot,\cdot)$ ($\hat{\bfw}$)     &
\cite{hendrickx2020minimaxpairwisebtl}      &    $\lesssim
e^{2\kappa}\frac{Tr(\mathcal{L}_{A}^{\dagger})}{L\|w\|_2^2}$                 \\
\hline
$\|\cdot\|_1$ ($\hat{\bfw}$)    & \cite{agarwal2018acceleratedspectralranking}
&    $\lesssim \frac{\eta e^{\kappa}
n_{avg}}{\lambda_2(D^{-1}A)n_{\min}}\sqrt{\frac{\log n}{L}}$                \\
\hline
\end{tabular}
\caption{Comparison of results for a fixed general comparison graph in literature.}
\label{tab:general}
\end{table}

In \cite{negahban2017rankcentralitypairwisecomparisons}, they establish an
$\ell_2$ upper bound for $\|\hat{\bbrpi} -
\tilde{\bbrpi}\|_2/\|\tilde{\bbrpi}\|_2$ in the order of
$\frac{e^{2.5\kappa}}{\lambda_2(\mathcal{L}_{rw})}\sqrt{\frac{n_{\max}\log
n}{L}}$ where $\tilde{\pi}(i) := w_i/\sum_{j}w_j$ with $w_i = \exp(\theta_i)$,
$\hat{\bbrpi}$ is the rank centrality estimamtor of $\tilde{\bbrpi}$,
$\lambda_2$ refers to the second smallest eigenvalue, $\mathcal{L}_{rw} =
D^{-1}A$ (which has the same spectrum as $D^{-1/2}\mathcal{L}_AD^{-1/2}$), and
$\mathcal{L}_A = D-A$. Recall that our $\ell_2$ upper bound is for
$\|\hat{\theta} - \theta\|_2$ and the order is
$\frac{e^{\kappa_E}}{\lambda_2(\mathcal{L}_A)}\sqrt{\frac{n_{\max}n}{L}}$. We
can now see that it's hard to give a general comparison between the two results
because 1. for a general graph, there is no precise relationship between
$\lambda_2(\mathcal{L}_{rw})$ and $\lambda_2(\mathcal{L}_A)$; 2. more
importantly, for a general model parameter $\theta$, there is no tight two-sided
bound between $\|\hat{\theta} - \theta\|_2$ and $\|\pi -
\tilde{\pi}\|_2/\|\tilde{\pi}\|_2$. Although, it would be a very interesting
future work to give a tight description of these two relevant pairs of
quantities and make a meaningful comparison.

\cite{agarwal2018acceleratedspectralranking} establish an $\ell_1$-norm upper
bound for the score parameter $\tilde{\pi}_i:= w_i/\sum_{j}w_j$ with $w_i =
\exp(\theta_i)$. Their bound is of the order $\frac{\eta e^{\kappa}
n_{avg}}{\lambda_2(D^{-1}A)n_{\min}}\sqrt{\frac{\log n}{L}}$, where $n_{avg} =
\sum_{i\in [n]}\tilde{\pi}_i n_i$, $D = {\rm diag}(n_1,\cdots,n_n)$, and $\eta
:= \log\left( \frac{n_{avg}}{n_{\min} \pi_{\min}} \right)$ with $\pi_{\min} =
\min_{i\in [n]}\tilde{\pi}_i$.

In \cite{hendrickx2020minimaxpairwisebtl}, they propose a novel weighted least
square method to estimate vector $w$, with $w_i = \exp(\theta_i)$, and provide
delicate theoretical analysis of their method. Their estimator shows a sharp
upper bound for $\mathbb{E}[\sin^2(\hat{w},w)]$ and equivalently for
$\mathbb{E}\|\hat{w}/\|\hat{w}\|_2 - w/\|w\|_2\|^2_2$, in the sense that the
upper bound for $\mathbb{E}[\sin^2(\hat{w},w)]$ matches a instance-wise lower
bound up to constant factors. Such a universal sharp/optimal bound for general
comparison graph (although the lower bound is not in the form of minimax rate)
is unique in literature. For convenience of comparison, here we assume $w_i$'s
are $O(1)$ (otherwise we can put a factor $e^{2B}$ in the bound
$\frac{Tr(\mathcal{L}_A^{\dagger})}{L}$). Then their upper bound is of the order
$\frac{Tr(\mathcal{L}_A^\dagger)}{L \|w\|_2^2}$, where $\mathcal{L}_A^\dagger$
refers to the Moore-Penrose pseudo inverse of the graph Laplacian of the
comparison graph. To correct for their different choice of metric, we need to
multiply $\|w\|^2_2$ to their bound, and it becomes
$\frac{Tr(\mathcal{L}_A^\dagger)}{L}$. On the other hand, the upper bound for
expected $\ell_2$ loss in \cite{shah2015estimationfrompairwisecomps} is
$\frac{n}{L\lambda_2(\mathcal{L}_A)}$. Since $\lambda_2(\mathcal{L}_A)$ is the
smallest positive eigenvalue of $\mathcal{L}_A$, it holds that
$Tr(\mathcal{L}_A^\dagger)< n/\lambda_2(\mathcal{L}_A)$, and hence the upper
bound for expected error in \cite{hendrickx2020minimaxpairwisebtl} is tighter
than the one in \cite{shah2015estimationfrompairwisecomps}.
Although, it should be noted that the loss function in
\cite{hendrickx2020minimaxpairwisebtl} is not directly comparable to a plain
$\ell_2$ loss, and we are just doing an approximate comparison.

In our paper, however, we provide a high probability bound for $\|\hat{\theta} -
\theta\|^2_2$ in the order of $\frac{n_{\max}n}{L\lambda^2_2(\mathcal{L}_A)}$.
It's usually hard to make a fair comparison between a high probability bound and
a bound for expected metrics, but in
\cite{hajek2014minimaxinferencepartialrank,shah2015estimationfrompairwisecomps},
they also provide a high probability bound in equation (8b) of Theorem 2. As we
discussed in Section 2.2 and Section 5, the $\ell_2$ bound is not the primary
focus of our paper, and although our theoretical analysis is not optimized for
$\ell_2$ error, our bound is still tighter than the bound in
\cite{hajek2014minimaxinferencepartialrank,shah2015estimationfrompairwisecomps}
for moderately dense and regular graphs, and is only worse for fairly sparse and
irregular graphs.


\clearpage


\subsection{Proof in Section 2}
\label{sec:prf_upper_bounds}
\bprfof{\Cref{nthm:thm1}} We will use a gradient descent sequence defined
recursively by $\theta^{(0)} = \theta^*$ and, for $t=1,2, \ldots$,
\begin{equation*}
\theta^{(t+1)}=\theta^{(t)}-\eta[\nabla  \ell_{\rho}(\theta^{(t)})+\rho \theta^{(t)}].
\end{equation*}
Our proof builds heavily on the ideas and techniques  developed by
\cite{chen2019spectralregmletopk} and further extended by
\cite{chen2020partialtopkranking}, and contains two key steps. In the first
step, we control $\|\theta^{(T)} - \hat{\theta}_\rho\|$  for $T$  large enough,
by leveraging the convergence property of gradient descent for strong convex
functions. In the second step, we control $\|\theta^{(T)} - \theta^*\|$ through
a leave-one-out argument. The proof can be sketched as follows:
\begin{enumerate}
    \item Bound $\|\theta^{(T)} - \hat{\theta}_\rho\|_\infty$, for large $T$
    using the linear convergence property of gradient descent for
    strongly-convex and smooth functions.
    \item Bound $\|\theta^{(T)} - \theta^*\|_\infty$ for large $T$ using the
    leave-one-out argument.
    \item Finally, $\|\hat{\theta}_\rho - \theta^*\|_\infty$ is controlled by
    triangle inequality.
\end{enumerate}

\textbf{Step 1.} Bound $\|\theta^{(T)} - \hat{\theta}_\rho\|_\infty$, for large
$T$.

\begin{enumerate}
    \item \textbf{Linear convergence, orthogonality to $\mathbf{1}_n$.} We say
that a function $\ell$ is $\alpha$-strongly convex if $\nabla^2\ell (x)\succeq
\alpha I_n$ and $\beta$-smooth if $\|\nabla \ell(x) - \nabla \ell(y)\|_2\leq
\beta\|x - y\|_2$ for all $x,y\in {\rm dom}(\ell)$. By Lemma \ref{lm:lem8}, we
know that $\ell_\rho(\cdot)$ is $\rho$-strongly convex and $(\rho +
n_{\max})$-smooth. By Theorem 3.10 in \cite{Bubeck2015}, we have
\begin{equation}
    \|\theta^{(t)} - \hat{\theta}_\rho\|_2\leq (1 - \frac{\rho}{\rho + n_{\max}})^t\|\theta^{(0)} - \hat{\theta}_\rho\|_2.
    \label{eq:linear_conv}
\end{equation}
Besides, as we start with $\theta^*$ that satisfies $1_n^\top \theta^* = 0$, it
holds that $1_n^\top \theta^{(t)}=0$ for all $t\geq 0$. To see this, just notice
that
\begin{equation*}
\nabla \ell_\rho(\theta) = \rho \theta + \sum_{(i,j)\in E}[-\bar{y}_{ij}+\psi(\theta_i - \theta_j)]({\bm e}_i - {\bm e}_j)
\end{equation*}
and $1_n^\top ({\bm e}_i - {\bm e}_j)$ for any $i,j$. Then by $1_n^\top
\theta^{(t)}=0,\forall t$ and \eqref{eq:linear_conv}, we have $1_n^\top
\hat{\theta}_\rho=0$.

    \item \textbf{Control $\|\theta^* - \hat{\theta}_\rho\|_2$.} By a Taylor
    expansion, we have that
    \begin{align*}
    \ell_\rho(\hat{\theta}_\rho;y)=\ell_\rho(\theta^*;y) &+ (\hat{\theta}_\rho - \theta^*)^\top \nabla \ell_\rho(\theta^*;y)
    \\ +&
    \frac{1}{2}(\hat{\theta}_\rho - \theta^*)^\top \nabla^2 \ell_\rho(\xi;y) (\hat{\theta}_\rho - \theta^*),
    \end{align*}
    where $\xi$ is a convex combination of $\theta^*$ and $\hat{\theta}_{\rho}$.
    By Cauchy-Schwartz inequality,
    \begin{equation*}
    |(\hat{\theta}_\rho - \theta^*)^\top \nabla \ell_\rho(\theta^*;y)|\leq \|\nabla \ell_\rho(\theta^*;y)\|_2 \|\theta^* - \hat{\theta}_\rho\|_2.
    \end{equation*}
    The two inequalities above and the fact that
    $\ell_\rho(\theta^*;y)\geq\ell_\rho(\hat{\theta}_\rho;y)$ yield that
    \begin{equation*}
    \|\theta^* - \hat{\theta}_\rho\|_2\leq \frac{2 \|\nabla \ell_\rho(\theta^*;y)\|_2 }{\rho_{\min} (\nabla^2 \ell_\rho(\xi;y))}.
    \end{equation*}
    By Lemma \ref{lm:lem7}, $\|\nabla \ell_\rho(\theta^*;y)\|_2\lesssim
    \sqrt{\frac{n_{\max} (n+r)}{L}}$. This fact, together with $\rho_{\min}
    (\nabla^2 \ell_\rho(\xi;y))\geq \rho$ and $\rho
    \asymp\frac{1}{\kappa}\sqrt{\frac{n_{\max}}{L}}$, gives that  $\|\theta^* -
    \hat{\theta}_\rho\|_2\leq c \kappa \sqrt{n + r}$, for some $c>0$.
    \item \textbf{Bound $\|\theta^{(T)} - \hat{\theta}_\rho\|_2$.} Take $T =
    \lfloor \kappa^2 e^{3\kappa_E} n^6 \rfloor$ and remember that $L\leq
    \kappa^2 e^{4\kappa_E}n^8$. The previous two steps imply that
        \begin{equation*}
        \|\theta^{(T)} - \hat{\theta}_\rho\|_2\leq c(1 - \frac{\rho}{\rho + n_{\max}})^T \kappa\sqrt{n + r}\leq c\exp\left( - \frac{T\rho}{\rho + n_{\max}}\right)\kappa\sqrt{n + r}.
        \end{equation*}
        Let $\tilde{f}_d =\frac{e^{2\kappa_E}}{\lambda_2(\mathcal{L}_A)}
        \frac{n_{\max}}{n_{\min}}\sqrt{\frac{n + r}{L}} +
        \frac{e^{\kappa_E}}{\lambda_2(\mathcal{L}_A)} \sqrt{\frac{n_{\max} (\log
        n + r)}{L}}$ and consider inequality $e^{-g}\kappa\sqrt{n + r}\leq
        \tilde{f}_d$. The solution is given by $g\geq \log \kappa +
        \frac{1}{2}\log (n+r) -\log\tilde{f}_d$ and the inequality holds as long
        as $g\geq \kappa + 6\log n + 3\log \kappa$ since $L\leq \max\{1,
        \kappa\} e^{3\kappa_E}n^8$. Take $g = \kappa + 5\log n + \log \kappa$,
        then as long as
        \begin{equation*}
            T\rho \geq 2n_{\max} n g,
        \end{equation*}
        it holds that $T \rho > \frac{1}{2}g(\rho + n_{\max})$, then
        $\|\theta^{(T)} - \hat{\theta}_\rho\|_2\leq c\exp(-g)\kappa\sqrt{n+r}$
        is smaller than $\tilde{C}_d\tilde{f}_d$ for some constant
        $\tilde{C}_d$. Since $T=\lfloor \kappa^2 e^{3\kappa_E} n^6 \rfloor$ and
        $\rho \geq \frac{c_{\rho}}{\kappa}\sqrt{\frac{n_{\max}}{\kappa^2
        e^{4\kappa_E}n^8}}$, we have
        \begin{equation*}
            T\rho \geq c_{\rho}\kappa e^{\kappa_E} n^2\sqrt{n_{\max}}\geq 2n_{\max}ng.
        \end{equation*}
        In conclusion, we have
        \begin{equation*}
        \|\theta^{(T)} - \hat{\theta}_\rho\|_\infty\leq \|\theta^{(T)} - \hat{\theta}_\rho\|_2\leq  \tilde{C}_d\tilde{f}_d.
        \end{equation*}
        The arguments above also hold with $\tilde{f}_a =
        \frac{e^{\kappa_E}}{\lambda_2(\mathcal{L}_A)}
        \sqrt{\frac{n_{\max}(n+r)}{L}}$, i.e., we have $\|\theta^{(T)} -
        \hat{\theta}_\rho\|_2\leq \tilde{C}_a\tilde{f}_a$ for some constant
        $\tilde{C}_a$.
\end{enumerate}

\par
\textbf{Step 2.} Bound $\|\theta^{(T)} - {\theta}^*\|_\infty$ by a leave-one-out
argument.
\par
Denote $\psi(x) = \frac{1}{1 + e^{-x}}$ and $r = \log \kappa + \kappa_E$, and
define the leave-one-out negative log-likelihood as
\begin{equation}
\begin{aligned}
\ell_{n}^{(m)}(\theta)=& \sum_{1 \leq i<j \leq n: i, j \neq m} A_{i j}\left[\bar{y}_{i j} \log \frac{1}{\psi\left(\theta_{i}-\theta_{j}\right)}+\left(1-\bar{y}_{i j}\right) \log \frac{1}{1-\psi\left(\theta_{i}-\theta_{j}\right)}\right] \\
&+\sum_{j \in[n] \backslash\{m\}}A_{mj} \left[\psi\left(\theta_{m}^{*}-\theta_{j}^{*}\right) \log \frac{1}{\psi\left(\theta_{m}-\theta_{j}\right)} + \psi\left(\theta_{j}^{*}-\theta_{m}^{*}\right) \log \frac{1}{\psi\left(\theta_{j}-\theta_{m}\right)}\right],
\end{aligned}
\end{equation}
so the leave-one-out gradient descent sequence is, for $t = 0, 1, \ldots$,
\begin{equation*}
\theta^{(t+1, m)}=\theta^{(t, m)}-\eta\left(\nabla \ell_{n}^{(m)}\left(\theta^{(t, m)}\right)+\rho \theta^{(t, m)}\right).
\end{equation*}
We initialize both sequences by $\theta^{(0)} = \theta^{(0,m)} = \theta^*$ and
use step size $\eta = \frac{1}{\rho + n_{\max}}$. By assumption 2,
$\lambda_2(\mclL_{\bfA})>0$, so we can let $f_a=
C_a\frac{e^{\kappa_E}}{\lambda_2(\mathcal{L}_A)}\left[
\sqrt{\frac{n_{\max}(n+r)}{L}} + \rho \kappa(\theta^*)\sqrt{n} \right]$, $f_b =
10e^{\kappa_E}\frac{\sqrt{n_{\max}}}{n_{\min}}f_a$, $f_c =
C_c\frac{e^{\kappa_E}}{\lambda_2(\mathcal{L}_A)}\sqrt{\frac{n_{\max}(\log n +
r)}{L}}$, $f_d = f_b + f_c$ with sufficiently large constant $C_c>0$ and $C_a\gg
C_c$. By assumption 1, we have $f_c + f_d\leq 0.1$. We will show in Lemma
\ref{lm:lem14}, \ref{lm:lem15}, \ref{lm:lem16}, \ref{lm:lem17} that for all
$0\leq t\leq T = \lfloor \kappa^2 e^{3\kappa_E} n^6 \rfloor$
\begin{equation}
\begin{split}
\|\theta^{(t)}-\theta^{*}\|_2 \leq f_a&, \\
\max _{m \in[n]}|\theta_{m}^{(t, m)}-\theta_{m}^{*}| \leq f_b &,\\
\max _{m \in[n]}\|\theta^{(t, m)}-\theta^{(t)}\|_2 \leq f_c &, \\
\|\theta^{(t)} - \theta^*\|_\infty \leq f_d&.
\end{split}
\label{eq:4bound}
\end{equation}
When $t = 0$, \eqref{eq:4bound} holds since $\theta^{(0)} = \theta^{(0,m)} =
\theta^*$. By Lemma \ref{lm:lem14}, \ref{lm:lem15}, \ref{lm:lem16},
\ref{lm:lem17}, and a union bound, we know that \eqref{eq:4bound} holds for all
$0\leq t\leq T= \lfloor \kappa^2 e^{3\kappa_E} n^6 \rfloor$ with probability at
least $1 - O(n^{-4})$. Therefore, using the result in step 1, we have
\begin{equation*}
\|\hat{\theta}_\rho - \theta^*\|_\infty\leq \|\hat{\theta}_\rho - \theta^{(T)}\|_\infty + \|\theta^{(T)} - \theta^*\|_\infty\leq 2f_d.
\end{equation*}
As a byproduct, we have
\begin{equation*}
\|\hat{\theta}_{\rho} - \theta^*\|_2\leq \|\hat{\theta}_{\rho} - \theta^{(T)}\|_2 + \|{\theta}^{(T)} - \theta^*\|_2\leq 2f_a
\end{equation*}
\eprfof

\bnlem \label{lm:lem7} With probability at least $1 -
O(\kappa^{-2}e^{-3\kappa_E}n^{-10})$ the gradient of the regularized
log-likelihood satisfies
\begin{equation*}
\|\nabla \ell_{\rho}(\theta^*)\|_2^2\lesssim \frac{n_{\max}(n+r)}{L} + \rho \kappa(\theta^*)\sqrt{n}.
\end{equation*}
In particular, for $\rho \asymp
\frac{1}{{\kappa(\theta^*)}}\sqrt{\frac{n_{\max}}{L}}$, we have $\|\nabla
\ell_{\rho}(\theta^*)\|_2^2\lesssim \frac{n_{\max}(n +r)}{L}$.
\enlem
\begin{proof}
Triangle inequality gives
\begin{equation*}
\|\nabla \ell_{\rho}(\theta^*)\|_2\leq \|\nabla \ell_{0}(\theta^*)\|_2 + \rho\|\theta^*\|_2.
\end{equation*}
By definition of $\kappa(\theta^*)$, we have $\|\theta^*\|_2\leq
\sqrt{n}\kappa(\theta^*)$. For the first term, by Lemma \ref{lem:lem7.4} we have
\begin{equation*}
\|\nabla \ell_{0}(\theta^*)\|_2^2 = \sum_{i=1}^n\left[\sum_{j\in \mathcal{N}(i)}[\bar{y}_{ij} - \psi(\theta_i^* - \theta_j^*)]\right]^2 \leq C_1\frac{n_{\max}(n + r)}{L}.
\end{equation*}
\end{proof}

\bnlem \label{lm:lem8} Let $\kappa_E(x) = \max_{(i,j)\in E}|x_i - x_j|$, then
$\forall \theta\in\mathbb{R}^n$,
\begin{equation}
    \begin{split}
        \lambda_{\max}(\nabla^2\ell_{\rho}(\theta;y))&\leq \rho + \frac{1}{2}n_{\max},\\
        \lambda_{2}(\nabla^2\ell_{\rho}(\theta;y))&\geq \rho + \frac{1}{4e^{\kappa_E(\theta)}}\lambda_{2}(\mathcal{L}_A).
    \end{split}
\end{equation}
In particular, we have
\begin{equation*}
\lambda_{2}(\nabla^2\ell_{\rho}(\theta;y))\geq \rho + \frac{1}{4e^{\kappa_E(\theta^*)}e^{2\|\theta - \theta^*\|_{\infty}}}\lambda_{2}(\mathcal{L}_A).
\end{equation*}

\enlem
\begin{proof}
Use the fact that
\begin{equation*}
\nabla^2\ell_{0}(\theta;y) = \sum_{(i,j)\in E}\frac{e^{\theta_i}e^{\theta_j}}{(e^{\theta_i} + e^{\theta_j})^2}({\bm e}_i - {\bm e}_j)({\bm e}_i - {\bm e}_j)^\top,
\end{equation*}
and $\forall (i,j)\in E$, $\frac{1}{4\exp({\kappa_E(\theta))}}\leq
\frac{e^{\theta_i}e^{\theta_j}}{(e^{\theta_i} + e^{\theta_j})^2} \leq
\frac{1}{4}$, $\kappa_E(x_1)\leq \kappa_E(x_2) + 2\|x_1 - x_2\|_\infty$. In
addition, the largest eigenvalue of graph Laplacian satisfies (Corollary 3.9.2
in \cite{spectra2012})
\[
\lambda_{\max}(\mathcal{L}_A) \leq \max_{(i,j)\in E}(n_i+n_j)\leq 2n_{\max}.
\]
\end{proof}

\bnlem
Provided that \eqref{eq:4bound} holds, then
\begin{equation}
\begin{split}
    \max_{m\in [n]}\|\theta^{(t+1,m)}-\theta^*\|_\infty &\leq f_c+f_d,\\
    \max_{m\in [n]}\|\theta^{(t+1,m)}-\theta^*\|_2 &\leq f_c+f_a.
\end{split}
\end{equation}
\label{lm:lem13}
\enlem
\begin{proof}
By triangle inequality, we have
\begin{align*}
    \max_{m\in [n]}\|\theta^{(t,m)}-\theta^*\|_\infty &\leq \max_{m\in [n]}\|\theta^{(t,m)}-\theta^{(t)}\|_\infty + \|\theta^{(t)}-\theta^{*}\|_\infty\leq f_c+f_d\\
    \max_{m\in [n]}\|\theta^{(t,m)}-\theta^*\|_2 &\leq \max_{m\in [n]}\|\theta^{(t,m)}-\theta^{(t)}\|_2 + \|\theta^{(t)}-\theta^{*}\|_2\leq f_c+f_a.
\end{align*}
\end{proof}

\bnlem \label{lm:lem14} Suppose \eqref{eq:4bound} holds, and the step size
satisfies $0<\eta \leq \frac{1}{\rho + n_{\max}}$. If
\begin{equation*}
f_d\leq 0.1 \text{ and } f_a\geq C_a\frac{e^{\kappa_E(\theta^*)}}{\lambda_2(\mathcal{L}_A)}\left[ \sqrt{\frac{n_{\max}(n+r)}{L}} + \rho \kappa(\theta^*)\sqrt{n} \right],
\end{equation*}
for some large constant $C_a$, then with probability at least $1 -
O(\kappa^{-2}e^{-3\kappa_E}n^{-10})$ we have
\begin{equation*}
\|\theta^{(t+1)}-\theta^{*}\|_2 \leq f_a.
\end{equation*}
\enlem
\begin{proof}
By the form of the gradient descent, we have that
\begin{align*}
    \theta^{(t+1)} - \theta^* &= \theta^{(t)} - \eta\nabla\ell_{\rho}(\theta^{(t)}) - \theta^*\\
    & = \theta^{(t)} - \eta\nabla\ell_{\rho}(\theta^{(t)}) - [\theta^{*} - \eta\nabla\ell_{\rho}(\theta^{*})] - \eta\nabla\ell_{\rho}(\theta^{*})\\
    &= \left[ I_n - \eta\int_0^1 \nabla^2\ell_\rho(\theta(\tau)){\rm d}\tau  \right](\theta^{(t)} - \theta^*) - \eta\nabla\ell_{\rho}(\theta^{*}),
\end{align*}
where $\theta(\tau) = \theta^* + \tau(\theta^{(t)} - \theta^*)$. Letting
$H=\int_0^1 \nabla^2\ell_\rho(\theta(\tau)){\rm d}\tau$, by the triangle
inequality,
\[
\|\theta^{(t+1)} - \theta^*\|_2\leq \|(I_n - \eta H)(\theta^{(t)} - \theta^*)\|_2 + \eta \|\nabla \ell_{\rho}(\theta^*)\|_2.
\]
Setting $\kappa_E(x) = \max_{(i,j)\in E}|x_i - x_j|$, then, for sufficiently
small $\epsilon$, we have that
\[
\kappa_E(\theta(\tau)) \leq \kappa_E(\theta^*) + 2\|\theta^{(t)} - \theta^*\|_\infty\leq \kappa_E(\theta^*) + \epsilon.
\]
as long as
\begin{equation}
2f_d\leq \epsilon.
\end{equation}

Then, by Lemma \ref{lm:lem8} and setting $\epsilon=0.2$, for any $\tau\in
[0,1]$,
\begin{equation}
    \rho + \frac{\lambda_2(\mathcal{L}_A)}{10 e^{\kappa_E(\theta^*)}}\leq \rho + \frac{\lambda_2(\mathcal{L}_A)}{8 e^{\kappa_E(\theta^*)}e^{\epsilon}} \leq \lambda_{2}(\nabla^2\ell_{\rho}(\theta(\tau)))\leq \lambda_{\max}(\nabla^2\ell_{\rho}(\theta(\tau)))\leq \rho + \frac{1}{2}n_{\max}.
    \label{eq:bound_hess}
\end{equation}
Since $1_n^\top (\theta^{(t)} - \theta^*) = 0$, we obtain that
\[
\|(I_n - \eta H)(\theta^{(t)} - \theta^*)\|_2\leq \max\{|1 -\eta\lambda_{2}(H)|,|1 -\eta\lambda_{\max}(H)|\}\|\theta^{(t)} - \theta^*\|_2.
\]
By \eqref{eq:bound_hess} and the fact that $\eta \leq \frac{1}{\rho +
n_{\max}}$, we get
\begin{equation}
    \|(I_n - \eta H)(\theta^{(t)} - \theta^*)\|_2\leq (1 - \frac{\eta \lambda_2(\mathcal{L}_A)}{10e^{\kappa_E(\theta^*)}})\|\theta^{(t)} - \theta^*\|_2.
    \label{eq:v1}
\end{equation}
By Lemma \ref{lm:lem7} and the induction hypothesis, we have
\[
\|\theta^{(t+1)} - \theta^*\|_2\leq (1 - \frac{\eta \lambda_2(\mathcal{L}_A)}{10e^{\kappa_E(\theta^*)}})f_a + C\eta\left[ \sqrt{\frac{n_{\max}(n+r)}{L}} + \rho \kappa(\theta^*)\sqrt{n} \right]\leq f_a
\]
as long as
\[
f_a\geq C_a\frac{e^{\kappa_E(\theta^*)}}{\lambda_2(\mathcal{L}_A)}\left[ \sqrt{\frac{n_{\max}(n + r)}{L}} + \rho \kappa(\theta^*)\sqrt{n} \right]
\]
for some large constant $C_a$.
\end{proof}

\bnlem \label{lm:lem15} Suppose \eqref{eq:4bound} holds and assume that
\begin{enumerate}
    \item  $f_a=
    C_a\frac{e^{\kappa_E(\theta^*)}}{\lambda_2(\mathcal{L}_A)}\left[
    \sqrt{\frac{n_{\max}(n+r)}{L}} + \rho \kappa(\theta^*)\sqrt{n} \right]$,
    $f_c =
    C_c\frac{e^{\kappa_E(\theta^*)}}{\lambda_2(\mathcal{L}_A)}\sqrt{\frac{n_{\max}(\log
    n + r)}{L}}$ with $C_a\gg C_c$.
    \item $\frac{n_{\min}}{10e^{\kappa_E(\theta^*)}}f_b\geq
    \frac{3\sqrt{n_{\max}}}{4}f_a$.
    \item $f_c + f_d\leq 0.1$.
\end{enumerate}
then as long as the step size satisfies $0<\eta \leq \frac{1}{\rho + n_{\max}}$,
with probability at least $1 - O(\kappa^{-2}e^{-3\kappa_E}n^{-10})$ we have
\begin{equation*}
\max _{m \in[n]}\left|\theta_{m}^{(t+1, m)}-\theta_{m}^{*}\right| \leq f_b.
\end{equation*}
\enlem
\begin{proof}
Recall that the gradient descent step for leave-one-out estimator $\theta^{(m)}$
is defined as
\begin{equation*}
    \theta^{(t+1, m)}=\theta^{(t, m)}-\eta\left(\nabla \ell_{n}^{(m)}\left(\theta^{(t, m)}\right)+\rho \theta^{(t, m)}\right),
\end{equation*}
where
\begin{equation*}
\begin{aligned}
\ell_{n}^{(m)}(\theta)=& \sum_{1 \leq i<j \leq n: i, j \neq m} A_{i j}\left[\bar{y}_{i j} \log \frac{1}{\psi\left(\theta_{i}-\theta_{j}\right)}+\left(1-\bar{y}_{i j}\right) \log \frac{1}{1-\psi\left(\theta_{i}-\theta_{j}\right)}\right] \\
&+\sum_{j \in[n] \backslash\{m\}}A_{mj} \left[\psi\left(\theta_{m}^{*}-\theta_{j}^{*}\right) \log \frac{1}{\psi\left(\theta_{m}-\theta_{j}\right)}+\psi\left(\theta_{j}^{*}-\theta_{m}^{*}\right) \log \frac{1}{\psi\left(\theta_{j}-\theta_{m}\right)}\right].
\end{aligned}
\end{equation*}
Direct calculations give
\begin{equation*}
\begin{split}
[\nabla \ell_{n}^{(m)}(\theta)]_m&= \sum_{j\in [n]\setminus \{m\}}A_{mj}\left[ \psi(\theta^*_m - \theta^*_j)(\psi(\theta^*_m - \theta^*_j) - 1) + (1 - \psi(\theta^*_m - \theta^*_j))\psi(\theta_m - \theta_j) \right]\\
&= \sum_{j\in [n]\setminus \{m\}}A_{mj}\left[ -\psi(\theta^*_m - \theta^*_j) + \psi(\theta_m - \theta_j) \right].
\end{split}
\end{equation*}
Thus, we have
\begin{equation*}
\theta^{(t+1,m)}_m - \theta_m^* = \left( 1 - \eta\rho -\eta \sum_{j\in [n]\setminus \{m\} }A_{mj}\psi{'} (\xi_j) \right)(\theta^{(t,m)}_m - \theta_m^*) - \rho \eta \theta_m^* + \eta  \sum_{j\in [n]\setminus \{m\}}A_{mj}\psi{'}(\xi_j) (\theta^{(t,m)}_j - \theta_j^*),
\end{equation*}
where $\xi_j$ is a scalar between $\theta_m^*-\theta_j^*$ and $\theta^{(t,m)}_m
- \theta^{(t,m)}_j$. Notice that $\psi'(x) = \frac{e^x}{(1 + e^x)^2}\leq
\frac{1}{4}$ for any $c\in \mathbb{R}$, thus by Cauchy-Schwartz inequality we
have
\begin{equation}
    |\sum_{j\in [n]\setminus \{m\}}A_{mj}\psi{'}(\xi_j) (\theta^{(t,m)}_j - \theta_j^*)|\leq \frac{1}{4} \sqrt{n_{\max}} \|\theta^{(t,m)}-\theta^*\|_2.
\label{eq:key_difference}
\end{equation}
Also, since $\eta\leq \frac{1}{\rho + n_{\max}}$,
\begin{equation*}
1 - \eta\rho -\eta  \sum_{j\in [n]\setminus \{m\} }A_{mj}\psi{'} (\xi_j)\geq 1 - \eta\rho - \eta \frac{n_{\max}}{4}\geq 0.
\end{equation*}
Therefore,
\begin{equation*}
0\leq 1 - \eta\rho -\eta  \sum_{j\in [n]\setminus \{m\} }A_{mj}\psi{'} (\xi_j)\leq 1-\eta n_{\min}\min_{j\in \mathcal{N}(m)}\psi'(\xi_j).
\end{equation*}
Since $\xi_j$ is a scalar between $\theta_m^*-\theta_j^*$ and $\theta^{(t,m)}_m
- \theta^{(t,m)}_j$, we have
\begin{align*}
    \max_{j\in \mathcal{N}(m)}|\xi_j| &\leq \max_{j\in \mathcal{N}(m)}|\theta_m^* - \theta_j^*| + \max_{j\in\mathcal{N}(m)}|\theta_m^* - \theta_j^* - (\theta_m^{(t,m)} - \theta_j^{(t,m)})|\\
    &\leq \kappa_E(\theta^*) + 2\|\theta^{(t,m)} - \theta^*\|_\infty \leq \kappa_E(\theta^*) + \epsilon
\end{align*}
as long as
\begin{equation}
    \|\theta^{(t,m)} - \theta^*\|_\infty\leq f_c + f_d\leq \epsilon/2.
\end{equation}
Let $\epsilon = 0.2$, then $e^{\epsilon}\leq 5/4$ and
\begin{equation*}
\psi'(\xi_j) = \frac{e^{\xi_j}}{(1+e^{\xi_j})^2} = \frac{e^{-|\xi_j|}}{(1+e^{-|\xi_j|})^2}\geq \frac{e^{-|\xi_j|}}{4}\geq \frac{1}{4e^{\epsilon + \kappa_E(\theta^*)}}\geq \frac{1}{5e^{ \kappa_E(\theta^*)}}.
\end{equation*}
By triangle inequality we get
\begin{equation}
    \begin{split}
        |\theta^{(t+1,m)}_m - \theta_m^*| &\leq \left( 1 - \frac{\eta n_{\min}}{10e^{\kappa_E(\theta^*)}} \right)|\theta^{(t,m)}_m - \theta_m^*| + \rho \eta \|\theta^*\|_\infty + \frac{\eta \sqrt{n_{\max}}}{4}\|\theta^{(t,m)} - \theta^*\|_2\\
        &\leq f_b - \frac{\eta n_{\min}}{10e^{\kappa_E(\theta^*)}}f_b + \eta \rho \kappa(\theta^*) + \eta\frac{\sqrt{n_{\max}}}{4}(f_a + f_c)\leq f_b
    \end{split}
\end{equation}
as long as
\begin{equation*}
\frac{n_{\min}}{10e^{\kappa_E(\theta^*)}}f_b \geq \rho \kappa(\theta^*) + \frac{\sqrt{n_{\max}}}{4}(f_a + f_c).
\end{equation*}
By assumption, $f_a=
C_a\frac{e^{\kappa_E(\theta^*)}}{\lambda_2(\mathcal{L}_A)}\left[
\sqrt{\frac{n_{\max}(n + r)}{L}} + \rho \kappa(\theta^*)\sqrt{n} \right]$, $f_c
=
C_c\frac{e^{\kappa_E(\theta^*)}}{\lambda_2(\mathcal{L}_A)}\sqrt{\frac{n_{\max}(\log
n + r)}{L}}$ with $C_a\gg \max\{C_c,1\}$, so
\begin{equation*}
f_a\gg f_c,\text{ and } \frac{\sqrt{n_{\max}}}{4}f_a \gg \frac{n_{\max}}{\lambda_2(\mathcal{L}_A)} \left[\sqrt{\frac{n+r}{L}} + \sqrt{\frac{n}{n_{\max}}}\rho\kappa(\theta^*)\right]\geq \rho\kappa(\theta^*).
\end{equation*}
Therefore, a sufficient condition for $|\theta^{(t+1,m)}_m - \theta_m^*|\leq
f_b$ is
\begin{equation*}
\frac{n_{\min}}{10e^{\kappa_E(\theta^*)}}f_b\geq \frac{3\sqrt{n_{\max}}}{4}f_a,
\end{equation*}
which is satisfied by our assumption.
\end{proof}

\bnlem \label{lm:lem16} Suppose \eqref{eq:4bound} holds with $f_c=
C_c\frac{e^{\kappa_E(\theta^*)}}{\lambda_2(\mathcal{L}_A)}\sqrt{\frac{n_{\max}(\log
n+r)}{L}}$ for some sufficiently large constant $C_c$, $f_d = f_b + f_c$, then
as long as the step size satisfies $0<\eta \leq \frac{1}{\rho + n_{\max}}$, with
probability at least $1 - O(\kappa^{-2}e^{-3\kappa_E}n^{-10})$ we have
\begin{equation*}
\max _{m \in[n]}\|\theta^{(t+1, m)}-\theta^{(t)}\|_2 \leq f_c.
\end{equation*}

\enlem

\begin{proof}
By the update rules, we have
\begin{align*}
    \theta^{(t+1)} - \theta^{(t+1,m)} =& \theta^{(t)} - \eta\nabla\ell_{\rho}(\theta^{(t)}) - \left[ \theta^{(t,m)} - \eta\nabla\ell_{\rho}^{(m)}(\theta^{(t,m)}) \right]\\
    =& \theta^{(t)} - \eta\nabla\ell_{\rho}(\theta^{(t)}) - \left[ \theta^{(t,m)} - \eta\nabla\ell_{\rho}(\theta^{(t,m)}) \right]\\
    &-\eta\left[ \nabla\ell_{\rho}(\theta^{(t,m)}) - \nabla\ell_{\rho}^{(m)}(\theta^{(t,m)}) \right]\\
    =& v_1 - v_2,
\end{align*}
where
\begin{equation*}
v_1 = \left[ I_n - \eta\int_0^1 \nabla^2\ell_\rho(\theta(\tau)){\rm d}\tau  \right](\theta^{(t)} - \theta^{(t,m)}),\quad v_2 = \eta\left[ \nabla\ell_{\rho}(\theta^{(t,m)}) - \nabla\ell_{\rho}^{(m)}(\theta^{(t,m)}) \right]
\end{equation*}
Now following the same arguments towards \eqref{eq:v1}, as long as $\eta \leq
\frac{1}{\rho + n_{\max}}$, we can get
\begin{equation*}
\|v_1\|_2\leq (1 - \frac{\eta \lambda_2(\mathcal{L}_A)}{10e^{\kappa_E(\theta^*)}})\|\theta^{(t)} - \theta^{(t,m)}\|_2.
\end{equation*}
For $v_2$, we know that
\begin{equation*}
\begin{split}
    \frac{1}{\eta}v_2 &= \sum_{i\in [n]\setminus\{m\}} \left\lbrace A_{mi}\left[ \psi(\theta_i^{(t,m)} - \theta_m^{(t,m)}) -\bar{y}_{im} \right] - A_{mi} \left[ \psi(\theta_i^{(t,m)} - \theta_m^{(t,m)}) - \psi(\theta_i^* - \theta_m^*) \right] \right\rbrace ({\bm e}_i - {\bm e}_m) \\
    &= \sum_{i\in [n]\setminus\{m\}} A_{mi}\left[ \psi(\theta_i^* - \theta_m^*) -\bar{y}_{im} \right] ({\bm e}_i - {\bm e}_m).
\end{split}
\end{equation*}
By the form of the derivatives and Lemma \ref{lem:lem7.4}, we know that with
probability at least $1 - O(n^{\kappa^{-2}e^{-3\kappa_E}n^{-10}})$,
\begin{equation}
\begin{aligned}
\|\frac{1}{\eta}v_2\|_2^2=&\left[\sum_{i \in[n] \backslash\{m\}} A_{i m}\left(\bar{y}_{i m}-\psi\left(\theta_{i}^{*}-\theta_{m}^{*}\right)\right)\right]^{2} + \sum_{i \in[n] \backslash\{m\}} A_{i m}\left(\bar{y}_{i m}-\psi\left(\theta_{i}^{*}-\theta_{m}^{*}\right)\right)^{2}\\
\lesssim & \frac{n_{\max}(\log n + r)}{L} + \frac{\log n + n_{\max} + r}{L} \lesssim  \frac{n_{\max}(\log n + r)}{L} .
\end{aligned}
\end{equation}
Therefore, we have
\begin{equation}
    \begin{split}
        \|\theta^{(t+1)} - \theta^{(t+1,m)}\|_2 &\leq \|v_1\|_2 + \|v_2\|_2\\
        &\leq (1 - \frac{\eta \lambda_2(\mathcal{L}_A)}{10e^{\kappa_E(\theta^*)}})\|\theta^{(t)} - \theta^{(t,m)}\|_2 + C\eta \sqrt{\frac{n_{\max}(\log n + r)}{L} }\\
        &\leq (1 - \frac{\eta \lambda_2(\mathcal{L}_A)}{10e^{\kappa_E(\theta^*)}})f_c + C\eta \sqrt{\frac{n_{\max}(\log n + r)}{L} }\leq f_c,
    \end{split}
\end{equation}
where the last inequality is due to the fact that $C_c$ is a sufficiently large
constant by our assumption and
\begin{equation}
    \frac{\eta \lambda_2(\mathcal{L}_A)}{30e^{\kappa_E(\theta^*)}}f_c \geq C\eta \sqrt{\frac{n_{\max}(\log n + r)}{L} }\Leftarrow f_c= C_c\frac{e^{\kappa_E(\theta^*)}}{\lambda_2(\mathcal{L}_A)}\sqrt{\frac{n_{\max}(\log n + r)}{L}}.
\end{equation}
\end{proof}

\bnlem \label{lm:lem17} Suppose \eqref{eq:4bound} holds and $f_d\geq f_b + f_c$,
then with probability at least $1 - O(\kappa^{-2}e^{-3\kappa_E}n^{-10})$ we have
\begin{equation*}
\|\theta^{(t+1)} - \theta^*\|_\infty \leq f_d.
\end{equation*}

\enlem
\begin{proof}
By Lemma \ref{lm:lem15} and Lemma \ref{lm:lem16} we have
\begin{align*}
|\theta_{m}^{(t+1)}-\theta_{m}^{*}| & \leq|\theta_{m}^{(t+1)}-\theta_{m}^{(t+1,m)}|+|\theta_{m}^{(t+1,m)}-\theta_{m}^{*}| \\
& \leq\|\theta^{(t+1)}-\theta^{(t+1,m)}\|_{2}+\mid \theta_{m}^{(t+1,m)}-\theta_{m}^{*}|\leq f_c + f_b\leq f_d,
\end{align*}
since $f_d\geq f_b + f_c$ by our assumption.
\end{proof}

\bnlem  \label{lem:lem7.4} With probability at least
$1-O(\kappa^{-2}e^{-3\kappa_E}n^{-10})$ it holds that
\begin{equation}
    \begin{split}
        \max_{i\in [n]}\left[\sum_{j\in \mathcal{N}(i)}[\bar{y}_{ij} - \psi(\theta_i^* - \theta_j^*)]\right]^2 &\leq C\frac{(\log n + r) \cdot n_{\max}}{L},\\
        \sum_{i=1}^n\left[\sum_{j\in \mathcal{N}(i)}[\bar{y}_{ij} - \psi(\theta_i^* - \theta_j^*)]\right]^2 &\leq C\frac{(n + r)\cdot n_{\max}}{L}.\\
        \max_{i\in [n]}\sum_{j\in \mathcal{N}(i)}[\bar{y}_{ij} - \psi(\theta_i^* - \theta_j^*)]^2 &\leq C\frac{\log n + n_{\max}+r}{L},
    \end{split}
\end{equation}
where $r = \log\kappa + \kappa_E$.
\enlem
\begin{proof}
To prove the first inequality, notice that by Hoeffding's inequality we have
\begin{equation*}
\mathbb{P}\left(\sum_{j< i}A_{ij}[\bar{y}_{ij} - \psi(\theta_i^* - \theta_j^*)]\geq \sqrt{\frac{8n_{\max}(\log n + r)}{L}}\right)\leq 2\exp(-\frac{2L}{n_{\max}}\cdot \frac{8n_{\max} (\log n + r)}{L}) = 2\kappa^{-2}e^{-3\kappa_E}n^{-12},
\end{equation*}
where $r = \log \kappa + \kappa_E$. By union bound we know that on an event $B$
with probability at least $1 - \kappa^{-2}e^{-3\kappa_E}n^{-10}$ we have
$\forall i\in [n]$, $\left[\sum_{j< i}A_{ij}[\bar{y}_{ij} - \psi(\theta_i^* -
\theta_j^*)]\right]^2<\frac{8n_{\max} (\log n +r)}{L}$.

Next we prove the second inequality. Consider the unit ball $\mclS = \{v\in
\mathbb{R}^n:\sum_{i\in [n]}v_i^2 = 1\}$ in $\mathbb{R}^n$. By Lemma 5.2 of
\cite{vershynin2011introduction}, we can pick a subset $\mclU\subset \mclS$ so
that $\log |\mclU|\leq cn$ and for any $v\in \mclS$, there exists a vector $u
\in \mclU$ such that $\|u - v\|_2\leq \frac{1}{2}$. For a given $v\in \mclS$,
pick $u\in \mclU$ such that $\|u - v\|_2\leq \frac{1}{2}$ and we have
\begin{align*}
    \sum_{i=1}^n v_i \left[\sum_{j\in \mathcal{N}(i)}[\bar{y}_{ij} - \psi(\theta_i^* - \theta_j^*)]\right] &= \sum_{i=1}^n u_i \left[\sum_{j\in \mathcal{N}(i)}[\bar{y}_{ij} - \psi(\theta_i^* - \theta_j^*)]\right] + \sum_{i=1}^n (v_i - u_i)\left[\sum_{j\in \mathcal{N}(i)}[\bar{y}_{ij} - \psi(\theta_i^* - \theta_j^*)]\right] \\
    &\leq \sum_{i=1}^n u_i\left[\sum_{j\in \mathcal{N}(i)}[\bar{y}_{ij} - \psi(\theta_i^* - \theta_j^*)]\right] + \frac{1}{2}\sqrt{\sum_{i=1}^n  \left[\sum_{j\in \mathcal{N}(i)}[\bar{y}_{ij} - \psi(\theta_i^* - \theta_j^*)]\right]^2}.
\end{align*}
Taking maximum over $v$ and the left hand side can achieve $\sqrt{\sum_{i=1}^n
\left[\sum_{j\in \mathcal{N}(i)}[\bar{y}_{ij} - \psi(\theta_i^* -
\theta_j^*)]\right]^2}$, thus we have
\begin{align*}
        \sqrt{\sum_{i=1}^n \left[\sum_{j\in \mathcal{N}(i)}[\bar{y}_{ij} - \psi(\theta_i^* - \theta_j^*)]\right]^2} &\leq 2\max_{u\in \mclU}\sum_{i=1}^n u_i\left[\sum_{j\in \mathcal{N}(i)}[\bar{y}_{ij} - \psi(\theta_i^* - \theta_j^*)]\right]\\
        &= 2\max_{u\in \mclU}\sum_{i<j} A_{ij}(u_i - u_j)\left[\bar{y}_{ij} - \psi(\theta_i^* - \theta_j^*)\right].
\end{align*}
To apply Hoeffding's inequality and union bound, we should account for
$|\mclU|\leq e^{cn}$, so we can get that with probability at least $1 -
\kappa^{-2}e^{-3\kappa_E}n^{-10}$ it holds that
\begin{align*}
    \sum_{i=1}^n \left[\sum_{j\in \mathcal{N}(i)}[\bar{y}_{ij} - \psi(\theta_i^* - \theta_j^*)]\right]^2 &\leq C_1\frac{1}{L}\left[ (\log n + n + r)\sum_{i<j}A_{ij}(u_i - u_j)^2\right]\\
    &\leq C_1\frac{(\log n + n + r)\lambda_{\max}(\mclL_A)}{L},
\end{align*}
where $r= \log \kappa + \kappa_E$. Since $\lambda_{\max}(\mclL_A)\leq
2n_{\max}$, the second inequality is proved.

Next we prove the third inequality. For each $i\in [n]$, let $\mclV_i:=\{v\in
\mathbb{R}^{n-1}:\sum_{j\neq i}A_{ij}v_j^2\leq 1\}$. By Lemma 5.2 of
\cite{vershynin2011introduction}, we can pick a subset $\mclU_i\subset \mclV_i$
so that $\log |\mclU_i|\leq 2\sum_{j\neq i}A_{ij}$ and for any $v\in \mclV_i$,
there exists a vector $u \in \mclU_i$ such that $\|u - v\|_2\leq \frac{1}{2}$.
For a given $v\in \mclV_i$, pick $u\in \mclU_i$ such that $\|u - v\|_2\leq
\frac{1}{2}$ and we have
\begin{equation*}
    \begin{split}
        \sum_{j\in \mathcal{N}(i)}v_{ij}[\bar{y}_{ij} - \psi(\theta_i^* - \theta_j^*)] &= \sum_{j\in \mathcal{N}(i)}u_{ij}[\bar{y}_{ij} - \psi(\theta_i^* - \theta_j^*)] + \sum_{j\in \mathcal{N}(i)}(v_{ij} - u_{ij})[\bar{y}_{ij} - \psi(\theta_i^* - \theta_j^*)]\\
        &\leq \sum_{j\in \mathcal{N}(i)}u_{ij}[\bar{y}_{ij} - \psi(\theta_i^* - \theta_j^*)] + \frac{1}{2}\sqrt{\sum_{j\in \mathcal{N}(i)}[\bar{y}_{ij} - \psi(\theta_i^* - \theta_j^*)]^2}
    \end{split}
\end{equation*}
Taking maximum over $v$ and the left hand side can achieve $\sqrt{\sum_{j\in
\mathcal{N}(i)}[\bar{y}_{ij} - \psi(\theta_i^* - \theta_j^*)]^2}$, thus we have
\begin{equation*}
    \sqrt{\sum_{j\in \mathcal{N}(i)}[\bar{y}_{ij} - \psi(\theta_i^* - \theta_j^*)]^2}\leq 2\max_{u\in \mclU_i}\sum_{j\in \mathcal{N}(i)}u_{ij}[\bar{y}_{ij} - \psi(\theta_i^* - \theta_j^*)].
\end{equation*}
Therefore,
\begin{equation*}
    \sqrt{\max_{i\in [n]}\sum_{j\in \mathcal{N}(i)}[\bar{y}_{ij} - \psi(\theta_i^* - \theta_j^*)]^2}\leq 2\max_{i\in [n]}\max_{u\in \mclU_i}\sum_{j\in \mathcal{N}(i)}u_{ij}[\bar{y}_{ij} - \psi(\theta_i^* - \theta_j^*)].
\end{equation*}
Now a straightforward application of Hoeffding's inequality and union bound
gives
\begin{equation*}
    \max_{i\in [n]}\sum_{j\in \mathcal{N}(i)}[\bar{y}_{ij} - \psi(\theta_i^* - \theta_j^*)]^2\leq C\frac{1}{L}\left[\log n + n_{\max} + \kappa_E + \log \kappa\right]
\end{equation*}
with probability at least $1 - O(\kappa^{-2}e^{-3\kappa_E}n^{-10})$.
\end{proof}

\textbf{Corollary 2} (Erd\"os-R\'enyi graph)\textbf{.} Suppose that the
comparison graph comes from an Erd\"os-R\'enyi graph $ER(n,p)$. Assume that
$\mathbf{1}_n^\top \bbrtheta^* = 0$, $\kappa \leq n$, $\kappa_E\leq\log n$,
$L\leq n^8e^{4\kappa_E}\max\{1,\kappa\}$, $np>C_1\log n$, and $L\geq C_2\
\max\{1,\kappa\} e^{4\kappa_E}n/\log^2 n$ for some sufficiently large constants
$C_1,C_2>0$. Set $\rho = c_{\rho}/\kappa \sqrt{n_{max}/L}$. Then
$\mathbf{1}_n^\top \hat{\bbrtheta}_\rho = 0$, and with probability at least $1 -
O(n^{-4})$, it holds that
\begin{equation}
    \begin{split}
        \|\hat{\bbrtheta}_{\rho} - \bbrtheta^*\|_{\infty} \lesssim & e^{2\kappa_E} \sqrt{\frac{1}{np^2L}} + e^{\kappa_E} \sqrt{\frac{\log n}{npL}}, \\
        \|\hat{\bbrtheta}_{\rho} - \bbrtheta^*\|_2 \lesssim &  {e^{\kappa_E}} \sqrt{\frac{1}{pL}}.
    \end{split}
\end{equation}
\bprfof{\Cref{cor:cor_ER}} For an $ER(n,p)$ graph $\mclG$ with $p\geq
c\frac{\log n}{n}$ for some larege $c>0$, it holds with probability at least $1
- O(n^{-10})$ that $\mclG$ is connected, and
\begin{equation*}
    \frac{1}{2} n p \leq n_{\min} \leq n_{\max} \leq 2 n p,
\end{equation*}
and
\begin{equation*}
    \lambda_{2}\left(\mathcal{L}_{A}\right)=\min _{u \neq 0: \mathbf{1}_{n}^{\top} u=0} \frac{u^{\top} \mathcal{L}_{A} u}{\|u\|^{2}} \geq \frac{n p}{2}.
\end{equation*}
The proof can be seen in either \cite{chen2019spectralregmletopk} or
\cite{chen2020partialtopkranking}. Thus, by a union bound, we can replace the
corresponding quantities in upper bounds in \Cref{nthm:thm1} and get the
high probability bounds in \Cref{cor:cor_ER}.
\eprfof

\textbf{Proposition 3} (Path graph)\textbf{.} Suppose the comparison graph is a
path graph $([n],E)$ with $E = \{(i, i+1)\}_{i\in [n - 1]}$ and for each $i$,
item $i$ and item $i+1$ are compared $L_{i,i+1}$ times such that
$\min_{i}L_{i,i+1}>c e^{2\kappa_E}n\log n$ for some universal constant $c$, then
with probability at least $1 - n^{-4}$, the vanilla MLE $\hat{\bbrtheta}_0$
satisfies
\begin{equation}
    \|\hat{\bbrtheta}_0 - \bbrtheta^*\|_{\infty}\lesssim \sqrt{\sum_{i=1}^{n-1}\frac{\exp(2|\theta^*_i - \theta^*_{i+1}|)\log n}{L_{i,i+1}}}.
\end{equation}
In particular, when $L_{i,i+1}=L$ for all $i\in [n - 1]$, we have
\begin{equation}
    \|\hat{\bbrtheta}_0 - \bbrtheta^*\|_{\infty}\lesssim e^{\kappa_E}\sqrt{\frac{n\log n}{L}}, \|\hat{\bbrtheta}_0 - \bbrtheta^*\|_{2}\lesssim e^{\kappa_E}n\sqrt{\frac{\log n}{L}}
\end{equation}
\bprfof{\Cref{prop:path}} Consider a path graph with edge set $E =
\{(i,i+1):i\in [n-1]\}$. Let $M_{ij}$ be the number of wins of item $i$ against
item $j$. For the ease of notations, in this proof we use $\hat{\theta}$ to
denote the vanilla MLE $\hat{\bbrtheta}_0$. We know that
\begin{equation*}
\nabla \ell\left(\theta\right)_i=\sum_{j \in\mathcal{N}(i) } \left(\frac{M_{ij}}{M_{ij} +M_{ji}}-\frac{\exp(\theta_{i}-\theta_{j})}{1 + \exp(\theta_{i}-\theta_{j})}\right).
\end{equation*}
Thus the vanilla MLE solving $\nabla \ell(\theta) = 0$ is given by
\begin{equation*}
\hat{\theta}_{i+1} - \hat{\theta}_{i} = \log M_{i+1,i} - \log M_{i,i+1} := \log R_{i+1,i},
\end{equation*}
where $M_{ij}:=\#\{i \text{ beats } j\}$ and $R_{ij}:= \frac{M_{ij}}{M_{ji}}$.

Let $\hat{\theta}_1 =0$ and we can get
\begin{equation*}
\hat{\theta}_1 = 0,\ \hat{\theta}_{i+1} = \sum_{j=1}^{i} (\log R_{j+1,j} - \log r_{j+1,j}),\ i = 1,\cdots,n-1,
\end{equation*}
where $\log R_{j+1,j}:= \log M_{j+1,j} - \log M_{j,j+1}$, $\log r_{j+1,j} =
\theta^*_{j+1} - \theta^*_{j}$ . Shifting $\theta^*$ to make $\theta^*_1 = 0$
and we have
\begin{equation}
|\hat{\theta}_{i+1} - {\theta}^*_{i+1}| =\left|\sum_{j=1}^{i} (\log R_{j+1,j} - \log r_{j+1,j}) \right|
\label{eq:err-closed_form}
\end{equation}
Let $F_{ij}:=\frac{M_{ij}}{L_{ij}}$ with $L_{ij} = M_{ij} + M_{ji}$ and we have
\begin{equation*}
\log R_{ij} - \log r_{ij} = \log \frac{1 - F_{ji}}{F_{ji}} - \log \frac{1 - p_{ji}}{p_{ji}}.
\end{equation*}
Now using the Talyor expansion of $f(x) = \log (\frac{1}{x} - 1)$ at $p_{ji}$,
we can get
\begin{equation*}
\log R_{ij} - \log r_{ij} = -v_{ij}(F_{ij} - p_{ij}) + \frac{1}{2}\frac{1 - 2z_{ji}}{z_{ji}^2(1 - z_{ji})^2}(F_{ij} - p_{ij})^2,
\end{equation*}
where $v_{ij} := [p_{ij}(1 - p_{ij})]^{-1} = 2 + r_{ij} + r_{ji}=v_{ji}\leq
4\exp(|\theta^*_i - \theta^*_j|)$, and $z_{ji}$ is a number between $p_{ji}$ and
$F_{ji}$. We know that $v_{ij}(F_{ij} - p_{ij})$ is a
subgaussian-$\frac{v_{ij}^2}{4L_{ij}}$ variable. In particular, for the path
graph, by a standard sub-Gaussian tail bound on the first term, we can show that
with probability at least $1 - n^{-9}$, it holds for $i = 1,\cdots,n-1$
simultaneously and some constant $C>0$ that
\begin{equation*}
    |\sum_{j=1}^i v_{j,j+1}(F_{j,j+1} - p_{j,j+1})|\leq C \sqrt{\sum_{j}\frac{v_{j,j+1}^2\log n}{L_{j,j+1}}}.
\end{equation*}
Using Chernoff's method we can show a slightly sharper bound that if $\min_j
L_{ij}> 25e^{2\kappa_E}\log n$ then with probability at least $1 - n^{-9}$,
\begin{equation*}
    \max_j |F_{ij} - p_{ij}| \leq C\sqrt{\frac{\log n}{L_{ij}v_{ij}}}.
\end{equation*}
Since $e^{\kappa_E}= \max_{(i,j)\in E}\frac{p_{ji}}{p_{ij}} \geq
\frac{1}{\frac{1}{p_{ij}p{ji}}p_{ij}^2} = \frac{1}{v_{ij}p_{ij}^2}$, the last
inequality and the condition on $L_{ij}$ imply that $|F_{ij} - p_{ij}|\leq
\min\{p_{ij},p_{ji}\}/5$ simultaneously, and consequently $z_{ji}\in
[0.8p_{ji},1.2p_{ji}]$. Therefore, the second term can be bounded as
\begin{equation*}
    \frac{1}{2}\frac{1 - 2z_{ji}}{z_{ji}^2(1 - z_{ji})^2}(F_{ij} - p_{ij})^2\leq c \frac{1}{p_{ij}^2p_{ji}^2}(F_{ij} - p_{ij})^2\leq c\frac{v_{ij}\log n}{L_{ij}}.
\end{equation*}
In particular, for the path graph, the summation of the second term is
controlled by $\sum_{j}\frac{v_{j,j+1}\log n}{L_{j,j+1}}$. Thus, when
\begin{equation*}
\min_j L_{j,j+1}> c e^{\kappa_E}n\log n,
\end{equation*}
for some sufficiently large constant $c$, the summation of the second term is
negligible compared to the bound on the summation of the first term
$\sqrt{\sum_{j}\frac{v_{j,j+1}^2\log n}{L_{j,j+1}}}$. Therefore, the error
\eqref{eq:err-closed_form} can be bounded as
\begin{equation*}
|\hat{\theta}_{i+1} - {\theta}^*_{i+1}|\leq C\sqrt{\sum_{j=1}^{i}\frac{\exp(2|\theta^*_j - \theta^*_{j+1}|)\log n}{L_{j,j+1}}},
\end{equation*}
which implies that $\|\hat{\theta} - \theta^*\|_{\infty}\leq
C\sqrt{\sum_{i=1}^{n-1}\frac{\exp(2|\theta^*_i - \theta^*_{i+1}|)\log
n}{L_{i,i+1}}}$.

In the special case $L_{i,i+1} = L$ for all $i$, we have $\|\hat{\theta} -
\theta^*\|_{\infty}\leq C\sqrt{\frac{\exp(2\kappa_E) n\log n}{L}}$ and
$\|\hat{\theta} - \theta^*\|_{2}\leq n\sqrt{\frac{\exp(2\kappa_E)\log n}{L}}$.
\eprfof

In Shah (2016) they prove the lower bound $e^{-\kappa} \frac{n}{\sqrt{L}}$ for
$\ell_2$ error under path graph, while in our paper we prove the lower bound
$e^{-\kappa} \sqrt{\frac{n}{L}}$ for $\ell_\infty$ error under path graph. Thus
the upper bound above on the closed-form MLE achieves the minimax lower bound up
to a $\sqrt{\log n}$ and $\exp(2\kappa_E)$ factor.

\textbf{Proposition 4} (General tree graph)\textbf{.} Suppose the graph is a
tree graph $([n],E)$ where item $i$ and $j$ are compared $L_{ij}$ times such
that $\min_{i,j}L_{ij}>c e^{2\kappa_E}n\log n$ for some universal constant $c$.
Then with probability at least $1 - n^{-4}$, the vanilla MLE $\hat{\bbrtheta}_0$
satisfies
\begin{equation}
\|\hat{\bbrtheta}_0 - \bbrtheta^*\|_{\infty}\lesssim \sqrt{\max_{i_1,i_2}\sum_{(i,j)\in {\rm path}(i_1,i_2)}\frac{\exp({2|\theta^*_{i} - \theta^*_{j}|})\log n}{L_{ij}}}.
\end{equation}
In particular, when $L_{i,j} = L$, we have
\begin{equation}
   \|\hat{\bbrtheta}_0 - \bbrtheta^*\|_{\infty}\lesssim e^{\kappa_E}\sqrt{\frac{D\log n}{L}}, \|\hat{\bbrtheta}_0 - \bbrtheta^*\|_{2}\lesssim e^{\kappa_E}\sqrt{\frac{Dn\log n}{L}}
\end{equation}
\bprfof{\Cref{prop:tree}} Suppose $\mclG = ([n],E)$ is a tree graph and let
$D$ be the diameter of the graph, i.e., $D = \max_{i,j\in [n]} |{\rm
path}(i,j)|$ where ${\rm
path}(i,j):=\{(i,i_1),(i_2,i_3),\cdots,(i_m,j)\}\subseteq E$ is the shortest
path from $i$ to $j$. For instance, for complete graph and star graph, $D = 1$;
for path graph, $D = n - 1$; for a complete binary tree with depth or height
$h$, $D = 2h$.

The key property of a tree graph is that it has three equivalent definitions: 1.
it is a maximal loop-free graph; 2. it is a minimal connected graph; 3. it is a
simple graph with $|E| = |V| - 1$. By the first definition, it can be seen that
for any two nodes $i$ and $j$, there is exactly one path from $i$ to $j$,
otherwise there will be a loop. Using this property, for a fixed item $i_0$, the
MLE equation
\begin{equation*}
\nabla \ell\left(\theta\right)_i=\sum_{j \in\mathcal{N}(i) } \left(\frac{M_{ij}}{M_{ij} +M_{ji}}-\frac{\exp(\theta_{i}-\theta_{j})}{1 + \exp(\theta_{i}-\theta_{j})}\right) = 0,\ \forall i\in [n]
\end{equation*}
can be solved by
\begin{equation*}
\hat{\theta}_{i_0} = 0,\quad \hat{\theta}_{l} = \sum_{(i,j)\in {\rm path}(i_0,l)}(\log R_{ji} - \log r_{ji}),\ l\neq i_0.
\end{equation*}
Note that, the choice of $i_0$ is not crucial for the $\ell_{\infty}$ bound, as
such shifting would at most lead to an inflation on the $\ell_{\infty}$ error by
2. In fact, if we let $\hat{\bbrepsilon}$ with $\hat{\epsilon}_{i_0} = 0$ be the
entry-wise error when we choose $\hat{\theta}_{i_0} = {\theta}_{i_0} = 0$, then
after shifting $\hat{\bbrtheta}$ to get $\tilde{\bbrtheta}$ with
$\tilde{\theta}_{i_1} = {\theta}_{i_1} = 0$ for some $i_1\neq i_0$ we will get
$\|\tilde{\bbrepsilon}\|_{\infty}\leq \|\hat{\bbrepsilon}\|_{\infty} +
\|\hat{\epsilon}_{i_1}\mathbf{1}_n\|_{\infty}\leq
2\|\hat{\bbrepsilon}\|_{\infty}$.

Following the same argument in the proof of \Cref{prop:path}, we can show
that when $\min_{(i,j)\in E} L_{i,j}> e^{2\kappa_E}n\log n$, we have
\begin{equation*}
\| \hat{\theta} - \theta^* \|_{\infty}\lesssim \sqrt{\max_{i_1,i_2}\sum_{(i,j)\in {\rm path}(i_1,i_2)}\frac{\exp({2|\theta^*_{i} - \theta^*_{j}|})\log n}{L_{ij}}}.
\end{equation*}
Again, the overall $\ell_{\infty}$ error is determined by the extreme values of
$|\theta_i^* - \theta_j^*|$ and $1/L_{ij}$.

To simplify the bound, we can consider the uniform sampling scheme under which
$L_{ij} = L$ for all $(i,j)\in E$. By the definition $D = \max_{i,j\in [n]}
|{\rm path}(i,j)|$, we have
\begin{equation}
\| \hat{\theta} - \theta^* \|_{\infty}\lesssim\sqrt{\frac{\exp(2\kappa_E) D \log n}{L}},
\end{equation}
and as a corollary,
\begin{equation}
\| \hat{\theta} - \theta^* \|_{2}\lesssim\sqrt{\frac{\exp(2\kappa_E) D n\log n}{L}}.
\end{equation}

In particular, for a path graph, $D = n-1$ and the two bounds become the same
bounds in \Cref{prop:path}.

For a star graph, $D=1$ and the bounds become
\begin{equation*}
\| \hat{\theta} - \theta^* \|_{\infty}\lesssim\sqrt{\frac{\exp(2\kappa_E)  \log n}{L}}, \quad \| \hat{\theta} - \theta^* \|_{2}\lesssim\sqrt{\frac{\exp(2\kappa_E) n\log n}{L}},
\end{equation*}
which are both sharp up to a $\log n$ factor according to the lower bounds in
\cite{shah2015estimationfrompairwisecomps} and our work. In addition, for the
star graph, one can actually show with much simpler arguments that the condition
on $L$ can be relaxed to $L> c e^{2\kappa_E}\log n$, because in a star graph
essentially we are just comparing $n-1$ pairs separately.
\eprfof


\clearpage


\subsection{Proof in Section 3}
\label{subsec:prf_lower_bounds}
Let $F(t) = \frac{1}{1 + e^{-t}}$, then the Bradley-Terry model can be written
as $p_{ij} = F[(w_i^* - w_j^*)/\sigma]$ with $\sigma = 1$. We have $\max_{t\in
[0,2\kappa/\sigma]} F'(t) = F'(0) = 1/4$, so $\zeta$ in
\Cref{nlem:kl-div-prob-score-vectors} satisfies
\begin{align*}
    \zeta = c e^{\frac{2\kappa}{\sigma}}.
\end{align*}
Moreover, we denote $\mclW:=\{\bbrtheta^*\mid \|{\bbrtheta}^*\|_\infty \leq
B\}$. Since $\mathbf{1}_{\btldima}^\top \bbrtheta^* = 0$, we have $B\asymp
\kappa:=\max_{i,j}|\theta^*_i - \theta^*_j|$ (in general, there is no more
information, but for some special cases, e.g., when entries of $\bbrtheta^*$ are
equal-spaced, we have $\kappa = 2B$). In what follows in this section, we still
use quantity $\sigma$ for generality, but keep in mind that in our setting
$\sigma=1$. We still use $B$ to differ it from $\kappa$.

\subsubsection{Background and Required
Results}\label{subsubsec:background-req-results} Before writing our proof of
\Cref{nthm:thm2-lb} we note down specific results from
\citep{shah2015estimationfrompairwisecomps} and other sources that we will use
repeatedly.

\bnumdefn[Pairwise $(\delta, \beta)-$packing set from
\citep{shah2015estimationfrompairwisecomps}]\label{ndefn:delta-beta-packing-set}
Let $\btldimc \in \reals^{\btldima}$ be a parameter to be estimated, as indexed
over a class of probability distributions $\mclP \defined
\thesetb{\mathbb{P}_{\btldimc}}{\btldimc \in \mclW}$, and let $\rho :
\reals^{\btldima} \times \reals^{\btldima} \to \reals_{\geq 0}$ be a
pseudo-metric. Suppose there exist a finite set of $M$ vectors
$\theseta{\btldimc^{1}, \ldots, \btldimc^{M}}$ such that the following
conditions hold:
\begin{align*}
\min_{j, k \in [M], j \neq k} \rho \parens{\btldimc^{j}, \btldimc^{k}} \geq \delta \text{ and } \frac{1}{\binom{M}{2}} \sum_{j, k \in[M], j \neq k} D_{\mathrm{KL}}\left(\mathbb{P}_{\btldimc^{j}} \| \mathbb{P}_{\btldimc^{k}}\right) \leq \beta
\end{align*}
Then we refer to $\theseta{\btldimc^{1}, \ldots, \btldimc^{M}}$ as $(\delta,
\beta)-$packing set.
\enumdefn

\bnlem[Fano minimax lower
bound]\label{nlem:pairwise-fano-minimax-lower-bound-mod} Suppose that we can
construct a $(\delta, \beta)-$packing set with cardinality $M$, then the minimax
risk is lower bounded as:

\begin{align*}
\inf_{\widehat{\btldimc}} \sup _{\btldimc^{*} \in \mathcal{W}} \mathbb{E}\brackets{\rho\parens{\widehat{\btldimc}, \btldimc^{*}}}
\geq \frac{\delta}{2}\left(1-\frac{\beta + \log 2}{\log M}\right)
\end{align*}
\enlem
\bprf\label{prf:pairwise-fano-minimax-lower-bound-mod} See
\citep[Lemma~3]{yu1997assouadfanolecamminimaxlb} for details.
\eprf

\bnlem[Equivalence of $\norm{\cdot}_{\infty}$ and $\norm{\cdot}_{2}$
norms]\label{nlem:equiv-ellinf-ell2-mod} Given any vector $\btldimc \in
\reals^{\btldima}$, with $\btldima \in \nats$ fixed, the following inequalities
holds:
\begin{align*}
    \frac{1}{\sqrt{\btldima}} \norm{\btldimc}_{2}
    \leq \norm{\btldimc}_{\infty}
    \leq \norm{\btldimc}_{2}
\end{align*}
\enlem

\bprf\label{prf:equiv-ellinf-ell2-mod} The result is standard \eg see
\citep[Proposition~2.10]{holger2018numericallinalg_book} for a more general
version and proof. For the sake of completeness we provide a direct proof of the
equivalent statement $\norm{\btldimc}_{\infty} \leq \norm{\btldimc}_{2} \leq
\sqrt{\btldima}\norm{\btldimc}_{\infty}$ as follows:
\begin{align*}
    \norm{\btldimc}_{\infty}
    &= \max_{i \in [\btldima]}{\abs{\btldimc_{i}}}
    = \sqrt{\max_{i \in [\btldima]}{\abs{\btldimc_{i}}^{2}}}
    \leq \sqrt{\sum_{i}^{\btldima} {\btldimc_{i}}^{2}}
    \leq \norm{\btldimc}_{2} \\
    \norm{\btldimc}_{2}
    &= \sqrt{\sum_{i}^{\btldima} {\btldimc_{i}}^{2}}
    \leq \sqrt{\sum_{i}^{\btldima} \max_{i \in [\btldima]}{\abs{\btldimc_{i}}^{2}}}
    = \sqrt{\btldima \max_{i \in [\btldima]}{\abs{\btldimc_{i}}^{2}}}
    = \sqrt{\btldima} \norm{\btldimc}_{\infty}
\end{align*}
This proves the lower and upper inequalities respectively, as required.
\eprf

\bnrmk\label{nrmk:equiv-ellinf-ell2-mod} We note that the above inequalities are
tight \ie have optimal constants. In the case of the upper bound consider
$\btldimc = \mathbf{1}_{\btldima}= \parens{1, \ldots, 1}^{\top} \in
\reals^{\btldima}$. So $\norm{\btldimc}_{\infty} = 1$ and $\norm{\btldimc}_{2} =
\sqrt{\btldima} = \sqrt{\btldima} \norm{\btldimc}_{\infty}$ showing the
tightness of the upper bound. In the case of the lower bound, consider $\btldimc
= \bfe_{1}$ \ie \WLOG the first standard basis vector in $\reals^{\btldima}$. We
then have that $\norm{\btldimc}_{\infty} = \norm{\btldimc}_{2} = 1$, which shows
tightness in the lower bound.
\enrmk

\bnlem[Lemma 7 in
\citep{shah2015estimationfrompairwisecomps}]\label{nlem:varshamov-gilbert-type-bnd-shah}
For any $\alpha \in\left(0, \frac{1}{4}\right),$ there exists a set of
$M(\alpha)$ binary vectors $\theseta{z^{1}, \ldots, z^{M(\alpha)}} \subset
\theseta{0,1}^{\btldima}$ such that

\begin{align}
\alpha \btldima
\leq \left\|z^{j}-z^{k}\right\|_{2}^{2}
&\leq \btldima \quad \text { for all } j \neq k \in [M(\alpha)], \text { and } \label{neqn:varshamov-gilbert-type-bnd-shah-1} \\
\angles{e_{1}, z^{j}}
&= 0 \quad \text { for all } j \in [M(\alpha)]
\label{neqn:varshamov-gilbert-type-bnd-shah-2}
\end{align}
\enlem

\bnlem[Lemma 8 in
\citep{shah2015estimationfrompairwisecomps}]\label{nlem:kl-div-prob-score-vectors}
For any pair of quality score vectors $\btldimc^{j}$ and $\btldimc^{k},$ and for
\begin{align*}
    \zeta &\defined \frac{\max_{x \in[0,2 B / \sigma]} F^{\prime}(x)}{F(2 B / \sigma)(1 - F(2 B / \sigma))}
\end{align*}
we have
\begin{align*}
    D_{\mathrm{KL}}\parens{\mathbb{P}_{\btldimc^{j}} \| \mathbb{P}_{\btldimc^{k}}}
    \leq \frac{\btldimb \zeta}{\sigma^{2}}\parens{\btldimc^{j}-\btldimc^{k}}^{\top} L\parens{\btldimc^{j}-\btldimc^{k}}
    \defines \frac{\btldimb \zeta}{\sigma^{2}}\norm{\btldimc^{j}-\btldimc^{k}}_{L}^{2}
\end{align*}
\enlem

\bnlem[Lemma 14 in
\citep{shah2015estimationfrompairwisecomps}]\label{nlem:laplacian-trace-constraints-shah}
The Laplacian matrix $\tilde{\mclL}_{\bfA}$ satisfies the trace constraints:
\begin{align}
    \tr{(\tilde{\mclL}_{\bfA})} &= 2 \\
    \tr{(\tilde{\mclL}_{\bfA}^{\dagger})} &\geq \frac{\btldima^{2}}{4}
\end{align}
\enlem
\bprf
See \citep[Lemma~14]{shah2015estimationfrompairwisecomps} for details.
\eprf

\nit The challenge - constructing a suitable pairwise packing set that meets
\Cref{ndefn:delta-beta-packing-set}. The main tool to use here is the
Varshamov-Gilbert Lemma

\subsubsection{Sketch of Lower Bound
Proof}\label{subsubsec:sketch-of-prf-lower-bound} In brief, we seek a minimax
lower bound proof of the same form as
\citep[Theorem~2(a)]{shah2015estimationfrompairwisecomps}, except in our case we
choose our packing set norm to the $\ell_{\infty}$ norm, rather than the
$\ell_{2}^{2}$ in \citep[Theorem~2(a)]{shah2015estimationfrompairwisecomps}. We
leverage their construction directly in two main ways. First, we use the
slightly modified version of Fano's Lemma that enables the $(\delta,
\beta)$-packing set to be constructed for the $\ell_{\infty}$ norm, not the
$\ell_{\infty}^{2}$ norm, consistent with our high probability upper bound per
\Cref{nlem:pairwise-fano-minimax-lower-bound-mod}. Second, we can switch out
their use of the $(\delta, \beta)$-packing set in the $\ell_{2}$ norm to the
$(\delta^{\prime}, \beta^{\prime})$-packing set in the $\ell_{\infty}$ norm.
This is done by using the topological equivalence of norms in finite dimensions
per \Cref{nlem:equiv-ellinf-ell2-mod}, which is shown to be tight in the
dimension per \Cref{nrmk:equiv-ellinf-ell2-mod}. For the sake of clarity, we use
much of the same wording as the proof from
\cite[Appendix~B]{shah2015estimationfrompairwisecomps}, for the convenience of
the reader.

\subsubsection{Lower Bound Proof - Part
I}\label{subsubsec:formal-prf-lower-bound-part-I}

Our proof follows directly the approach taken from
\citep[Section~B.1]{shah2015estimationfrompairwisecomps}. The normalized
Laplacian $\tilde{\mclL}_{\bfA}$ of the comparison graph is symmetric and
positive-semidefinite. We can thus decompose this via diagonalization as $L =
U^{\trp} \Lambda U$ where $U \in \reals^{\btldima \times \btldima}$ is an
orthonormal matrix, and $\Lambda$ is a diagonal matrix of nonnegative
eigenvalues $\Lambda_{jj} = \lambda_{j}(L)$ for each $j \in [\btldima]$. Similar
to \citep[Section~B.1]{shah2015estimationfrompairwisecomps} we first prove that
the minimax risk is lower bounded by $c \sigma^{2}
\frac{\btldima^{2}}{\btldimb}$.

Fix scalars $\alpha \in (0, \frac{1}{4})$ and $\delta > 0$, with values to be
specified later. Obtain set of vectors on the Boolean Hypercube $\theseta{0,
1}^{\btldima}$ \ie $\theseta{z^{1}, \ldots, z^{M(\alpha)}}$ given by
\Cref{nlem:varshamov-gilbert-type-bnd-shah}, where $M(\alpha)$ is set to be
\begin{align}\label{neqn:m-alpha-lower-bound-I}
    M(\alpha) &\defined \floor{\exp \braces{\frac{\btldima}{2}(\log 2+2 \alpha \log 2 \alpha+(1-2 \alpha) \log (1-2 \alpha))}}.
\end{align}
Define another set of vectors of the same cardinality $\thesetb{\btldimc^{j}}{j
\in [M(\alpha)]}$ via $\btldimc^{j} \defined \frac{\delta}{\sqrt{\btldima}}
U^{\top} P z^{j}$, where $P$ is a permutation matrix. The permutation matrix $P$
has the constraint that it keeps the first coordinate constant \ie $P_{11} = 1$.
By construction for each $j \neq k$ we have that

\begin{align}\label{neqn:formal-prf-lower-bound-delta-prime}
\norm{\btldimc^{j} - \btldimc^{k}}_{\infty}
\stackrel{(i)}{\geq} \frac{1}{\sqrt{\btldima}} \norm{\btldimc^{j} - \btldimc^{k}}_{2}
\stackrel{(ii)}{=} \frac{1}{\sqrt{\btldima}} \parens{\frac{\delta}{\sqrt{\btldima}}\norm{z^{j} - z^{k}}_{2}}
\stackrel{(iii)}{\geq} \delta\sqrt{\frac{\alpha}{\btldima}}
\end{align}

Here $(i)$ follows from \Cref{nlem:equiv-ellinf-ell2-mod}. Additionally $(ii)$
follows since $\btldimc^{j} \defined \frac{\delta}{\sqrt{\btldima}} U^{\top} P
z^{j}$. In the case of the final inequality $(iii)$, we have
$\frac{\delta^{2}}{\btldima^2}\norm{z^{j} - z^{k}}_{2}^{2} \geq \frac{\alpha
\btldima \delta^{2}}{\btldima^{2}} = \frac{\alpha \delta^{2}}{\btldima}$. Since
the set $\theseta{z^{1}, \ldots, z^{M(\alpha)}}$ are binary vectors with a
minimum Hamming distance at least $\alpha \btldima$ using
\Cref{neqn:varshamov-gilbert-type-bnd-shah-1}. Consider, any distinct $j, k \in
[M(\alpha)]$, then for some subset $\theseta{i_{1}, \ldots, i_{r}} \subseteq
\theseta{2, \ldots, \btldima}$ with $\alpha \btldima \leq r \leq \btldima$ it
must follow that

\begin{align*}
\norm{\btldimc^{j}-\btldimc^{k}}_{\tilde{\mclL}_{\bfA}}^{2}
= \frac{\delta^{2}}{\btldima}\norm{U^{\top} P z^{j}-U^{\top} P z^{k}}_{\tilde{\mclL}_{\bfA}}^{2}
= \frac{\delta^{2}}{\btldima}\norm{z^{j}-z^{k}}_{\Lambda}^{2}
= \frac{\delta^{2}}{\btldima} \sum_{m=1}^{r} \lambda_{i_{m}}(\tilde{\mclL}_{\bfA})
\end{align*}

The last part follows since $\Lambda$ is a diagonal matrix of non-negative
eigenvalues with $\Lambda_{ii} = \lambda_{j}(L)$. Now for given $\theseta{a_{2},
\ldots, a_{\btldima}}$ such that $\alpha \btldima \leq \sum_{i = 2}^{\btldima}
a_{i} \leq \btldima$ we have that
\begin{align*}
    \frac{1}{\binom{M(\alpha)}{2}} \sum_{j \neq k}\norm{\btldimc^{j}-\btldimc^{k}}_{\tilde{\mclL}_{\bfA}}^{2}
    &= \frac{\delta^{2}}{\btldima} \sum_{i=2}^{\btldima} a_{i} \lambda_{i}(\tilde{\mclL}_{\bfA})
\end{align*}
The permutation matrix $P$ is chosen such that the last $\btldima - 1$
coordinates are permuted to have $a_{1} \geq \ldots \geq a_{\btldima}$ and keep
the $\btldima^{\text{th}}$ coordinate fixed. By this particular choice, and
using the fact that $\tr{(\tilde{\mclL}_{\bfA})} = 2$ we have that:
\begin{align*}
    \frac{1}{\binom{M(\alpha)}{2}} \sum_{j \neq k}\norm{\btldimc^{j}-\btldimc^{k}}_{\tilde{\mclL}_{\bfA}}^{2}
    = \frac{\delta^{2}}{\btldima} \frac{\btldima}{\btldima - 1} \tr{(\tilde{\mclL}_{\bfA})}
    \leq \frac{2 \delta^{2}}{\btldima} \tr{(\tilde{\mclL}_{\bfA})}
    = \frac{4 \delta^{2}}{\btldima}
\end{align*}
Now by the choice of $P$ above, we have that for every choice of $j \in
\brackets{M(\alpha)}$
\begin{align*}
    \angles{\tilde{\mclL}_{\bfA}, \btldimc^{j}} = \frac{\delta}{\sqrt{\btldima}} \bfe_{1}^{\top} P z^{j}=\bfe_{1}^{\top} z^{j}=0
\end{align*}
where the last equality follows from
\Cref{neqn:varshamov-gilbert-type-bnd-shah-2}. Now the basic condition needs to
be verified \ie did the $\btldimc^{j}$ we chose satisfy the boundedness
constraint, to ensure that $\btldimc^{j} \in \mclW_{B}$? Setting $\delta^{2} =
0.01 \frac{\sigma^{2} \btldima^{2}}{4 \btldimb \zeta}$, it indeed follows that
$\norm{\btldimc^{j}}_{\infty} \leq
\frac{\delta}{\sqrt{\btldima}}\norm{z^{j}}_{2} \stackrel{(i)}{\leq} \delta
\stackrel{(ii)}{\leq} B$. Here $(i)$ follows since $z^{j} \in \theseta{0,
1}^{\btldima}$. Furthermore $(ii)$ follows from our choice of $\delta$ and our
assumption that $\btldimb \geq \frac{c \sigma^{2}
\operatorname{tr}\parens{\tilde{\mclL}_{\bfA}^{\dagger}}}{\zeta B^{2}}$ with $c
= 0.002$, where \Cref{nlem:laplacian-trace-constraints-shah} guarantees that
$\btldimb \geq \frac{c \sigma^{2} \btldima^{2}}{4 \zeta B^{2}}$. We have thus
verified that each vector $\btldimc_{j}$ also satisfies the boundedness
constraint $\norm{\btldimc^{j}}_{\infty} \leq B$, which is required for
membership in $\mclW_{B}$. Finally by \Cref{nlem:kl-div-prob-score-vectors} we
have that:
\begin{align}\label{neqn:formal-prf-lower-bound-beta-prime}
    D_{\mathrm{KL}}\left(\mathbf{P}_{\btldimc^{j}} \| \mathbf{P}_{\btldimc^{k}}\right)
    &\leq \frac{\btldimb \zeta}{\sigma^{2}} \frac{4 \delta^{2}}{\btldima}=0.01 \btldima
\end{align}

To summarize, we have now constructed a $(\delta^{\prime},
\beta^{\prime})$-packing set with respect to the norm $\rho
\parens{\btldimc^{j}, \btldimc^{k}} \defined \norm{\btldimc^{j} -
\btldimc^{k}}_{\infty}$, where $\delta^{\prime} =
\delta\sqrt{\frac{\alpha}{\btldima}}$ from
\Cref{neqn:formal-prf-lower-bound-delta-prime}, and $\beta^{\prime} = 0.01d$
from \Cref{neqn:formal-prf-lower-bound-beta-prime}.

Finally we have by substituting $(\delta^{\prime}, \beta^{\prime})$ into the
pairwise Fano's lower bound (\Cref{nlem:pairwise-fano-minimax-lower-bound-mod})
that:
\begin{align*}
    \sup _{\btldimc^{*} \in \mathcal{W}_{B}} \mathbb{E}\left[\left\|\widetilde{\btldimc}-\btldimc^{*}\right\|_{\infty}\right] \geq \delta\sqrt{\frac{\alpha}{\btldima}}\left(1-\frac{0.01 \btldima+\log 2}{\log M(\alpha)}\right) = c\sigma\sqrt{\frac{\btldima}{\zeta \btldimb}}\left(1-\frac{0.01 \btldima+\log 2}{\log M(\alpha)}\right)
\end{align*}
which yields the claim, after appropriate substitution of $\delta$ and setting
$\alpha = 0.01$.

For the case of $\btldima \leq 9$, consider the set of the three
$\btldima$-length vectors $z^{1} = \parens{0, \ldots, -1}$, $z^{2} = \parens{0,
\ldots, 0}$ and $z^{3} = \parens{0, \ldots, 0}$. Construct the packing set
$\parens{\btldimc^{1}, \btldimc^{2}, \btldimc^{3}}$ from these three vectors
$\parens{z^{1}, z^{2}, z^{3}}$ as done above for the case of $\btldima>9$. From
the calculations made for the general case above, we have for all pairs $\min
_{j \neq k} \norm{\btldimc^{j}-\btldimc^{k}}_{\infty}^{2}\geq \frac{1}{9}\min
_{j \neq k} \norm{\btldimc^{j}-\btldimc^{k}}_{2}^{2} \geq \frac{\delta^{2}}{81}$
and $\max _{j, k}\norm{\btldimc^{j}-\btldimc^{k}}_{\tilde{\mclL}_{\bfA}}^{2}
\leq 4 \delta^{2}$, and as a result $\max_{j, k}
D_{\mathrm{KL}}\left(\mathbf{P}_{\btldimc^{j}} \|
\mathbf{P}_{\btldimc^{k}}\right) \leq \frac{4 \btldimb \zeta
\delta^{2}}{\sigma^{2}}$. Choosing $\delta^{2}=\frac{\sigma^{2} \log 2}{8
\btldimb \zeta}$ and applying the pairwise Fano's lower bound
(\Cref{nlem:pairwise-fano-minimax-lower-bound-mod}) yields the claim.

For the other case, the lower bound in terms of
$\lambda_i(\tilde{\mclL}_{\bfA})$, the argument is similar, and we end up with
an extra factor of $\frac{1}{\btldima^{\prime}}$.

\subsubsection{Lower Bound Proof - Part
II}\label{subsubsec:formal-prf-lower-bound-part-II}

Given an integer $\btldima^{\prime} \in \theseta{2, \ldots, \btldima}$, and
constants $\alpha \in (0, \frac{1}{4})$, $\delta > 0$, define the integer:

\begin{align}\label{neqn:m-alpha-lower-bound}
    M^{\prime}(\alpha)
    &\defined \floor{\exp \braces{\frac{\btldima^{\prime}}{2}(\log 2+2 \alpha \log 2 \alpha+(1-2 \alpha) \log (1-2 \alpha))}}
\end{align}

Applying \Cref{nlem:varshamov-gilbert-type-bnd-shah} using $\btldima^{\prime}$
as the dimension results in a subset $\theseta{z^{1}, \ldots,
z^{M^{\prime}(\alpha)}}$ of the Boolean hypercube $\theseta{0,
1}^{\btldima^{\prime}}$, with specified properties. We then define a finite set
of size $M^{\prime}(\alpha)$, of $\btldima$-length vectors
$\theseta{\widetilde{\btldimc}^{1}, \ldots,
\widetilde{\btldimc}^{M^{\prime}(\alpha)}}$ using:

\begin{align*}
    \widetilde{\btldimc}^{j}
    &= \brackets{0\left(z^{j}\right)^{\top} 0 \cdots 0}^{\top} \quad \text { for each } j \in [M(\alpha)]
\end{align*}

For each $j \in \brackets{M(\alpha)}$, let us define $\btldimc^{j} \defined
\frac{\delta}{\sqrt{\btldima^{\prime}}} U^{\top} \sqrt{\Lambda^{\dagger}}
\tilde{\btldimc}^{j}$. For the first standard basis vector $\bfe_{1} \in
\reals^{\btldima}$, we then have that $\angles{\mathbf{1}_{\btldima},
\btldimc^{j}} = \frac{\delta}{\sqrt{\btldima^{\prime}}}
\tilde{\mclL}_{\bfA}^{\top} U^{\top} \sqrt{\Lambda^{\dagger}}
\widetilde{\btldimc}^{j} = 0$. Here the main fact used is
$\tilde{\mclL}_{\bfA}\mathbf{1}_{\btldima} = 0$. Additionally we have that for
any $j \neq k$, we have that:

\begin{align}
\norm{\btldimc^{j}-\btldimc^{k}}_{\infty}^{2}\geq \frac{1}{\btldima}\norm{\btldimc^{j}-\btldimc^{k}}_{2}^{2}
= \frac{1}{\btldima}\frac{\delta^{2}}{\btldima^{\prime}}\left(\widetilde{\btldimc}^{j}-\widetilde{\btldimc}^{k}\right)^{\top} \Lambda^{\dagger}\left(\widetilde{\btldimc}^{j}-\widetilde{\btldimc}^{k}\right)
\geq \frac{1}{\btldima}\frac{\delta^{2}}{\btldima^{\prime}} \sum_{i=\left\lceil(1-\alpha) \btldima^{\prime}\right\rceil}^{\btldima^{\prime}} \frac{1}{\lambda_{i}}
\end{align}

Now, setting $\delta^{2}=0.01 \frac{\sigma^{2} \btldima^{\prime}}{\btldimb
\zeta}$ results in:

\begin{align}
\norm{\btldimc^{j}}_{\infty}
\leq \frac{\delta}{\sqrt{\btldima^{\prime}}}\norm{\sqrt{\Lambda^{\dagger}} \tilde{\btldimc}^{j}}_{2}
\stackrel{(i)}{\leq} \frac{\delta}{\sqrt{\btldima^{\prime}}} \sqrt{\operatorname{tr}\left(\Lambda^{\dagger}\right)}
\stackrel{(ii)}{=} \frac{\delta}{\sqrt{\btldima^{\prime}}} \sqrt{\operatorname{tr}\left(L^{\dagger}\right)}
\stackrel{(iii)}{\leq} B
\end{align}

where inequality $(i)$ follows from the fact that $z^{j}$ has entries in
$\{0,1\} ;$ step $(ii)$ follows because the matrices $\sqrt{\Lambda^{\dagger}}$
and $\sqrt{\tilde{\mclL}_{\bfA}^{\dagger}}$ have the same eigenvalues; and
inequality $(iii)$ follows from our choice of $\delta$ and our assumption
$\btldimb \geq \frac{c \sigma^{2}
\operatorname{tr}\left(\tilde{\mclL}_{\bfA}^{\dagger}\right)}{\zeta B^{2}}$ on
the sample size with $c = 0.01$. We have thus verified that each vector
$\btldimc^{j}$ also satisfies the boundedness constraint
$\norm{\btldimc^{j}}_{\infty} \leq B$, as required for membership in
$\mclW_{B}$. Furthermore, for any pair of distinct vectors in this set, we have:

\begin{align*}
    \norm{\btldimc^{j}-\btldimc^{k}}_{\tilde{\mclL}_{\bfA}}^{2}
    = \frac{\delta^{2}}{\btldima^{\prime}}\norm{z^{j}-z^{k}}_{2}^{2}
    \leq \delta^{2}
\end{align*}
By \Cref{nlem:kl-div-prob-score-vectors} we have that:
\begin{align}
    D_{\mathrm{KL}}\left(\mathbf{P}_{\btldimc^{j}} \| \mathbf{P}_{\btldimc^{k}}\right)
    &\leq \frac{\btldimb \zeta}{\sigma^{2}} \norm{\btldimc^{j} - \btldimc^{k}}_{\tilde{\mclL}_{\bfA}}^{2} = 0.01 \btldima^{\prime}
\end{align}
Finally we have by substituting $(\delta^{\prime}, \beta^{\prime})$ into the
pairwise Fano's lower bound (\Cref{nlem:pairwise-fano-minimax-lower-bound-mod})
that:
\begin{align*}
    \sup _{\btldimc^{*} \in \mathcal{W}_{B}} \mathbb{E}\left[\left\|\widetilde{\btldimc}-\btldimc^{*}\right\|_{\infty}\right]
    \geq \frac{\frac{\delta}{\sqrt{\btldima^{\prime} d}} \sqrt{\sum_{i=\ceil{(1-\alpha) \btldima^{\prime}}}^{\btldima^{\prime}} \frac{1}{\lambda_{i}}}}{2}
         \parens{1-\frac{0.01 \btldima^{\prime}+\log 2}{\log M^{\prime}(\alpha)}}
\end{align*}

Substituting our choice of $\delta$ and setting $\alpha = 0.01$ proves the claim
for $\btldima^{\prime} > 9$.

For the case of $\btldima^{\prime} \leq 9$. Consider the packing set of the
three $\btldima^{\prime}$-length vectors $\btldimc^{1} = \delta U
\sqrt{\Lambda^{\dagger}}\parens{0, 1, \ldots, 0}$, $\btldimc^{2} =
-\btldimc^{1}$ and $\btldimc^{3} = \parens{0, \ldots, 0}$. Then we have for all
pairs $\min _{j \neq k} \norm{\btldimc^{j}-\btldimc^{k}}_{\infty}^{2} \geq
\frac{1}{9}\min _{j \neq k} \norm{\btldimc^{j}-\btldimc^{k}}_{2}^{2} \geq
\frac{\delta^{2}}{9\lambda_{2}(L)}$ and $\max _{j,
k}\norm{\btldimc^{j}-\btldimc^{k}}_{\tilde{\mclL}_{\bfA}}^{2} \leq 4
\delta^{2}$, and as a result $\max_{j, k}
D_{\mathrm{KL}}\left(\mathbf{P}_{\btldimc^{j}} \|
\mathbf{P}_{\btldimc^{k}}\right) \leq \frac{4 \btldimb \zeta
\delta^{2}}{\sigma^{2}}$. Choosing $\delta^{2}=\frac{\sigma^{2} \log 2}{8
\btldimb \zeta}$ and applying the pairwise Fano's lower bound
(\Cref{nlem:pairwise-fano-minimax-lower-bound-mod}) yields the claim.


\clearpage


\subsection{Additional Experiments} \label{sec:additional-experiments} In this
section, we show some additional results of experiments. The section has three
parts:
\begin{enumerate}
    \item Experiments related to \Cref{sec:implications-of-work}.
    \item Extra comparisons as a supplement to \Cref{sec:simulations}.
    \item Experiments to illustrate that in some cases $\kappa_E$ is still loose
    and point out a potential future direction.
\end{enumerate}

\subsubsection{Subadditivity}
\label{subsec:additivity}
In this section, we show some simulation results illustrating the advantage of
using subadditivity property in the estimation of the \btl{} model, as we
discuss in \Cref{sec:implications-of-work}.

\paragraph{Island graph.} In this setting, we consider the Island graph with $n$
nodes described in \Cref{exa:n_ij=0} and denote the node set of the $k$-th
island as $V_k$. Suppose we get the estimator $\hat{\bbrtheta}^{(k)}\in
\mathbb{R}^{|V_k|}$ ($k>1$) based on $k$-th Island $V_k$ with the augmented
version $\tilde{\bbrtheta}^{(k)}\in \mathbb{R}^n$ such that the subvector
$\tilde{\bbrtheta}^{(k)}(V_k) = \hat{\bbrtheta}^{(k)}$. We can define the
ensemble estimator add-MLE in the following way: first shift
$\tilde{\bbrtheta}^{(k)}$ and get
\begin{equation*}
    \check{\bbrtheta}^{(k)} = \hat{\bbrtheta}^{(k)} + s_k,
    s_k = s_{k-1} + \tilde{\bbrtheta}^{(k - 1)}(n_{\rm island}) - \tilde{\bbrtheta}^{(k)}(n_{\rm overlap}), s_0 = 0.
\end{equation*}
Then construct $\tilde{\bbrtheta}^{add}\in \mathbb{R}^n$ such that for all $k$,
$\tilde{\bbrtheta}^{add}(V_k) = \check{\bbrtheta}^{(k)}$. At last, we centerize
$\check{\bbrtheta}^{(k)}$ and get
\begin{equation*}
    \hat{\bbrtheta}^{add} = \tilde{\bbrtheta}^{add} -\mathbf{1}_n\cdot \frac{1}{n}\mathbf{1}_n^\top \tilde{\bbrtheta}^{add}.
\end{equation*}

For the ease of implementation and precise description of the performance, we
construct the true parameter $\bbrtheta^*$ by first set $\bbrtheta^*(i) =
\bbrtheta^*(1) + (i - 1)\delta$ such that ${\rm avg}(\bbrtheta^*) = 0$, and then
shift $\bbrtheta_{(k)}^*:=\bbrtheta^*(V_k)$ by $s_k := -(k-1)s$ and call $s\in
\mathbb{R}$ the shifting coefficient. \Cref{fig:island_additivity} shows that
add-MLE outperforms the joint-MLE, where ``shift'' in the right panel is the
shifting coefficient $s$. Notice that to save space and for the ease of
understanding, we transform ``shift'' to ``diff'' in the
\Cref{sec:simulations} so that ``diff'' shows the difference in the
average: ${\rm avg}(\bbrtheta_{(k-1)}^*)-{\rm avg}(\bbrtheta_{(k)}^*)$. We set
$L =10$ for all pairs.

\begin{figure}[!ht]
\centering
\includegraphics[width=0.45\textwidth]{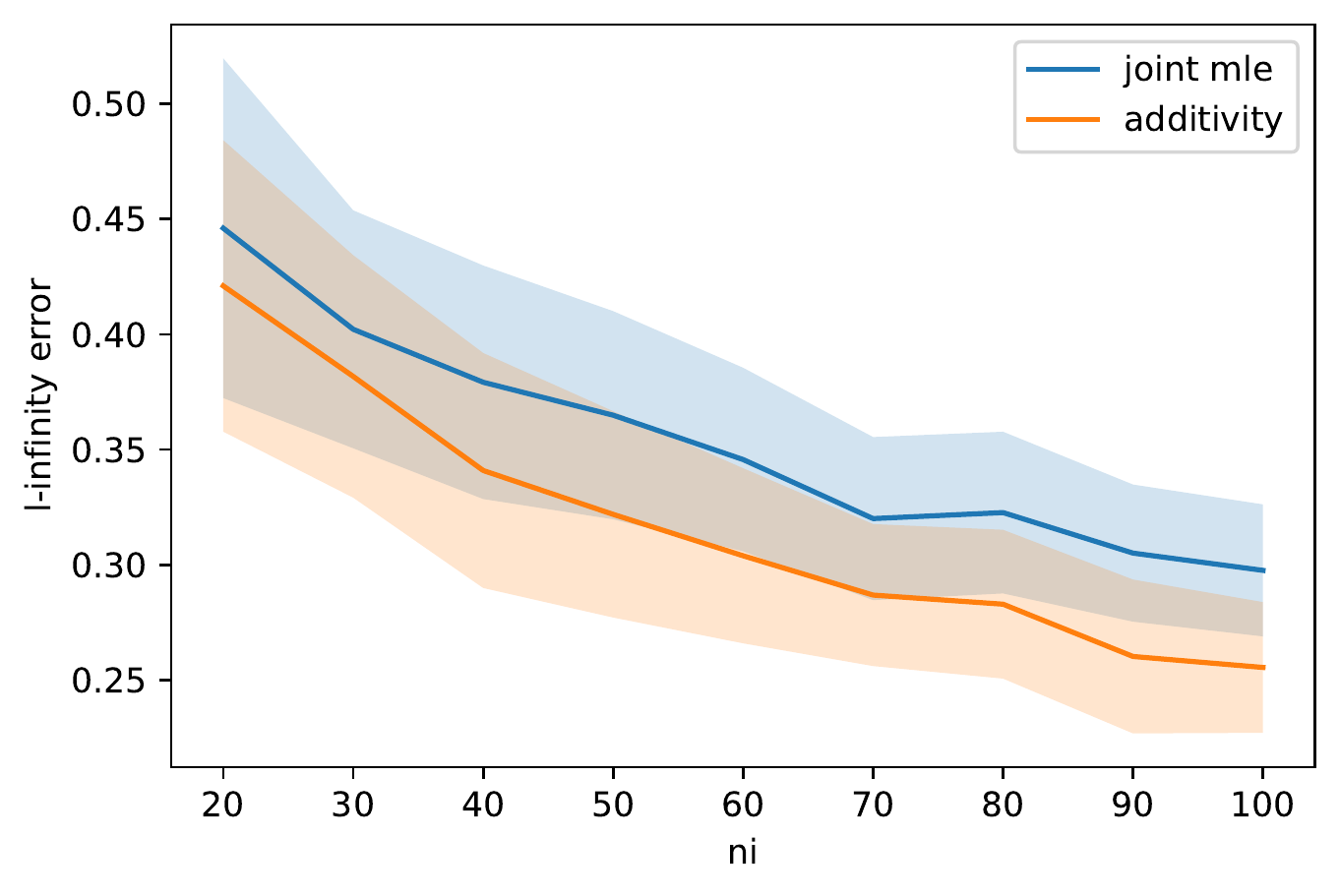}
\includegraphics[width=0.45\textwidth]{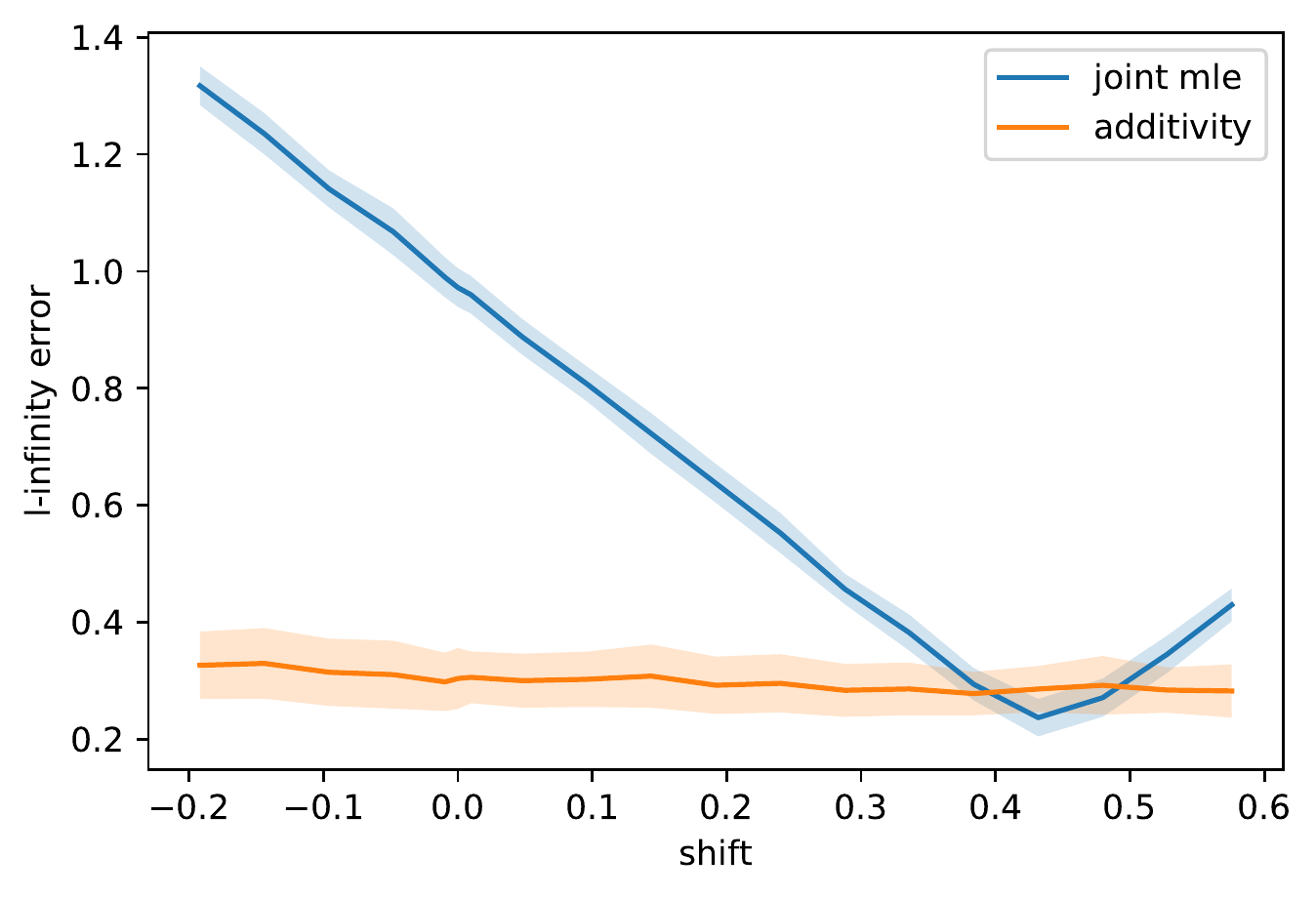}
\caption{Estimation errors given by joint-MLE and add-MLE on Island graphs. The
y-axis is $\|\hat{\bbrtheta} - \bbrtheta^*\|_{\infty}$. In the left panel, the
x-axis is the size of islands $n_{island}$ while $n_{overlap} = 5$, $k = 3$,
$L=10$ in all cases. In the right panel, the x-axis is $s$, the shifting
coefficient while $n_{island} = 50$, $n_{overlap} = 5$, $k = 5$, $L = 10$. The
lines show the average of 100 trials with the standard deviation shown as the
shaded area.}
\label{fig:island_additivity}
\end{figure}

\paragraph{Barbell graph.} In this setting, a barbell graph
$\mclG_{barbell}([n],E)$ of two equal-sized complete subgraphs $\mclG_1,
\mclG_2$ linked by a set of bridge edges $E_{bridge} = \{(i,j):i\in \mclG_1,
j\in \mclG_2, (i,j)\in E\}$ are generated as is discussed in
\Cref{exa:bridge}. Note that the vertex sets $V_1, V_2$ of
$\mclG_1,\mclG_2$ are disjoint and $V_1\cup V_2 = [n]$. We set $L =10$ for all
pairs. Again, two methods of estimation are compared:
\begin{itemize}
    \item Joint-MLE. A single regularized MLE $\hat{\bbrtheta}^{joint}$ is
    fitted on $\mclG_{barbell}$.
    \item Add-MLE. Two regularized MLE's $\hat{\bbrtheta}^{(1)}$ and
    $\hat{\bbrtheta}^{(2)}$ are fitted separately on $\mclG_1$ and $\mclG_2$.
    Then for each $e = (i,j)\in E_{bridge}$, we calculate
    \begin{equation*}
        \hat{d}_e := \log(\frac{\hat{p}_e}{1 - \hat{p}_e}),
    \end{equation*}
    where $\hat{p}_{e} = {\rm clip}(\frac{win_{ij}}{win_{ij} + loss_{ij}},
    p_{up}, p_{lb})$ for $e = (i,j)$ with two constants $0<p_{lb}<p_{up}<1$ for
    regularity. We take $p_{lb} = 0.1, p_{up} = 0.9$. Then for $e = (i,j)\in
    E_{bridge}$, the shifting constant $\hat{s}_e$ is defined as
    \begin{equation*}
        \hat{s}_e := \hat{d}_e - (\tilde{\theta}^{(1)}_i - \tilde{\theta}^{(2)}_j),
    \end{equation*}
    where $\tilde{\bbrtheta}^{(i)}$ is the augmented version of
    $\hat{\bbrtheta}^{(1)}$ satisfying $\tilde{\bbrtheta}^{(i)}(V_i) =
    \hat{\bbrtheta}^{(1)}$. Then the average $\hat{s}_E$ is calculated via
    $\hat{s}_E := \frac{1}{|E|}\sum_{e\in E}\hat{s}_e$ and the add-MLE
    $\hat{\bbrtheta}^{add}$ is constructed via
        \begin{equation*}
            \tilde{\theta}^{add}_i :=\begin{cases}
            \tilde{\theta}^{(1)}_i,i\in\mclG_1,\\
            \tilde{\theta}^{(2)}_i - \hat{s}_E,i\in\mclG_2,
            \end{cases}
        \end{equation*}
    and $\hat{\bbrtheta}^{add} = \tilde{\bbrtheta}^{add} -\mathbf{1}_n\cdot
    \frac{1}{n}\mathbf{1}_n^\top \tilde{\bbrtheta}^{add}$.
\end{itemize}
Notice that we slightly change the way of constructing the add-MLE in
\Cref{sec:implications-of-work} to exploit multiple random bridge edges in
this setting.

Similar to the Island graph case, we let $\mclG_1$ be the complete graph of node
$\{1,\ldots, \frac{n}{2}\}$ and $\mclG_2$ the complete graph of nodes
$\{\frac{n}{2} + 1,\ldots, n\}$. To set the true parameter, we let
\begin{equation*}
    \theta^*_i = \begin{cases}
    \theta_1^* + (i - 1)\delta,i\leq \frac{n}{2},\\
    -s + \theta_1^* + (i - 1)\delta,i> \frac{n}{2},
    \end{cases}
\end{equation*}
where $s\in \mathbb{R}$ is the shifting coefficient.

As we can see from \Cref{fig:barbell_additivity}, the  performance of the
add-MLE is more stable than that of the joint-MLE, while ensuring better or
similar $\ell_{\infty}$ estimation error.
\begin{figure}[!ht]
\centering
\includegraphics[width=0.45\textwidth]{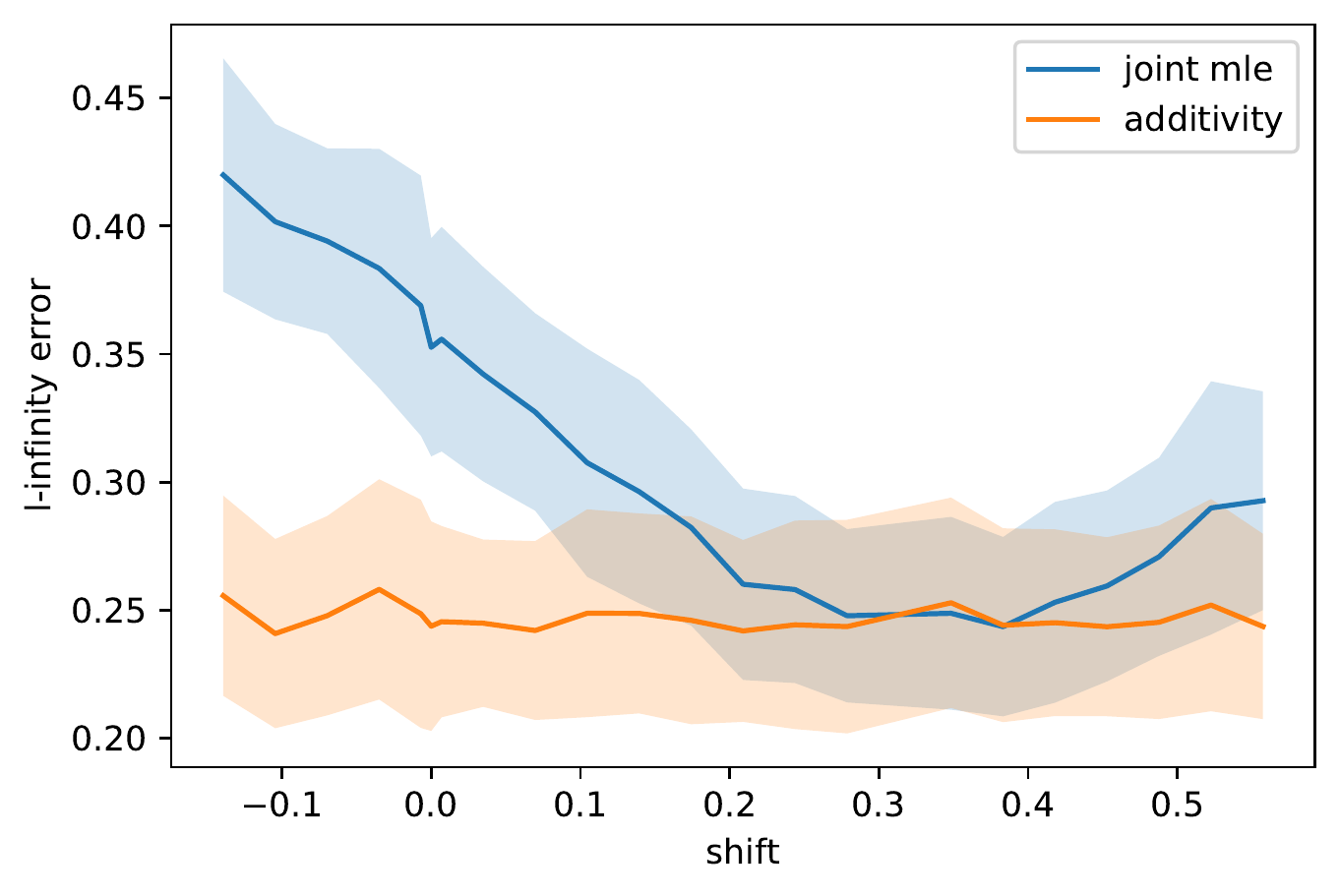}
\includegraphics[width=0.45\textwidth]{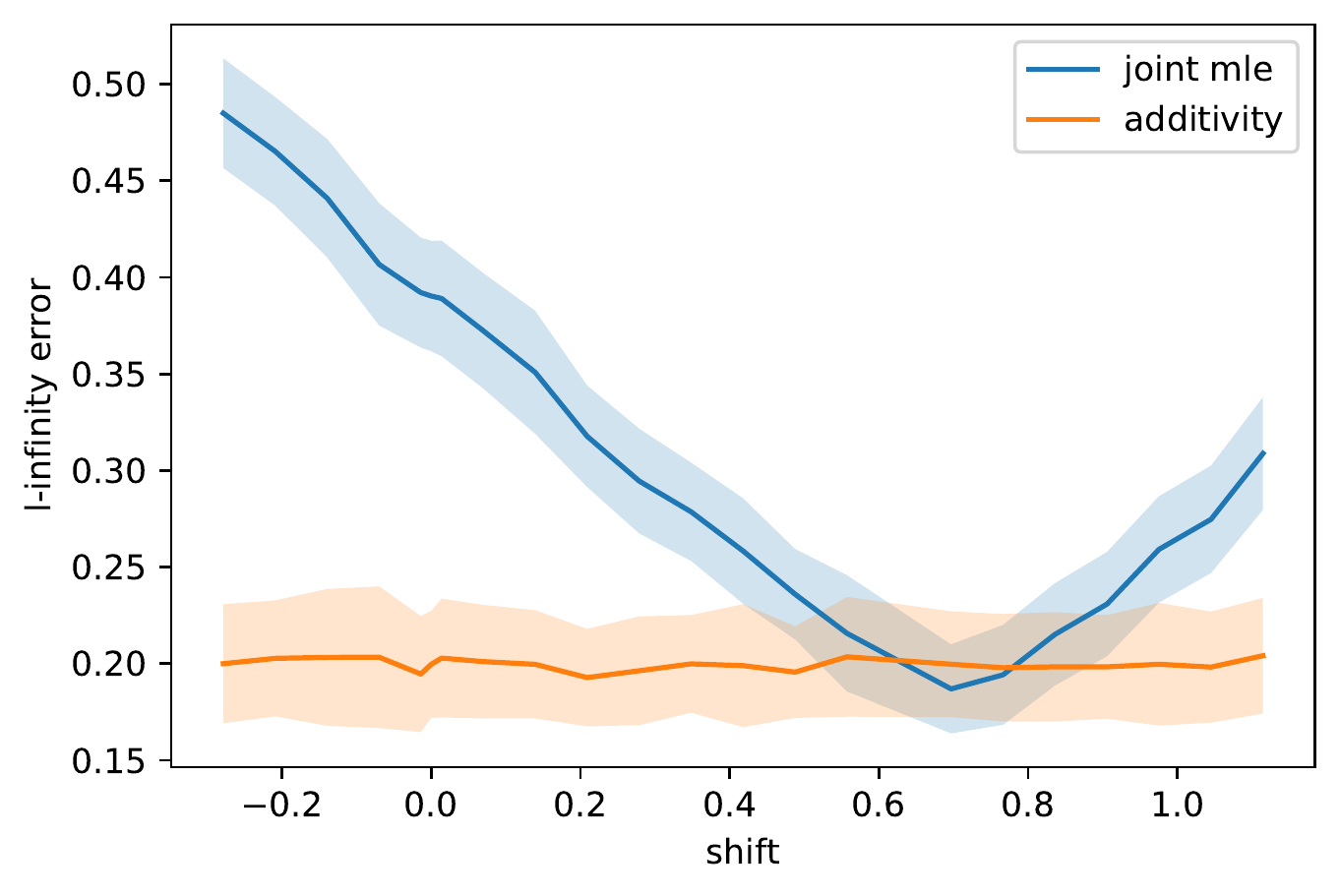}
\caption{Estimation errors given by joint-MLE and add-MLE on Barbell graphs with
100 nodes(left) and 200 nodes (right). Two equal-sized complete subgraphs are
linked by 10 (left) and 20 (right) randomly sampled bridge edges. The y-axis is
$\|\hat{\bbrtheta} - \bbrtheta^*\|_{\infty}$ and the x-axis is $s$, the shifting
coefficient. $L = 100$ for bridge edges and $L = 10$ for non-bridge edges. The
lines show the average of 100 trials with the standard deviation shown as the
shaded area.}
\label{fig:barbell_additivity}
\end{figure}

\subsubsection{Extra comparisons}

A one-sentence summary of this section is, we demonstrate our discussions in
\Cref{sec:simulations} by real cases that $\kappa_E$ can give much tighter
upper bounds than $\kappa$.

We consider a $k$-banded graph where comparisons are made only for pairs with
difference in indices smaller or equal to $k$. That is, the edge set of the
comparison graph is $E_{k}:= \{(i,j): |i-j|\leq k\}$. We consider two settings,
$k = \sqrt{n}$ and $k = n/\log n$, and in each setting, we set $\kappa = \log
(n)$, $L = 10$ and $\theta^*_i = \theta_1^* + (i-1)\delta$ for $i>1$ with
$\delta = \kappa / (n-1)$.

We compare the real $\ell_2$-error of the regularized MLE and three upper
bounds: the upper bound for the $\ell_2$-error provided by our paper using
$\kappa$ and $\kappa_E$, and the one provided by
\cite{shah2015estimationfrompairwisecomps}. As is shown in
\Cref{fig:banded_compare}, $\kappa_E$ gives tighter upper bounds than $\kappa$
in the setting of banded comparison graph. It should be noted that all curves
are going up because for a banded graph, $L$ needs to be sufficiently large to
guarantee $o(1)$ $\ell_{\infty}$ error, and we set $L = 10$ simply for
illustration of the effectiveness of $\kappa_E$ versus $\kappa$ here.

\begin{figure}[!ht]
\centering
\includegraphics[width=0.45\textwidth]{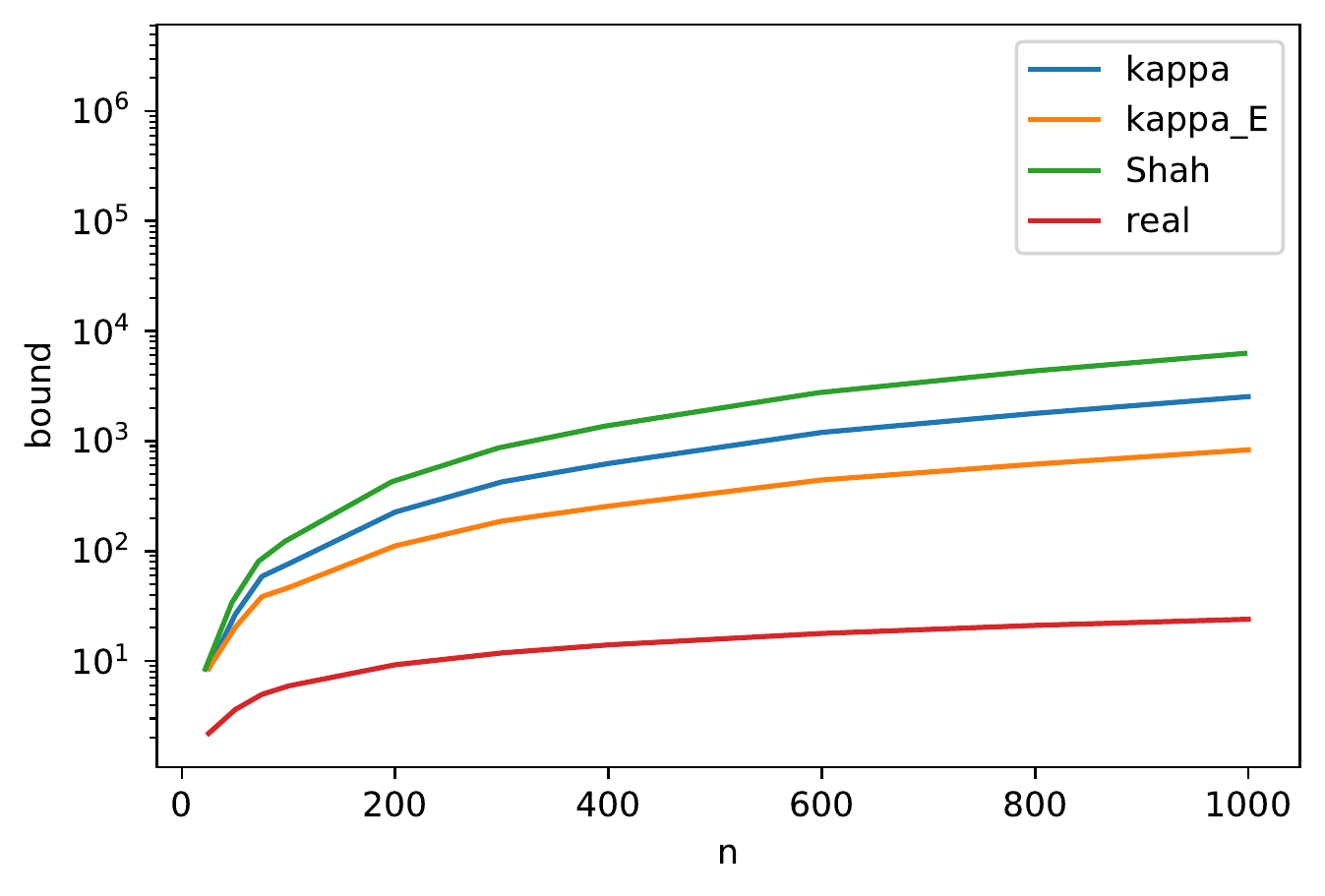}
\includegraphics[width=0.45\textwidth]{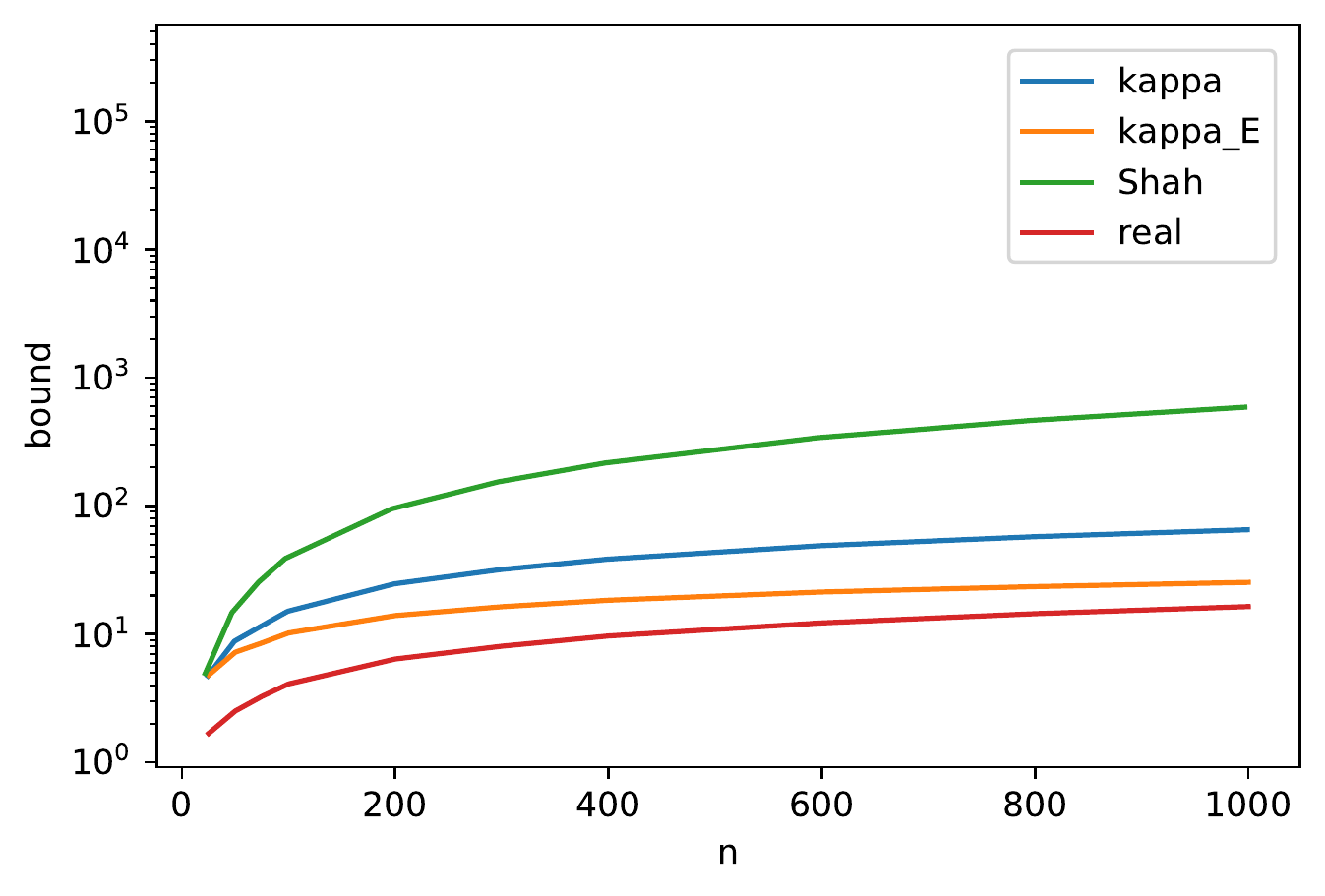}
\caption{Comparison under the $k$-banded graph with $k = \sqrt{n}$ (left) and $k
= n/\log(n)$ (right). Each point on lines is an average of 20 trials. The two
curves of upper bounds \texttt{kappa} and \texttt{Shah} are manually shifted
downwards so that they started at the same level with the curve for
\texttt{kappa\_E}, to remove the affect of the choice of constant (though all
leading constants are set to be 1 here) and make it easier to compare the
increasing rate of the bound.}
\label{fig:banded_compare}
\end{figure}

\subsubsection{Cases where $\kappa_E$ is loose}
Consider a path graph of $n$ nodes and edge set $\{(i,i+1):i\in [n-1]\}$. Assume
$\theta_i^*= \theta_1^* + (i - 1)\delta$, then $\kappa_E = \delta$ and $\kappa =
(n-1)\delta$. In this case, a factor of $e^{\kappa_E}$ gives tighter control
than $e^{\kappa}$. However, if we add one edge $(1,n)$ into the graph,
$\kappa_E$ becomes $(n-1)\delta$ and our upper bound will increase a lot, which
is counter-intuitive because the newly-added 1 out of $n$ edges should not
affect the estimation accuracy too much. In other words, the bound gets looser
after the new edge is added.

\begin{figure}[!ht]
\centering
\includegraphics[width=0.45\textwidth]{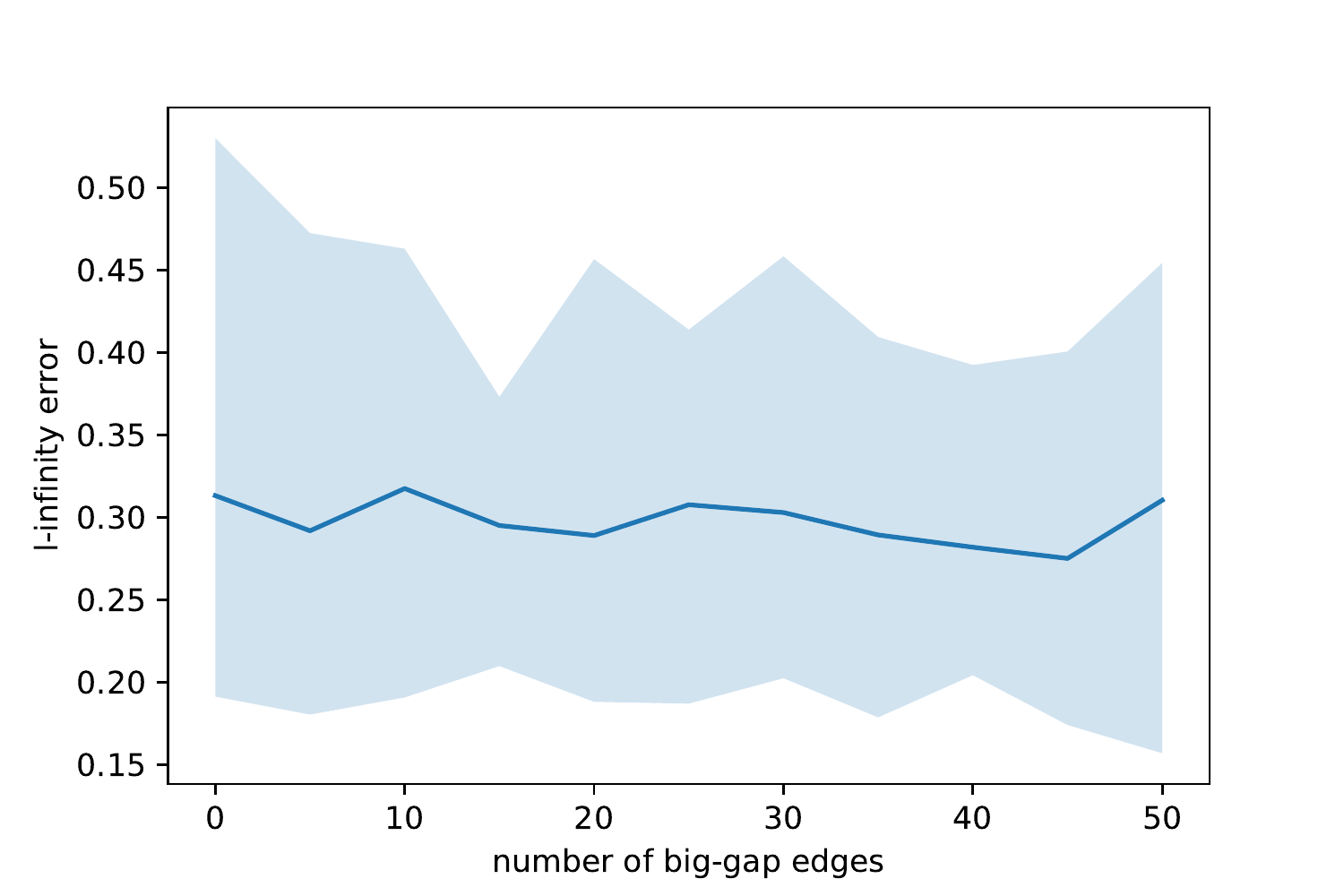}
\includegraphics[width=0.45\textwidth]{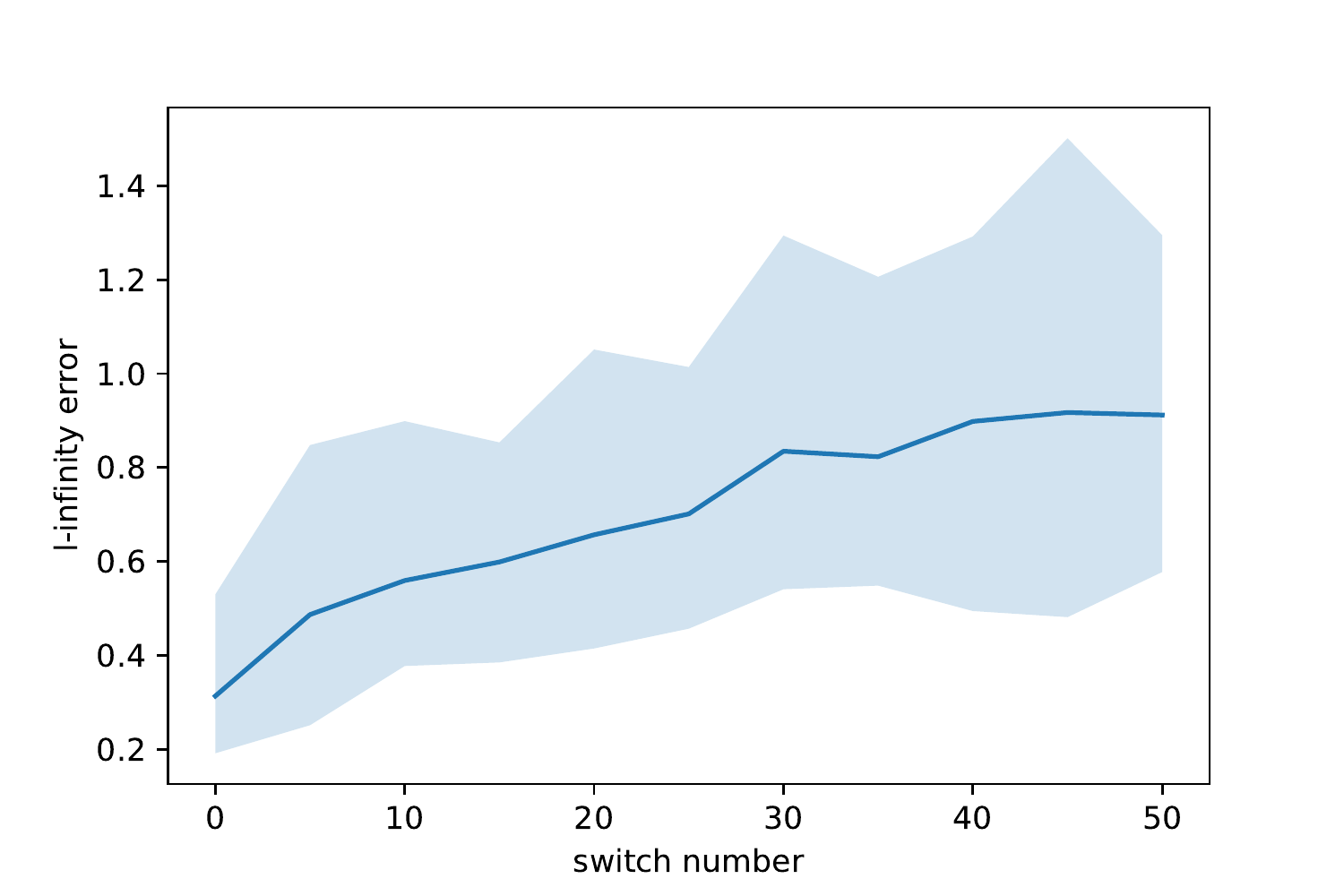}
\caption{Experiment on the path graph. Left: the $\ell_{\infty}$ error of the
regularized MLE when adding new edges with big performance gaps to the path
graph. Right: the $\ell_{\infty}$ error of the regularized MLE when switch some
pairs of the performance parameter $\{\theta_i^*\}$ so that there are more
big-gap edges while keeping the algebraic connectivity. The results show that
when the proportion of big-gap edges is small, the estimation error would not be
affected a lot. Both experiments are under $\kappa\approx 6.9$, equal-gap
$\theta^*$, $n = 200$ and $L = 5000$ and based on 40 trials for each
hyperparameter, with 0.05 and 0.95 quantiles shown by the shaded area.}
\label{fig:path_graph_big_gap}
\end{figure}

The reason is that $\kappa_E$ itself is not enough to provide tight control
across all comparison graphs. For instance, to avoid such loose cases, one also
needs to take into account the proportion of edges with big performance gaps
compared to edges with small performance gaps.

\begin{figure}[!ht]
\centering
\includegraphics[width=0.6\textwidth]{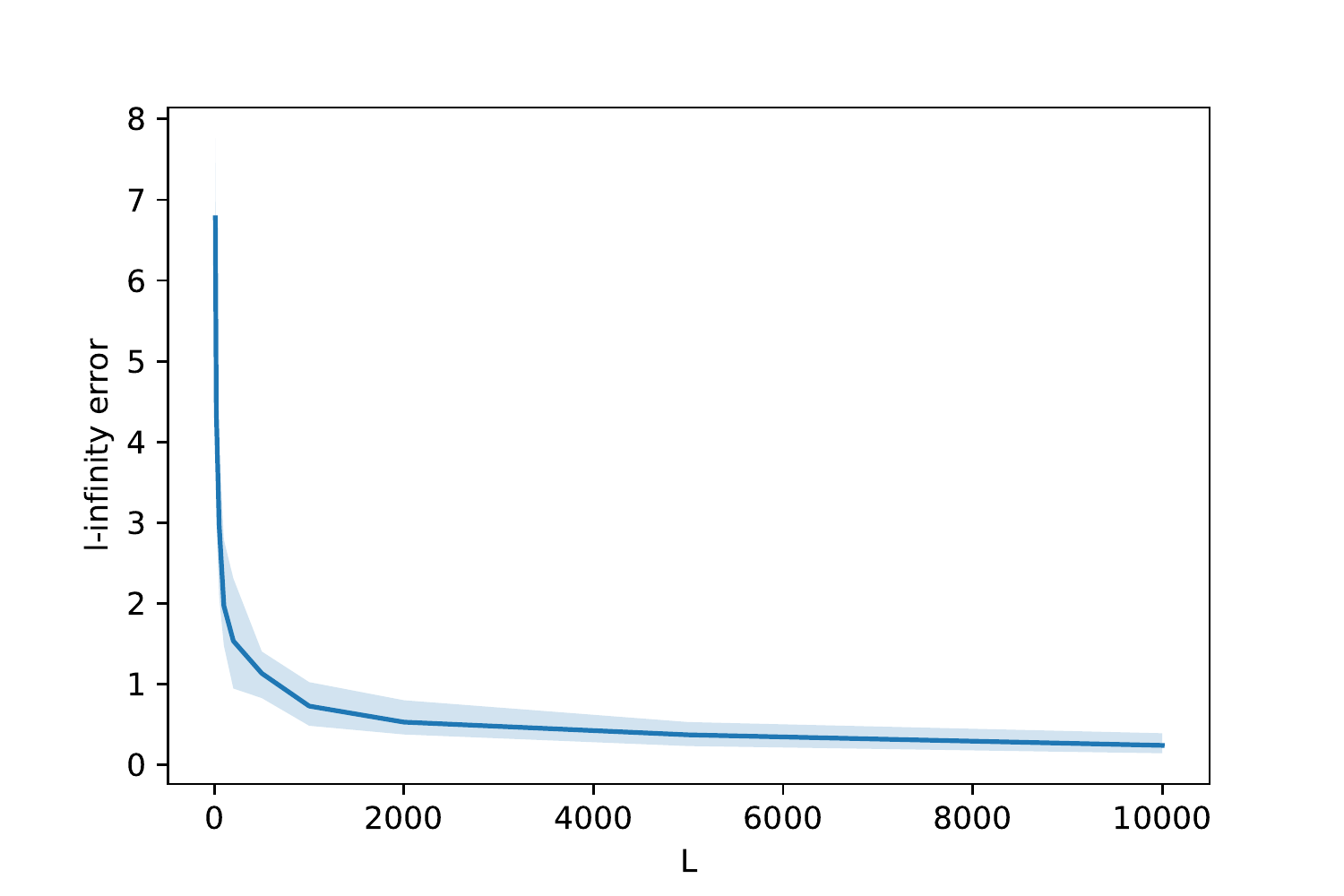}
\caption{Experiment on the path graph, under $\kappa\approx 6.9$, equal-gap $\theta^*$, $n = 200$ and based on 40 trials for each hyperparameter $L$, with 0.05 and 0.95 quantiles shown by the shaded area.}
\label{fig:path_graph}
\end{figure}

\Cref{fig:path_graph_big_gap} shows some numerical results on the impact of
big-gap edges. In the left panel, we add edges to the top-right and bottom-left
corner of the adjacency matrix so that their performance gaps are big. In the
right panel, we switch some pairs of parameters $\{\theta_i^*\}$ so that there
are more big-gap edges while keeping the algebraic connectivity. The first
switch will switch $\theta_1^*,\theta^*_{[n/2] + 1}$. The $i$-th switch will
switch the pair $\theta^*_{2i - 1}, \theta^*_{[{n}/{2}] + 2i - 1}$. As an
example, taking $n=8$, $\theta_1^* = 0$, and $\delta = 1$, 2 switches will make
parameters change as
\[
(1,2,3,4,5,6,7,8) \longrightarrow (5, 2, 7, 4, 1, 6, 3, 8).
\]
In the left panel, as new big-gap edges come in, the algebraic connectivity of
the graph gets larger as well, keeping estimation errors in the same level. In
the right panel, as we keep the algebraic connectivity constant, the impact of
the proportion of big-gap edges is shown more clearly.

Another thing we need to point out is that for cases like path graphs, even if
we ignore the factor of $\kappa$, bounding estimation error itself is hard due
to the poor connectivity of the graph, as is argued in
\cite{shah2015estimationfrompairwisecomps} (their upper bound for both path and
cycle graphs is not optimal). As we have shown in \Cref{sec:special-case},
even for $d$-regular graphs with relatively small $d$, the algebraic
connectivity is not big enough and our upper bound cannot match the lower bound.
But as the comparison graph gets denser and more regular, $\kappa_E$ will get
closer to $\kappa$, making the difference between $e^{\kappa_E}$ and
$e^{\kappa}$ not as dramatic, although in finite sample phase the difference can
still be big because it is an exponential factor.

The last thing we want is that the estimation itself (without providing a
theoretical tight upper bound for the estimation error) under the path graph is
hard: one need a huge $L$ to make accurate estimation when $n$ is big, as is
shown in \Cref{fig:path_graph}.


\clearpage


\subsection{Other Supporting Results} \label{sec:others}

\textbf{Proposition 6.} (Subadditivity) Let $I_1,I_2,I_3$ be three subsets of
$[n]$ such that $\cup_{j=1}^3 I_j = [n]$ and, for each $j \neq k$, $I_j
\not\subseteq I_k$ and for $i=1,2$, $I_i \cap I_3 \neq \emptyset$. Assume that
the sub-graphs induced by the $I_j$'s are connected. Let $\bbrtheta^*$ be the
vector of preference scores in the BTL model over $n$ items and
$\hat{\bbrtheta}_{(j)}$ be the MLE of $\bbrtheta^*_{(j)}$ for the BTL model
involving only items in $I_j$, $j=1,2,3$, with augmented versions
$\tilde{\bbrtheta}_{(j)}\in \mathbb{R}^n$ such that
$\tilde{\bbrtheta}_{(j)}(I_j) = \hat{\bbrtheta}_{(j)}$. Take two nodes $t_1\in
I_1\cap I_3$, $t_2\in I_2\cap I_3$, and let $\delta_3 =
\tilde{\bbrtheta}_{(1)}(t_1) - \tilde{\bbrtheta}_{(3)}(t_1)$, $\delta_2 =
\tilde{\bbrtheta}_{(3)}(t_2) - \tilde{\bbrtheta}_{(2)}(t_2)$. An ensemble MLE
$\hat{\bbrtheta}\in \mathbb{R}^n$ is a vector such that $\hat{\bbrtheta}(I_1) =
\hat{\bbrtheta}_{(1)}$, $\hat{\bbrtheta}(S_2) = \tilde{\bbrtheta}_{(2)}(S_2) +
\delta_3 + \delta_2$, and $\hat{\bbrtheta}(S_3) = \tilde{\bbrtheta}_{(3)}(S_3) +
\delta_3$, where $S_2 = I_2\setminus I_1$ and $S_3 = I_3\setminus (I_1\cup
I_2)$. It holds for any ensemble MLE $\hat{\bbrtheta}$ that
\begin{equation}
    \frac{1}{4}d_{\infty}(\hat{\bbrtheta},\bbrtheta^*) \leq d_{\infty}(\hat{\bbrtheta}_{(1)},\bbrtheta^*_{(1)}) + d_{\infty}(\hat{\bbrtheta}_{(2)},\bbrtheta^*_{(2)}) +
d_{\infty}(\hat{\bbrtheta}_{(3)},\bbrtheta^*_{(3)}),
\end{equation}
where $d_{\infty}(\bfv_1,\bfv_2) \defined \|(\bfv_1 - \mathbf{1}^\top{\rm
avg}(\bfv_1))-(\bfv_2 - \mathbf{1}^\top{\rm avg}(\bfv_2))\|_{\infty}$, where
${\rm avg}(\bfx):=\frac{1}{n}\textbf{1}_n^\top \bfx$.
\bprfof{\Cref{nlem:subbadditivity-ellinf-norm}} First, for each $i =
1,2,3$, we may shift $\hat{\bbrtheta}_{(i)}$ and  ${\bbrtheta}^*_{(i)}$
separately by constant vectors to ensure that so that
\begin{equation*}
    {\rm avg}(\hat{\bbrtheta}_{(i)}) = \frac{1}{|I_i|}\mathbf{1}_{|I_i|}^\top \hat{\bbrtheta}_{(i)} = 0\quad \text{and} \quad {\rm avg}({\bbrtheta}^*_{(i)})=0,
\end{equation*}
so that
\begin{equation}\label{eq:d_infty}
    d_{\infty}(\hat{\bbrtheta}_{(i)}, {\bbrtheta}^*_{(i)}) = \| \hat{\bbrtheta}_{(i)} - {\bbrtheta}^*_{(i)} \|_{\infty}.
\end{equation}
Next, for each $i = 1,2,3$, let $\tilde{\bbrtheta}_{(i)}\in \mathbb{R}^n$ be the
augmented version of $\hat{\bbrtheta}_{(i)}\in \mathbb{R}^{|I_i|}$, given by
\begin{equation*}
    \tilde{\bbrtheta}_{(i)}(I_i) =
    \hat{\bbrtheta}_{(i)}
    ,\quad  \tilde{\bbrtheta}_{(i)}(j) = 0,\ j\notin I_{i},
\end{equation*}
where $\bfv(j)$ refers to the $j$-th entry of vector $\bfv$. Similarly, we
define
\begin{equation*}
    \tilde{\bbrtheta}^*_{(i)}(I_i) =
    {\bbrtheta}^*_{(i)}
    ,\quad  \tilde{\bbrtheta}^*_{(i)}(j) = 0,\ j\notin I_{i},
\end{equation*}
\textbf{Step 1}. We now define the ensemble MLE $\hat{\bbrtheta}$ and show the
subadditivity property. The idea is to first fix $\tilde{\bbrtheta}_{(1)}$, and
then shift the entries of $\tilde{\bbrtheta}_{(2)}$ and
$\tilde{\bbrtheta}_{(3)}$ with coordinates in $I_2$ and $I_3$ respectively to
comply with the differences in the common entries of $I_1,I_3$ and $I_2,I_3$.

Let $S_1 = I_1$, $S_2 = I_2\setminus I_1$ (note that we don't put any constraint
on $I_1\cap I_2$), $S_3= I_3\setminus (I_1\cup I_2)$. We allow $S_3$ to be
$\emptyset$, but by the assumption that $I_j\not\subseteq I_k$, we have $S_2\neq
\emptyset$. Since $\cup_{j=1}^3 I_j = [n]$ and $I_i \cap I_3 \neq \emptyset$ for
$i=1,2$, we have $\cup_{j=1}^3 S_j = [n]$ and $S_i\cap S_k = \emptyset$ for any
$i\neq k$. Pick an arbitrary $t_1\in I_1\cap I_3$, $t_2\in I_2\cap I_3$, and let
\begin{equation*}
    \delta_3 = \tilde{\bbrtheta}_{(1)}(t_1) - \tilde{\bbrtheta}_{(3)}(t_1),\ \delta_2 = \tilde{\bbrtheta}_{(3)}(t_2) - \tilde{\bbrtheta}_{(2)}(t_2).
\end{equation*}
Notice that since there is no constraint on $I_1\cap I_2$, $t_1$ can be equal to
$t_2$.
Moreover, define $\check{\bbrtheta}_{(3)}$ and  $\check{\bbrtheta}_{(2)}$ by
\begin{equation*}
    \check{\bbrtheta}_{(3)}(j) = \begin{cases}
    \tilde{\bbrtheta}_{(3)}(j) + \delta_3, \quad j\in S_{3},\\
    0,\quad j\notin S_{3},
    \end{cases} \quad \text{and} \quad
    \check{\bbrtheta}_{(2)}(j) = \begin{cases}
    \tilde{\bbrtheta}_{(2)}(j) + \delta_3 + \delta_2, \quad j\in S_{2},\\
    0,\quad j\notin S_{2},
    \end{cases}
\end{equation*}
 respectively. Letting
\begin{equation}
    \hat{\bbrtheta} = \tilde{\bbrtheta}_{(1)} + \check{\bbrtheta}_{(2)} + \check{\bbrtheta}_{(3)},
\end{equation}
 it can be seen that
\begin{equation*}
    \hat{\bbrtheta}(j) = \begin{cases}
    \tilde{\bbrtheta}_{(1)}(j), \quad j\in S_{1}, \\
    \tilde{\bbrtheta}_{(2)}(j) + \delta_3 + \delta_2, \quad j\in S_{2},\\
    \tilde{\bbrtheta}_{(3)}(j) + \delta_3, \quad j\in S_{3}.\\
    \end{cases}
\end{equation*}
\textbf{Step 2}. Let $\delta_3^* = \tilde{\bbrtheta}^*_{(1)}(t_1) -
\tilde{\bbrtheta}^*_{(3)}(t_1)$, $\delta_2^* = \tilde{\bbrtheta}^*_{(3)}(t_2) -
\tilde{\bbrtheta}^*_{(2)}(t_2)$. Define a new true parameter
$\check{\bbrtheta}^*$ by shifting $\bbrtheta^*$ via
\begin{equation*}
    \check{\bbrtheta}^*(j) = \begin{cases}
    \tilde{\bbrtheta}^*_{(1)}(j), \quad j\in S_{1}, \\
    \tilde{\bbrtheta}^*_{(2)}(j) + \delta_3^* + \delta_2^*, \quad j\in S_{2},\\
    \tilde{\bbrtheta}^*_{(3)}(j) + \delta_3^*, \quad j\in S_{3}.\\
    \end{cases}
\end{equation*}
It can be verified that $\check{\bbrtheta}^* = {\bbrtheta}^* +
(\tilde{\bbrtheta}^*(t_1) - {\bbrtheta}^*(t_1))\mathbf{1}_n$, i.e.,
$\check{\bbrtheta}^*$ is a shift of ${\bbrtheta}^*$.
To show this, it suffices to show that for all $j\in [n]$,
\begin{equation}
    \check{\bbrtheta}^*(j) - \bbrtheta^*(j) = \tilde{\bbrtheta}^*(t_1) - \bbrtheta^*(t_1).
    \label{tmp_eq:shifting}
\end{equation}
In fact, for $j\in S_1$, since $\tilde{\bbrtheta}_{(1)}^*(I_1)$ is a shift of
$\bbrtheta^*({I_1})$, \Cref{tmp_eq:shifting} holds immediately by the definition
of $ \check{\bbrtheta}^*$. For $j\in S_3$, we have
\begin{equation*}
     \check{\bbrtheta}^*(j) - \bbrtheta^*(j) = \tilde{\bbrtheta}_{(3)}^*(j) + \tilde{\bbrtheta}^*_{(1)}(t_1) - \tilde{\bbrtheta}^*_{(3)}(t_1) - \bbrtheta^*(j).
\end{equation*}
Since $t_1\in I_3$ and $\tilde{\bbrtheta}_{(1)}^*(I_3)$ is a shift of
$\bbrtheta^*({I_3})$, it holds that $\tilde{\bbrtheta}_{(3)}^*(j) -
\bbrtheta^*(j) = \tilde{\bbrtheta}_{(3)}^*(t_1) - \bbrtheta^*(t_1)$ and
\Cref{tmp_eq:shifting} follows. For $j\in S_2$, we have
\begin{align*}
        \check{\bbrtheta}^*(j) - \bbrtheta^*(j) &= \tilde{\bbrtheta}_{(2)}^*(j) + \tilde{\bbrtheta}^*_{(1)}(t_1) - \tilde{\bbrtheta}^*_{(3)}(t_1) + \tilde{\bbrtheta}^*_{(3)}(t_2) - \tilde{\bbrtheta}^*_{(2)}(t_2) - \bbrtheta^*(j)\\
        & = \tilde{\bbrtheta}_{(2)}^*(j) - \bbrtheta^*(j) +  \tilde{\bbrtheta}^*_{(3)}(t_2) - \tilde{\bbrtheta}^*_{(2)}(t_2) + \tilde{\bbrtheta}^*_{(1)}(t_1) - \tilde{\bbrtheta}^*_{(3)}(t_1).
\end{align*}
Since $t_2\in I_2$ and $\tilde{\bbrtheta}_{(2)}^*(I_2)$ is a shift of
$\bbrtheta^*({I_2})$, it holds that $\tilde{\bbrtheta}_{(2)}^*(j) -
\bbrtheta^*(j) = \tilde{\bbrtheta}_{(2)}^*(t_2) - \bbrtheta^*(t_2)$ and thus
\begin{equation*}
    \check{\bbrtheta}^*(j) - \bbrtheta^*(j) = \tilde{\bbrtheta}^*_{(3)}(t_2) - {\bbrtheta}^*(t_2) + \tilde{\bbrtheta}^*_{(1)}(t_1) - \tilde{\bbrtheta}^*_{(3)}(t_1)
\end{equation*}
Again, since $t_1,t_2\in I_3$, we have $\tilde{\bbrtheta}^*_{(3)}(t_2) -
{\bbrtheta}^*(t_2) = \tilde{\bbrtheta}^*_{(3)}(t_1) - {\bbrtheta}^*(t_1)$ and
\Cref{tmp_eq:shifting} follows.

\textbf{Step 3}. Now we are ready to show the conclusion of the proposition. To
analyze the error, we first notice that for any $\bfv,\bfu\in \mathbb{R}^n$ and
$a\in \mathbb{R}$,
\begin{equation*}
    d_{\infty}(\bfu,\bfv) = d_{\infty}(\bfu,\bfv + a\mathbf{1}_n)\leq \|\bfu - (\bfv + a\mathbf{1}_n)\|_\infty + \frac{1}{n}\sum_{i=1}^n|u_i - (v_i + a)|\leq 2\|\bfu - (\bfv + a\mathbf{1}_n)\|_\infty.
\end{equation*}
Since $\check{\bbrtheta}^*$ is a shift of ${\bbrtheta}^*$, we have
\begin{equation*}
    d_{\infty}(\hat{\bbrtheta},{\bbrtheta}^*) \leq 2\|\hat{\bbrtheta} - \check{\bbrtheta}^*\|_{\infty}.
\end{equation*}
For $j\in S_1$,
\begin{equation*}
    |\hat{\bbrtheta}(j) - \check{\bbrtheta}^*(j)| = |\tilde{\bbrtheta}_{(1)}(j) - \tilde{\bbrtheta}^*_{(1)}(j)|\leq \|\hat{\bbrtheta}_{(1)} - {\bbrtheta}^*_{(1)}\|_{\infty}.
\end{equation*}
For $j\in S_3$,
\begin{equation*}
    |\hat{\bbrtheta}(j) - \check{\bbrtheta}^*(j)| = |\tilde{\bbrtheta}_{(3)}(j) - \tilde{\bbrtheta}^*_{(3)}(j)| + |\delta_3 - \delta_3^*|\leq 2\|\hat{\bbrtheta}_{(1)} - {\bbrtheta}^*_{(1)}\|_{\infty}+\|\hat{\bbrtheta}_{(3)} - {\bbrtheta}^*_{(3)}\|_{\infty}.
\end{equation*}
For $j\in S_2$.
\begin{align*}
    |\hat{\bbrtheta}(j) - \check{\bbrtheta}^*(j)| &= |\tilde{\bbrtheta}_{(2)}(j) - \tilde{\bbrtheta}^*_{(2)}(j)| + |\delta_3 - \delta_3^*| + |\delta_2 - \delta_2^*|\\
    &\leq \|\hat{\bbrtheta}_{(1)} - {\bbrtheta}^*_{(1)}\|_{\infty} + 2\|\hat{\bbrtheta}_{(2)} - {\bbrtheta}^*_{(2)}\|_{\infty}+2\|\hat{\bbrtheta}_{(3)} - {\bbrtheta}^*_{(3)}\|_{\infty}.
\end{align*}
Therefore by definition of $\|\cdot\|_\infty$ and \Cref{eq:d_infty},
\begin{align*}
    \|\hat{\bbrtheta} - \check{\bbrtheta}^*\|_{\infty}&\leq 2(\|\hat{\bbrtheta}_{(1)} - {\bbrtheta}^*_{(1)}\|_{\infty} + \|\hat{\bbrtheta}_{(2)} - {\bbrtheta}^*_{(2)}\|_{\infty}+\|\hat{\bbrtheta}_{(3)} - {\bbrtheta}^*_{(3)}\|_{\infty})\\
    &=2(d_{\infty}(\hat{\bbrtheta}_{(1)},\bbrtheta^*_{(1)}) + d_{\infty}(\hat{\bbrtheta}_{(2)},\bbrtheta^*_{(2)}) +
d_{\infty}(\hat{\bbrtheta}_{(3)},\bbrtheta^*_{(3)})),
\end{align*}
and we have
\begin{equation*}
    \frac{1}{4}d_{\infty}(\hat{\bbrtheta},{\bbrtheta}^*)\leq d_{\infty}(\hat{\bbrtheta}_{(1)},\bbrtheta^*_{(1)}) + d_{\infty}(\hat{\bbrtheta}_{(2)},\bbrtheta^*_{(2)}) +
d_{\infty}(\hat{\bbrtheta}_{(3)},\bbrtheta^*_{(3)}).
\end{equation*}
\eprfof

\bnlem\label{lem:cayley} Given a $d$-Cayley graph, where $(i,j)\in E$ if and
only if $i - j\equiv k({\rm mod}\ n)$ with $-d\leq k\leq d$, $k\neq 0$. It
satisfies $\lambda_2(\mclL_{\bfA}) \asymp d^3/n^2$.
\enlem
\begin{proof}
By definition, a $d$-Cayley graph is a $2d$-regular graph, and it's well known
that \citep[see, e.g.][]{spectra2012} the spectra of its adjacency matrix $\bfA$
is given by
\begin{equation*}
    \lambda_j(\bfA) = \sum_{k=1}^d (\zeta_j^k+\zeta_{-j}^k),\quad \zeta_j \defined \cos\frac{2\pi j}{n} + \sqrt{-1}\sin\frac{2\pi j}{n},
\end{equation*}
for $j = 1,\cdots,n$. Thus, $\lambda_2(\mclL_{\bfA}) = 2d -
2\sum_{k=1}^d\cos(\frac{2\pi k}{n})$. Since $\cos(\frac{2\pi k}{n}) = 1 -
2\sin^2 \frac{\pi}{n}k$ and for $k\leq d< 0.5n$, $\sin \frac{\pi}{n}k\in
(0.5\frac{\pi}{n}k,\frac{\pi}{n}k)$, we have
\begin{equation*}
    \lambda_2(\mclL_{\bfA}) =c' 4\pi^2/n^2\sum_{k =1}^d k^2 = cd^3/n^2,
\end{equation*}
for some factor $c\in (2\pi^2/3,4\pi^2/3)$.
\end{proof}


\clearpage


\subsection{Special cases of comparison graphs}\label{sec:appendix-special-case}

By \Cref{nthm:thm1}, for the estimator $\hat{\bbrtheta}_{\rho}$ to be
consistent, $L$ needs to be sufficiently large. We can check some common types
of comparison graph topologies and see in what order the necessary sample
complexity $N_{\textnormal{comp}} =|E|L$ needs to be, to achieve consistency. To
simplify results, we assume $e^{2\kappa_E}\lesssim \log n$. Spectral properties
of graphs listed here can be found in well-known textbooks \citep{spectra2012}.
\cite{shah2015estimationfrompairwisecomps} provide analogous comparisons. Per
\Cref{sec:special-case} we include results for the $d$-Caley graphs, and
expander graphs here. For reader convenience other results from
\Cref{sec:special-case} are also noted below.
\par
\textbf{Complete graph: } In this case, $\lambda_2(\mclL_{\bfA}) =n_{\max}
=n_{\min}= n-1$. Thus, we need $e^{2\kappa_E}\log n/n =o(L)$. Hence $L =
\Omega(1)$ and $N_{\textnormal{comp}}=\Omega (n^2)$.
\par
\textbf{Expander graph: } If the comparison graph is a $d$-regular expander
graph with edge expansion (or Cheeger number) coefficient $\phi$, then
$\lambda_2(\mclL_{\bfA})\geq \phi^2/(2d)$ \citep{expander2008}, $n_{\max} =
n_{\min} =d$. Here $\phi$ is defined by $\phi\defined\min_{|S|}{e(S,T)}/{|S|}$
where $\{S,T\}$ is a partition of the vertex set and $|S|\leq |T|$. We need
$d^2e^{2\kappa_E}/\phi^4\cdot (e^{2\kappa_E}n\vee d\log n) =o(L)$, so
$N_{\textnormal{comp}} = \Omega(nd^3e^{2\kappa_E}/\phi^4\cdot
(e^{2\kappa_E}n\vee d\log n))$.
\par
\textbf{Complete bipartite graph: } If the comparison graph has two partitioned
sets of size $m_1$ and $m_2$ such that $m_1\leq m_2$, then
$\lambda_2(\mclL_{\bfA}) = m_1$, $n_{\max} = m_2$, $n_{\min} = m_1$. We need
$e^{2\kappa_E}m_2/m_1^2\cdot [(e^{2\kappa_E}nm_2/m_1^2)\vee \log n] = o(L)$.
When $m_1=\Omega(n)$, we have $N_{\textnormal{comp}} = \Omega(n^2)$.
\par
\textbf{$d$-Cayley graph: } $(i,j)\in E$ if and only if $i - j\equiv k({\rm
mod}\ n)$ with $-d\leq k\leq d$, $k\neq 0$. It is a $2d$-regular graph, and
$\lambda_2(\mclL_{\bfA}) \asymp d^3/n^2$ (see \Cref{sec:others}),
$n_{\max}=n_{\min} = 2d$. Thus, we need $e^{2\kappa_E} n^4/d^6\cdot
(e^{2\kappa_E}n\vee 2d\log n)=o(L)$, so $N_{\textnormal{comp}} =
\Omega(e^{2\kappa_E} n^5/d^5\cdot (e^{2\kappa_E}n\vee 2d\log n))$ for $d = o(n)$
and $N_{\textnormal{comp}} = \Omega(n^2)$ for $d = \Omega(n)$.
\par
\textbf{Path or Cycle graph: } When comparisons occur based on a path or cycle
comparison graph, then by \Cref{prop:path}, $\|\hat{\bbrtheta}_0 -
\bbrtheta^*\|_{\infty}\lesssim e^{\kappa_E}\sqrt{\frac{n\log n}{L}}$. Thus, we
need $e^{2\kappa_E}n\log n = o(L)$ and $N_{\rm comp} =
\Omega(e^{2\kappa_E}n^2\log n)$.
\par
\textbf{Star graph: } A star graph on $n$ node is a tree graph with diameter $D
= 1$. By \Cref{prop:tree}, we have $\|\hat{\bbrtheta}_0 -
\bbrtheta^*\|_{\infty}\lesssim e^{\kappa_E}\sqrt{\frac{\log n}{L}}$. Thus, we
need $e^{2\kappa_E}\log n = o(L)$ and $N_{\rm comp} = \Omega(e^{2\kappa_E}n\log
n)$.
\par
\textbf{Barbell graph:} It contains two size-$n/2$ complete sub-graphs connected
by 1 edge, so $\lambda_2(\mclL_{\bfA}) \asymp 1/n$, $n_{\max} = n/2$, $n_{\min}
= n/2-1$. We need $e^{2\kappa_E}n^3\log n = o(L)$, and $N_{\textnormal{comp}} =
e^{2\kappa_E}n^5\log n$.



\clearpage


\subsection{Upper bound for unregularized/vanilla MLE}\label{sec:vanilla-mle}
\subsubsection{Main theorem}
The unregularized or vanilla MLE is defined as
\begin{equation}\label{eq:def_vanilla_mle}
    \hat{\bbrtheta} \defined \argmin_{\bbrtheta\in \mathbb{R}^n \colon {\bf 1}_n^\top \boldsymbol{\theta} = 0 } \ell(\bbrtheta;\bfy),
\end{equation}
where $\ell(\bbrtheta;\bfy)$ is the negative log-likelihood, given by
\begin{equation}
    \ell(\bbrtheta;\bfy) \defined -\sum_{1\leq i<j\leq n} A_{ij}\left\lbrace\bar{y}_{ij}\log {\psi(\theta_i - \theta_j)} + (1 - \bar{y}_{ij})\log [{1 - \psi(\theta_i - \theta_j)}]\right\rbrace,
\end{equation}
and $t \in \mathbb{R} \mapsto \psi(t) = {1}/{[1 + e^{-t}}]$ the sigmoid
function.  To make the expressions of results simpler, we consider the parameter
range $\kappa\leq n$ and $\kappa_E\leq \log n$ as discussed in the comments
after \Cref{nthm:thm1}. \bnthm[Vanilla MLE] \label{nthm:vanilla_mle} Assume
the BTL model with parameter $\bbrtheta^* = (\theta^*_1,\ldots,\theta^*_n)^\top$
such that ${\bf 1}_n^\top \bbrtheta^* =0$ and a comparison graph $\mathcal{G} =
\mclG([n],E)$ with adjacency matrix ${\bf A}$, algebraic connectivity
$\lambda_2(\mclL_\bfA)$ and maximum and minimum degrees $n_{\max}$ and
$n_{\min}$, respectively. Suppose that each pair of items $(i,j)\in E$ are
compared $L$ times. Let $\kappa = \max_{i,j}|\theta^*_i - \theta^*_j|$ and
$\kappa_E = \max_{(i,j)\in E}|\theta^*_i - \theta^*_j|$. Assume that $\mclG$ is
connected, or equivalently, $\lambda_2(\mclL_A) >0$. In addition, assume that 1.
$\lambda_2(\mclL_A)^2L>Ce^{2\kappa_E}\max\{n_{\max}\log n, e^{2\kappa_E}n\frac{n_{\max}^2}{n_{\min}^2}\}$ for some large constant $C>0$, and 2.
$\lambda_2(\mclL_{\bfA})\geq 2e^{2\kappa_E}n_{\max}/n_{\min}$. Then, with
probability at least $1 - O(n^{-5})$, the unregularized MLE  $\hat{\bbrtheta}$
from  \eqref{eq:def_vanilla_mle} satisfies
\begin{equation}\label{eq:vanilla_mle}
\|\hat{\theta} - \theta^*   \|_{\infty} \lesssim e^{\kappa_E}\sqrt{\frac{n_{\max}\log n}{Ln_{\min}^2}} + s\sqrt{\frac{n_{\max}}{L}}\left[ 1 + \frac{e^{\kappa_E}\sqrt{\log n} }{n_{\min}} +  s\sqrt{\frac{n}{n_{\max}}} \right],
\end{equation}
where $s\defined \frac{e^{2\kappa_E}n_{\max}}{\lambda_2 n_{\min}}$, and
\begin{equation}\label{eq:l2.rate_vanilla}
\|\hat{\bbrtheta}- \bbrtheta^*\|_2 \lesssim  \frac{e^{\kappa_E}}{\lambda_2(\mclL_\bfA)} \sqrt{\frac{n_{\max}n}{L}},
\end{equation}
provided that $L\lesssim n^5$, $\kappa < n$, $\kappa_E\leq \log n$, and the
right hand side of \Cref{eq:vanilla_mle} is smaller than a sufficiently small
constant $C>0$.
\enthm
\bnrmk
The expression of the $\ell_{\infty}$ bound looks messy, but in the $ER(n,p)$
case it reduces to the same form as the upper bound of the regularized MLE in
\Cref{cor:cor_ER}. Moreover, for vanilla MLE we require an additional pure
topological assumption $\lambda_2(\mclL_{\bfA})\geq
2e^{2\kappa_E}n_{\max}/n_{\min}$, which is stringent in some sense as it exclude
many graphs with small $n_{\min}$. But for such graphs with high degree
heterogeneity $n_{\max}/n_{\min}$, it's reasonable that we need some regularity
in our objective function. We believe that this condition can be weakened, and
the $\ell_\infty$ upper bound can be improved or tightened by improving the
proof techniques, which can be a good future direction for researchers.
\enrmk

To prove \Cref{nthm:vanilla_mle}, we need two lemmas, \Cref{lm:l_infty_const}
anbd \ref{lm:lem7.6}. \bnlem\label{lm:l_infty_const} Under the setting of
\Cref{nthm:vanilla_mle}, it holds with probability at least $1 - O(n^{-5})$ that
\begin{equation}
    \|\hat{\bbrtheta} - \bbrtheta^*\|_{\infty}\leq 5.
\end{equation}
\enlem
Following \cite{chen2020partialtopkranking}, we decompose the full negative
loglikelihood function as
\begin{equation}\label{eq:decompose_ell}
    \ell_n(\bbrtheta) = \ell_{n}^{(-m)}\left(\bbrtheta_{-m}\right) + \ell_{n}^{(m)}\left(\bbrtheta_{-m}\right),
\end{equation}
where $\theta_m\in \mathbb{R}$ is the $m$-th entry of $\bbrtheta$ and
$\bbrtheta_{-m}\in \mathbb{R}^{n-1}$ is the subvector containing the rest of
entries, and the two functions are given by
\begin{align*}
    \ell_{n}^{(-m)}\left(\bbrtheta_{-m}\right) &=\sum_{1 \leq i<j \leq n: i, j \neq m} A_{i j}\left[\bar{y}_{i j} \log \frac{1}{\psi\left(\theta_{i}-\theta_{j}\right)}+\left(1-\bar{y}_{i j}\right) \log \frac{1}{1-\psi\left(\theta_{i}-\theta_{j}\right)}\right],\\
    \ell_{n}^{(m)}\left(\theta_m|\bbrtheta_{-m}\right) &=\sum_{j\in [n]\setminus\{m\}} A_{m j}\left[\bar{y}_{m j} \log \frac{1}{\psi\left(\theta_{m}-\theta_{j}\right)}+\left(1-\bar{y}_{m j}\right) \log \frac{1}{1-\psi\left(\theta_{m}-\theta_{j}\right)}\right].
\end{align*}
Let $H^{(-m)}:=\nabla^2\ell_n^{(-m)}(\ell_{-m})$, and
\begin{equation*}
    \bbrtheta_{-m}^{(m)}=\underset{\bbrtheta_{-m}: \| \bbrtheta_{-m}-\bbrtheta_{-m}^{*}\|_{\infty \leq 5}}{\operatorname{argmin}} \ell_{n}^{(-m)}\left(\bbrtheta_{-m}\right).
\end{equation*}

\bnlem\label{lm:lem7.6} Under the setting of \Cref{nthm:vanilla_mle}, it holds
with probability at least $1 - O(n^{-9})$ that
\begin{equation}
    \max _{m \in[n]}\|\theta_{-m}^{(m)}-\theta_{-m}^{*}-a_{m} \textbf{1}_{n-1}\|^{2}_2 \leq C \frac{e^{2\kappa_E}n_{\max}n}{L\lambda_2(\mathcal{L}_{A_{}})^2}
\end{equation}
for some constant $C>0$, where
$a_{m}=\operatorname{avg}\left(\theta_{-m}^{(m)}-\bbrtheta_{-m}^{*}\right):=\frac{1}{n-1}\mathbf{1}_{n-1}^\top
\left(\theta_{-m}^{(m)}-\bbrtheta_{-m}^{*}\right)$.
\enlem

\bprfof{\Cref{nthm:vanilla_mle}} Again we define the leave-one-out negative
log-likelihood as
\begin{equation}
\begin{aligned}
\ell_{n}^{(m)}(\theta)=& \sum_{1 \leq i<j \leq n: i, j \neq m} A_{i j}\left[\bar{y}_{i j} \log \frac{1}{\psi\left(\theta_{i}-\theta_{j}\right)}+\left(1-\bar{y}_{i j}\right) \log \frac{1}{1-\psi\left(\theta_{i}-\theta_{j}\right)}\right] \\
&+\sum_{i \in[n] \backslash\{m\}} A_{mj}\left[\psi\left(\theta_{i}^{*}-\theta_{m}^{*}\right) \log \frac{1}{\psi\left(\theta_{i}-\theta_{m}\right)}+\psi\left(\theta_{m}^{*}-\theta_{i}^{*}\right) \log \frac{1}{\psi\left(\theta_{m}-\theta_{i}\right)}\right],
\end{aligned}
\end{equation}
For the $\ell_2$ bound, notice that by Taylor expansion
\begin{equation*}
    \ell_{n}(\hat{\theta})=\ell_{n}(\theta^{*})+(\hat{\theta}-\theta^{*})^{T} \nabla \ell_{n}(\theta^{*})+\frac{1}{2}(\hat{\theta}-\theta^{*})^{T} H(\xi)(\hat{\theta}-\theta^{*}),
\end{equation*}
where $\xi$ is a convex combination of $\hat{\bbrtheta}$ and $\bbrtheta^*$. By
\Cref{lm:l_infty_const}, we have $\|\hat{\bbrtheta}-\bbrtheta^*\|_{\infty}\leq
5$ and hence $\|\xi-\bbrtheta^*\|_{\infty}\leq 5$. By \Cref{lm:lem8}, we have
$\frac{1}{2}(\hat{\theta}-\theta^{*})^{T} H(\xi)(\hat{\theta}-\theta^{*})\geq
ce^{-\kappa_E}\lambda_2\|\hat{\bbrtheta}-\bbrtheta^*\|^2_2$ for some constant
$c>0$. By the fact that $\ell_{n}(\bbrtheta^*)\geq \ell_n(\hat{\bbrtheta})$,
Cauchy-Schwartz inequality, and \Cref{lm:lem7}, we have
\begin{equation*}
    \|\hat{\bbrtheta}-\bbrtheta^*\|_2\lesssim \frac{e^{\kappa_E}}{\lambda_2}\|\nabla\ell(\bbrtheta^*)\|_2\leq \frac{e^{\kappa_E}}{\lambda_2}\sqrt{\frac{n_{\max} n}{L}}.
\end{equation*}
For the $\ell_{\infty}$ bound, the proof can be sketched as following steps.
\begin{enumerate}
    \item[0.] By \Cref{lm:l_infty_const}, $\|\hat{\theta} -
    \theta^*\|_{\infty}\leq 5$ with probability at least $1 - O(n^{-5})$.
    \item By \Cref{lm:lem7.6}, it holds with probability exceeding $1 -
    O(n^{-9})$ that,
\begin{equation*}
        \max _{m \in[n]}\|\theta_{-m}^{(m)}-\theta_{-m}^{*}-a_{m} \textbf{1}_{n-1}\|^{2}_2 \leq C \frac{e^{2\kappa_E}n_{\max}n}{L\lambda_2(\mathcal{L}_A)^2}
\end{equation*}
    where
    $a_{m}=\operatorname{avg}\left(\theta_{-m}^{(m)}-\bbrtheta_{-m}^{*}\right):=\frac{1}{n-1}\mathbf{1}_{n-1}^\top
    \left(\theta_{-m}^{(m)}-\bbrtheta_{-m}^{*}\right)$.
    \item Show that on the same event,
\begin{equation}
\begin{aligned}
\max _{m \in[n]}\|\theta_{-m}^{(m)}-\widehat{\theta}_{-m}-a_{m} \textbf{1}_{n-1}\|^{2}_2 \leq & C_{1} \frac{\max _{m \in[n]} \sum_{i \in \mathcal{N}(m)} \left(\bar{y}_{m i}-\psi\left(\theta_{m}^{*}-\theta_{i}^{*}\right)\right)^{2}}{\lambda_2(\mathcal{L}_A)^2 e^{-2\kappa_E}} \\
&+C_{1} \frac{e^{2\kappa_E}n_{\max}}{\lambda_2(\mathcal{L}_A)^2}{\|\widehat{\theta}-\theta^{*}\|_{\infty}^{2}}.
\end{aligned}
\label{eq:l_infty_g1}
\end{equation}
Following the same arguments towards \Cref{eq:bound_l2_subset}, we can get
\begin{equation*}
    \|\bbrtheta_{-m}^{(m)}-\hat{\bbrtheta}_{-m}-\bar{a}_{m} \mathbf{1}_{n-1}\|_2^2\leq \frac{e^{2\kappa_E}}{c^2\lambda_2(\mathcal{L}_{A})^2}\|\nabla\ell_n^{(-m)}(\hat{\bbrtheta}_{-m})\|_2^2.
\end{equation*}
By \Cref{eq:decompose_ell}, for each $i\in [n]\setminus \{m\}$, we have
\begin{equation*}
    \frac{\partial}{\partial \bbrtheta_{i}} \ell_{n}^{(-m)}\left(\bbrtheta_{-m}\right)=\frac{\partial}{\partial \theta_{i}} \ell_{n}(\bbrtheta)-\frac{\partial}{\partial \theta_{i}} \ell_{n}^{(m)}\left(\theta_{m} \mid \bbrtheta_{-m}\right),
\end{equation*}
Using the fact that $\nabla \ell_n(\hat{\bbrtheta}) = 0$, we get
\begin{equation*}
    \frac{\partial}{\partial \theta_{i}} \ell_{n}^{(-m)}\left(\bbrtheta_{-m}\right)_{\mid \bbrtheta=\widehat{\bbrtheta}}=-\frac{\partial}{\partial \theta_{i}} \ell_{n}^{(m)}\left(\theta_{m} \mid \bbrtheta_{-m}\right)\mid_{\theta=\widehat{\bbrtheta}}=-A_{m i}\left[\bar{y}_{m i}-\psi(\widehat{\theta}_{m}-\widehat{\theta}_{i})\right].
\end{equation*}
Therefore, we have
\begin{align*}
\|\nabla \ell_{n}^{(-m)}(\widehat{\bbrtheta}_{-m})\|_2^{2}=& \sum_{i \in[n] \backslash\{m\}} A_{m i}\left[\bar{y}_{m i}-\psi(\widehat{\theta}_{m}-\widehat{\theta}_{i})\right]^{2} \\
\leq & 2 \sum_{i \in[n] \backslash\{m\}} A_{m i}\left(\bar{y}_{m i}-\psi\left(\theta_{m}^{*}-\theta_{i}^{*}\right)\right)^{2} \\
&+2 \sum_{i \in[n] \backslash\{m\}} A_{m i}\left[\psi\left(\theta_{m}^{*}-\theta_{i}^{*}\right)-\psi(\widehat{\theta}_{m}-\widehat{\theta}_{i})\right]^{2} \\
\leq & 2 \sum_{i \in[n] \backslash\{m\}} A_{m i}\left(\bar{y}_{m i}-\psi\left(\theta_{m}^{*}-\theta_{i}^{*}\right)\right)^{2}+2\|\widehat{\theta}-\theta^{*}\|_{\infty}^{2} \sum_{i \in[n] \backslash\{m\}} A_{m i} \\
\leq & 2 \sum_{i \in[n] \backslash\{m\}} A_{m i}\left(\bar{y}_{m i}-\psi\left(\theta_{m}^{*}-\theta_{i}^{*}\right)\right)^{2}+4 n_{\max}\|\widehat{\theta}-\theta^{*}\|_{\infty}^{2}.
\end{align*}
Now use the fact that $\mathbf{1}_n^\top \bbrtheta^*=\mathbf{1}_n^\top
\hat{\bbrtheta} = 0$, we have
\begin{equation*}
    \left\|a_{m} \mathbf{1}_{n-1}-\bar{a}_{m} \mathbf{1}_{n-1}\right\|^{2}_2=(n-1)[\operatorname{avg}(\widehat{\bbrtheta}_{-m}-\bbrtheta_{-m}^{*})]^{2}=\frac{(\hat{\bbrtheta}_{m}-\theta_{m}^{*})^{2}}{n-1} \leq \frac{\|\hat{\bbrtheta}-\bbrtheta^{*}\|_{\infty}^{2}}{n-1} .
\end{equation*}
These results, together with the fact that $\lambda_2(\mclL_{\bfA})\leq
2n_{\max}\leq 2n$, give \Cref{eq:l_infty_g1}.
    \item Show that on the same event,
\begin{equation}
\begin{aligned}
\|\widehat{\theta}-\theta^{*}\|_{\infty}\cdot C_{4}e^{-\kappa_E}n_{\min} \leq &  \max _{m \in[n]} \mid \sum_{i \in\mathcal{N}(m)} \left(\bar{y}_{m i}-\psi\left(\theta_{m}^{*}-\theta_{i}^{*}\right)\right)| \\
&+\sqrt{n_{\max}}\max _{m \in[n]}\|\theta_{-m}^{(m)}-\theta_{-m}^{*}-a_{m} \textbf{1}_{n-1}\|_2 \\
&+\sqrt{n_{\max}}\max_{m\in[n]}\|\theta_{-m}^{(m)}-\widehat{\theta}_{-m}-a_{m} \textbf{1}_{n-1}\|_2.
\end{aligned}
\label{eq:l_infty_g2}
\end{equation}
First, define two univariate functions as some proxy of gradient and hessian:
\begin{align*}
g^{(m)}\left(\theta_{m} \mid \bbrtheta_{-m}\right) &= \frac{\partial}{\partial \theta_{m}} \ell_{n}^{(m)}\left(\theta_{m} \mid \bbrtheta_{-m}\right)=-\sum_{i \in[n] \backslash\{m\}} A_{m i}\left(\bar{y}_{m i}-\psi\left(\theta_{m}-\theta_{i}\right)\right) \\
h^{(m)}\left(\theta_{m} \mid \bbrtheta_{-m}\right) &= \frac{\partial^{2}}{\partial \theta_{m}^{2}} \ell_{n}^{(m)}\left(\theta_{m} \mid \bbrtheta_{-m}\right)=\sum_{i \in[n] \backslash\{m\}} A_{m i} \psi\left(\theta_{m}-\theta_{i}\right) \psi\left(\theta_{i}-\theta_{m}\right) .
\end{align*}
By the definition of $\hat{\bbrtheta}$ and the shift invariance of $\ell_n$, we
have $\ell_n(\hat{\bbrtheta})\leq \ell_n(\bbrtheta)$ for any $\bbrtheta\in
\mathbb{R}^n$, thus
\begin{equation*}
    \ell_{n}^{(m)}(\theta_{m}^{*} \mid \hat{\bbrtheta}_{-m})+\ell_{n}^{(-m)}(\hat{\bbrtheta}_{-m}) \geq \ell_{n}(\hat{\bbrtheta}).
\end{equation*}
This implies
\begin{align*}
\ell_{n}^{(m)}(\theta_{m}^{*} \mid \hat{\theta}_{-m}) & \geq \ell_{n}^{(m)}(\hat{\theta}_{m} \mid \hat{\theta}_{-m}) \\
&=\ell_{n}^{(m)}(\theta_{m}^{*} \mid \hat{\bbrtheta}_{-m})+(\hat{\theta}_{m}-\theta_{m}^{*}) g^{(m)}(\theta_{m}^{*} \mid \hat{\bbrtheta}_{-m})+\frac{1}{2}(\hat{\theta}_{m}-\theta_{m}^{*})^{2} h^{(m)}(\xi \mid \hat{\bbrtheta}_{-m}),
\end{align*}
where $\xi$ is a convex combination of $\theta_m^*$ and $\hat{\theta}_m$. By
\Cref{lm:l_infty_const}, we have $|\xi - \theta_m^*|\leq |\hat{\theta}_m -
\theta_m^*|\leq 5$. Thus for any $i\neq m$ it holds that $|\xi -
\hat{\theta}_i|\leq
\left|\xi-\theta_{m}^{*}\right|+\left|\theta_{m}^{*}-\theta_{i}^{*}\right|+\left|\widehat{\theta}_{i}-\theta_{i}^{*}\right|
\leq 10+\kappa$. By definition of $h^{(m)}$, we have
$\frac{1}{2}h^{(m)}(\xi|\hat{\bbrtheta}_{-m})\geq c_2e^{-\kappa_E}n_{\min}$ for
some constant $c_2>0$. Therefore, we get
\begin{equation}\label{eq:gao_48}
    (\hat{\theta}_{m}-\theta_{m}^{*})^{2} \leq \frac{e^{2\kappa_E}}{(c_{2} n_{\min})^{2}}|g^{(m)}(\theta_{m}^{*} \mid \hat{\bbrtheta}_{-m})|^{2}.
\end{equation}
To bound $|g^{(m)}(\theta_{m}^{*} \mid \hat{\bbrtheta}_{-m})|$, we decompose it
as
\begin{align}
|g^{(m)}(\theta_{m}^{*} \mid \hat{\theta}_{-m})|=&\left|\sum_{i \in[n] \backslash\{m\}} A_{m i}(\bar{y}_{m i}-\psi(\theta_{m}^{*}-\widehat{\theta}_{i}))\right| \nonumber\\
\leq &\left|\sum_{i \in[n] \backslash\{m\}} A_{m i}(\bar{y}_{m i}-\psi(\theta_{m}^{*}-\theta_{i}^{*}))\right|   \label{eq:gao_49} \\
&+\left|\sum_{i \in[n] \backslash\{m\}} A_{m i}(\psi(\theta_{m}^{*}-\theta_{i}^{*})-\psi(\theta_{m}^{*}-\theta_{i}^{(m)}+a_{m}))\right|   \label{eq:gao_50}\\
&+\left|\sum_{i \in[n] \backslash\{m\}} A_{m i}(\psi(\theta_{m}^{*}-\theta_{i}^{(m)}+a_{m})-\psi(\theta_{m}^{*}-\widehat{\theta}_{i}))\right|  \label{eq:gao_51}.
\end{align}
By Cauchy-Schwartz inequality, we can bound \eqref{eq:gao_50} and
\eqref{eq:gao_51} by
\begin{align*}
        \left|\sum_{i \in[n] \backslash\{m\}} A_{m i}(\psi(\theta_{m}^{*}-\theta_{i}^{*})-\psi(\theta_{m}^{*}-\theta_{i}^{(m)}+a_{m}))\right|^2 &\leq n_{\max}\|\bbrtheta_{-m}^{(m)}-\bbrtheta_{-m}^{*} - a_m\mathbf{1}_{n-1}\|^2_2,\\
        \left|\sum_{i \in[n] \backslash\{m\}} A_{m i}(\psi(\theta_{m}^{*}-\theta_{i}^{(m)}+a_{m})-\psi(\theta_{m}^{*}-\widehat{\theta}_{i}))\right|^2 &\leq n_{\max}\|\bbrtheta_{-m}^{(m)}-\hat{\bbrtheta}_{-m} - a_m\mathbf{1}_{n-1}\|^2_2.
\end{align*}
Plugging these bounds into \Cref{eq:gao_48} and taking maximum over $m\in [n]$
give the desired bound \eqref{eq:l_infty_g2}.
    \item Plug \eqref{eq:l_infty_g2} back into \eqref{eq:l_infty_g1} and get
\begin{equation*}
\begin{split}
    \max _{m \in[n]}\|\theta_{-m}^{(m)}-\widehat{\theta}_{-m}-a_{m} \mathbf{1}_{n-1}\|_2^2 &\lesssim \frac{e^{2\kappa_E}}{\lambda_2(\mathcal{L}_A)^2}\max_m \sum_{i\in \mathcal{N}(m)} (\bar{y}_{mi} - \psi^*_{mi})^2\\
    +\frac{e^{4\kappa_E}}{\lambda_2(\mathcal{L}_A)^2}\frac{n_{\max}}{n_{\min}^2} & \max_{m\in [n]}|\sum_{i\in \mathcal{N}(m)}(\bar{y}_{mi} - \psi^*_{mi})|^2\\
    +\frac{e^{4\kappa_E}}{\lambda_2(\mathcal{L}_A)^2}\frac{n_{\max}^2}{n_{\min}^2}
    & \max _{m \in[n]}\|\theta_{-m}^{(m)}-\theta_{-m}^{*}-a_{m} \textbf{1}_{n-1}\|^{2}_2\\
    +\frac{e^{4\kappa_E}}{\lambda_2(\mathcal{L}_A)^2}\frac{n_{\max}^2}{n_{\min}^2} & \max _{m \in[n]}\|\theta_{-m}^{(m)}-\hat{\theta}_{-m}-a_{m} \textbf{1}_{n-1}\|^{2}_2.
\end{split}
\end{equation*}
as we assume
$\frac{e^{2\kappa_E}}{\lambda_2(\mathcal{L}_A)}\frac{n_{\max}}{n_{\min}} \leq
\frac{1}{2}$, we have
\begin{equation} \label{eq:l_infty_g3}
    \begin{split}
        \max _{m \in[n]}\|\theta_{-m}^{(m)}-\widehat{\theta}_{-m}-a_{m} \mathbf{1}_{n-1}\|_2^2\lesssim &\frac{e^{2\kappa_E}(\log n + n_{\max})}{\lambda_2(\mathcal{L}_A)^2 L} + \frac{e^{4\kappa_E}}{\lambda_2(\mathcal{L}_A)^2}\frac{n_{\max}^2}{n_{\min}^2}\frac{\log n}{L} +\\
&\frac{e^{4\kappa_E}}{\lambda_2(\mathcal{L}_A)^2}
\frac{n_{\max}^2}{n_{\min}^2}
\frac{e^{2\kappa_E}n_{\max}n}{\lambda_2(\mathcal{L}_A)^2 L}.
    \end{split}
\end{equation}
    \item Plug \eqref{eq:l_infty_g3} back into \eqref{eq:l_infty_g2} and we can
    get
\begin{equation}
\begin{split}
    \|\widehat{\theta}-\theta^{*}\|_{\infty}^2 &\lesssim  \frac{e^{2\kappa_E}}{n_{\min}^2}\frac{n_{\max}\log n }{ L} + \frac{e^{2\kappa_E}n^2_{\max}}{n_{\min}^2n}\frac{e^{2\kappa_E}n_{\max} n }{\lambda_2(\mathcal{L}_A)^2 L} \\
    + & \frac{e^{2\kappa_E}n_{\max}}{n_{\min}^2\lambda_2(\mathcal{L}_A)^2 L}\left[ e^{2\kappa_E}(\log n + n_{\max}) + \frac{e^{4\kappa_E}n_{\max}^2 \log n}{n_{\min}^2} +
    \frac{e^{6\kappa_E}n_{\max}^3n}{\lambda_2(\mathcal{L}_A)^2 n_{\min}^2} \right]
\end{split}
\end{equation}
One term can be reduced and the inequality becomes
\begin{equation}
\|\hat{\theta} - \theta^*   \|_{\infty} \lesssim e^{\kappa_E}\sqrt{\frac{n_{\max}\log n}{Ln_{\min}^2}} + e^{2\kappa_E}\sqrt{\frac{n^3_{\max}}{L\lambda_2(\mathcal{L}_A)^2n_{\min}^2}}\left[ 1 + e^{\kappa_E}\sqrt{\frac{\log n }{n_{\min}^2}} + e^{2\kappa_E} \sqrt{\frac{ n_{\max}n}{\lambda_2(\mathcal{L}_A)^2 n_{\min}^2}} \right].
\end{equation}
\end{enumerate}
\eprfof

\subsubsection{Proof of Lemmas}
\bprfof{\Cref{lm:l_infty_const}} We will use a gradient descent sequence defined
by
\begin{equation*}
    \theta^{(t+1)}=\theta^{(t)}-\eta[\nabla \ell_{n}(\theta^{(t)})+\rho \theta^{(t)}].
\end{equation*}
The leave-one-out negative log-likelihood is defined as
\begin{equation*}
\begin{aligned}
\ell_{n}^{(m)}(\theta)=& \sum_{1 \leq i<j \leq n: i, j \neq m} A_{i j}\left[\bar{y}_{i j} \log \frac{1}{\psi\left(\theta_{i}-\theta_{j}\right)}+\left(1-\bar{y}_{i j}\right) \log \frac{1}{1-\psi\left(\theta_{i}-\theta_{j}\right)}\right] \\
&+\sum_{i \in[n] \backslash\{m\}} A_{mi}\left[\psi\left(\theta_{i}^{*}-\theta_{m}^{*}\right) \log \frac{1}{\psi\left(\theta_{i}-\theta_{m}\right)}+\psi\left(\theta_{m}^{*}-\theta_{i}^{*}\right) \log \frac{1}{\psi\left(\theta_{m}-\theta_{i}\right)}\right],
\end{aligned}
\end{equation*}
so the leave-one-out gradient descent sequence is defined as
\begin{equation*}
    \theta^{(t+1, m)}=\theta^{(t, m)}-\eta[\nabla \ell_{n}^{(m)}(\theta^{(t, m)})+\rho\theta^{(t, m)}].
\end{equation*}
We initialize both sequences by $\theta^{(0)} = \theta^{(0,m)} = \theta^*$ and
set $\rho = \frac{1}{\kappa}\sqrt{\frac{n_{\max}}{L}}$ and step size $\eta =
\frac{1}{\lambda + n_{\max}}$. We will show that under the assumption
$\lambda_2(\mclL_A)^2L>Ce^{2\kappa_E}\max\{n_{\max}\log n, e^{2\kappa_E} n
n_{\max}^2/n_{\min}^2\}$ for some large constant $C>0$, we have
\begin{equation}
\begin{split}
\max _{m \in[n]}\|\theta^{(t, m)}-\theta^{(t)}\|_2 \leq & f_1:=C_1\frac{e^{\kappa_E}}{\lambda_2}\sqrt{\frac{n_{\max}\log n}{L}}\leq 1 \\
\|\theta^{(t)}-\theta^{*}\|_2 \leq&f_2:= C_2\frac{e^{\kappa_E}}{\lambda_2}\sqrt{\frac{n_{\max} n}{L}} \leq \sqrt{\frac{n}{\log n}} \\
\max _{m \in[n]}|\theta_{m}^{(t, m)}-\theta_{m}^{*}| \leq&f_3:= C_3\frac{e^{2\kappa_E}}{\lambda_2}\frac{n_{\max}}{n_{\min}}\sqrt{\frac{n}{L}}\leq 1.
\end{split}
\label{eq:3bound}
\end{equation}
A useful fact given that \eqref{eq:3bound} holds is that
\begin{equation}
    \|\theta^{(t,m)}-\theta^*\|_\infty \leq f_1 + f_2.
\label{eq:eq36}
\end{equation}
We again have the Taylor expansion
\begin{equation*}
    \theta^{(t+1)} - \theta^{(t+1,m)}  =[(1-\eta \rho) I_{n}-\eta H(\xi)](\theta^{(t)}-\theta^{(t, m)})-\eta[\nabla \ell_{n}(\theta^{(t, m)})-\nabla \ell_{n}^{(m)}(\theta^{(t, m)})].
\end{equation*}
Now by the fact that  $\lambda_{\min , \perp}(H(\xi))\geq c_0
e^{-\kappa}\lambda_2$, we have
\begin{equation*}
    \|\left((1-\eta \rho) I_{n}-\eta H(\xi)\right)(\theta^{(t)}-\theta^{(t, m)})\|_2 \leq(1-\eta \rho-c_{1} \eta \lambda_2)\|\theta^{(t)}-\theta^{(t, m)}\|_2
\end{equation*}
for some constant $c_1>0$ and the other term can be bounded as
\begin{equation*}
\begin{aligned}
&\|\nabla \ell_{n}(\theta^{(t, m)})-\nabla \ell_{n}^{(m)}(\theta^{(t, m)})\|_2^{2} \\
=&\left[\sum_{j \in[n] \backslash\{m\}} A_{j m}\left(\bar{y}_{j m}-\psi\left(\theta_{j}^{*}-\theta_{m}^{*}\right)\right)\right]^{2} + \sum_{j \in[n] \backslash\{m\}} A_{j m}\left(\bar{y}_{j m}-\psi\left(\theta_{j}^{*}-\theta_{m}^{*}\right)\right)^{2}\\
\leq& C_1 \frac{1}{L} n_{\max}\log n + C_1 \frac{1}{L}(\log n + n_{\max}).
\end{aligned}
\end{equation*}
Therefore, for
\begin{equation*}
    \|\theta^{(t+1)} - \theta^{(t+1,m)}\|_2\leq (1 - c_1 \eta \lambda_2)f_1 + \eta \sqrt{ 2 C_1 \frac{1}{L}n_{\max}\log n}\leq f_1
\end{equation*}
to hold, we need $f_1>C\frac{e^{\kappa_E}}{\lambda_2}\sqrt{\frac{n_{\max}\log
n}{L}}$ for some sufficiently large positive constant $C>0$.
\par
Next, we bound $\|\theta^{(t+1)} - \theta^*\|$. By Tarlor expansion,
\begin{equation*}
    \theta^{(t+1)} - \theta^* = \left((1-\eta \lambda) I_{n}-\eta H(\xi)\right)(\theta^{(t)}-\theta^{*})-\eta \lambda \theta^{*}-\eta \nabla \ell_{n}\left(\theta^{*}\right).
\end{equation*}
Equation (41) becomes
\begin{equation*}
    \left((1-\eta \lambda) I_{n}-\eta H(\xi)\right)(\theta^{(t)}-\theta^{*}) \leq\left(1-\eta \lambda-c_{2} \eta \lambda_2\right)\|\theta^{(t)}-\theta^{*}\|_2,
\end{equation*}
and equation (42) becomes
\begin{equation*}
    \left\|\nabla \ell_{n}\left(\theta^{*}\right)\right\|_2^{2}=\sum_{i=1}^{n}\left[\sum_{j \in[n] \backslash\{i\}} A_{i j}\left(\bar{y}_{i j}-\psi\left(\theta_{i}^{*}-\theta_{j}^{*}\right)\right)\right]^{2} \leq C_{2} \frac{nn_{\max}}{L} ,
\end{equation*}
for some constants $c_2,C_2>0$. Therefore,
\begin{equation*}
    \|\theta^{(t+1)} - \theta^*\|_2 \leq (1 - c_2\eta\lambda_2)f_2+ \eta \sqrt{C_2\frac{n n_{\max}}{L}} + \eta\lambda \|\theta^*\|_2.
\end{equation*}
For $\|\theta^{(t+1)} - \theta^*\|\leq f_2$ to hold, we need
$\frac{\lambda_2}{e^{\kappa_E}}f_2>C\sqrt{\frac{n_{\max}n}{L}}$ for some
sufficiently large constant $C$, which is guaranteed by the definition of $f_2$.
\par
Next, we bound $|\theta^{(t+1,m)}_m - \theta_m^*|$. Note that by the definition
of the gradient descent
\begin{equation*}
    \theta^{(t+1,m)}_m - \theta_m^* = \left[ 1 - \eta\lambda -\eta \sum_{j\in [n]\backslash\{m\} }A_{mj}\psi{'} (\xi_j) \right](\theta^{(t,m)}_m - \theta_m^*) - \lambda \eta \theta_m^* + \eta \sum_{j\in [n]\backslash\{m\}}A_{mj}\psi{'}(\xi_j) (\theta^{(t,m)}_j - \theta_j^*)
\end{equation*}
where $\xi_j$ is a scalar between $\theta_m^*-\theta_j^*$ and $\theta^{(t,m)}_m
- \theta^{(t,m)}_j$. Since $\|\theta^{(t,m)}-\theta^*\|_\infty \leq 3$, we have
$|\xi_j - \theta_m^* + \theta_j^*|\leq |\theta_m^* - \theta_j^*
-\theta_m^{(t,m)} + \theta_j^{(t,m)}|\leq 6$ and $\|\xi\|_{\infty}$ is bounded.
Therefore,
\begin{equation*}
    \sum_{j\in [n]\backslash\{m\} }A_{mj}\psi{'} (\xi_j)\geq c_3\min_{i\in [n]}n_i
\end{equation*}
and
\begin{equation*}
\begin{split}
        |\sum_{j\in \mathcal{N}(m)} A_{mj}\psi'(\xi_j)(\theta_j^{(t,m)} - \theta_j^*)| &\leq \sqrt{\sum_{j\in \mathcal{N}(m)} [A_{mj}\psi'(\xi_j)]^2}\sqrt{\sum_{j\in \mathcal{N}(m)}(\theta_j^{(t,m)} - \theta_j^*)^2}\\
        &\leq c_4 \sqrt{n_{\max}} (f_1 + f_2).
\end{split}
\end{equation*}
for some constant $c_3, c_4>0$. Thus we have
\begin{equation*}
    |\theta^{(t+1,m)}_m - \theta_m^*|\leq (1 - c_3 \eta n_{\min}) + \eta \sqrt{n_{\max}}(f_1 + f_2) + \lambda \eta |\theta_m^*|.
\end{equation*}
For $|\theta^{(t+1,m)}_m - \theta_m^*|\leq f_3$ to hold, we need
$\frac{n_{\min}}{e^{\kappa_E}}f_3>C\sqrt{n_{\max}}f_2$ for some sufficiently
large constant $C>0$, which is ensured by the definition of $f_2,f_3$.

As the last step, we again use the fact that $\ell_\rho(\cdot)$ is
$\rho$-strongly convex and $(\rho + n_{\max})$-smooth (see definition in the
paragraph before \Cref{eq:linear_conv}), so by Theorem 3.10 in
\cite{Bubeck2015}, we have
\begin{equation}
    \|\theta^{(t)} - \hat{\theta}_\rho\|_2\leq (1 - \frac{\rho}{\rho + n_{\max}})^t\|\theta^{*} - \hat{\theta}_\rho\|_2.
\end{equation}
By a union bound, \Cref{eq:3bound} holds for all $t\leq T$ with probability at
leaast $1 - O(Tn^{-10})$. Triangle inequality implies that
\begin{equation*}
    \|\hat{\theta}_{\rho}-\theta^{*}\|_{\infty} \leq\|\theta^{(T)}-\hat{\theta}_{\rho}\|_2+\|\theta^{(T)}-\theta^{*}\|_{\infty} \leq(1-\frac{\rho}{\rho+n_{\max}})^{T} \sqrt{n}\|\hat{\theta}_{\rho}-\theta^{*}\|_{\infty}+2
\end{equation*}

Take $T = n^5$ and remember that $L\lesssim n^5$. If $\rho > n_{\max}$, then $(1
- \frac{\rho}{\rho + n_{\max}})^T \sqrt{n} \leq 2^{-n^5}\sqrt{n}\leq 1/2$.
Otherwise, since
\begin{equation*}
(1 - \frac{\rho}{\rho + n_{\max}})^T \sqrt{n}\leq \exp\left( - \frac{T\rho}{\rho + n_{\max}}\right)\sqrt{n},
\end{equation*}
using the fact that $\kappa < n$, we have
\begin{equation*}
(1 - \frac{\rho}{\rho + n_{\max}})^T \sqrt{n}\leq \exp\left(-\frac{T}{c\kappa}\sqrt{\frac{1}{n_{\max} L}}\right)\kappa \sqrt{n}\leq ce^{-n}n^{3/2}\leq  \frac{1}{2}.
\end{equation*}
In conclusion, we have $\|\hat{\theta}_{\rho}-\theta^{*}\|_{\infty}\leq
\frac{1}{2}\|\hat{\theta}_{\rho}-\theta^{*}\|_{\infty} + 2$, thus
$\|\hat{\theta}_{\rho}-\theta^{*}\|_{\infty}\leq 4$, with probability at least
$1 - O(n^{-5})$.
\eprfof

\bprfof{\Cref{lm:lem7.6}} By definition, $\bbrtheta_{-m}^{(m)}$ is a constrained
MLE on a subset of the data, thus by Taylor expansion, for $\xi$ given by a
convex combination of $\bbrtheta_{-m}^{*}$ and $\bbrtheta_{-m}^{(m)}$, we have
\begin{align*}
\ell_{n}^{(-m)}\left(\bbrtheta_{-m}^{*}\right) \geq & \ell_{n}^{(-m)}(\bbrtheta_{-m}^{(m)}) \\
=& \ell_{n}^{(-m)}\left(\bbrtheta_{-m}^{*}\right)+(\bbrtheta_{-m}^{(m)}-\bbrtheta_{-m}^{*}-a_{m} \mathbf{1}_{n-1})^{T} \nabla \ell_{n}^{(-m)}\left(\bbrtheta_{-m}^{*}\right) \\
&+\frac{1}{2}(\bbrtheta_{-m}^{(m)}-\bbrtheta_{-m}^{*}-a_{m} \mathbf{1}_{n-1})^{T} H^{(-m)}(\xi)(\bbrtheta_{-m}^{(m)}-\bbrtheta_{-m}^{*}-a_{m} \mathbf{1}_{n-1}),
\end{align*}
where we use the invariant property of
$\ell_{n}^{(-m)}\left(\bbrtheta_{-m}\right)$, i.e.,
$\ell_{n}^{(-m)}\left(\bbrtheta_{-m}\right) =
\ell_{n}^{(-m)}\left(\bbrtheta_{-m} + c\mathbf{1}_{n-1}\right)$. By the fact
that $\|\xi - \bbrtheta^*_{-m}\|_{\infty}\leq \|\bbrtheta^{(m)}_{-m} -
\bbrtheta^*_{-m}\|_{\infty}\leq 5$ and \Cref{lm:lem8}, we have
\begin{equation*}
    (\bbrtheta_{-m}^{(m)}-\bbrtheta_{-m}^{*}-a_{m} \mathbf{1}_{n-1})^{T} H^{(-m)}(\xi)(\bbrtheta_{-m}^{(m)}-\bbrtheta_{-m}^{*}-a_{m} \mathbf{1}_{n-1})\geq ce^{-\kappa_E}\lambda_2(\mathcal{L}_{A_{-m}})\|\bbrtheta_{-m}^{(m)}-\bbrtheta_{-m}^{*}-a_{m} \mathbf{1}_{n-1}\|_2^2.
\end{equation*}
Applying Cauchy-Schwartz inequality to the expansion and we can get
\begin{equation*}
    \|\bbrtheta_{-m}^{(m)}-\bbrtheta_{-m}^{*}-a_{m} \mathbf{1}_{n-1}\|_2^2\leq \frac{e^{2\kappa_E}}{c^2\lambda_2(\mathcal{L}_{A_{-m}})^2}\|\nabla\ell_n^{(-m)}(\bbrtheta_m^*)\|_2^2.
\end{equation*}
Where $A_{-m}$ is the adjacency matrix of the comparison graph with node $m$
excluded. By the interlacing property of the eigenvalue sequences of Laplacians
of graph and its induced subgraph \citep[see, e.g.][Proposition
3.2.1]{spectra2012}, we have $\lambda_2(\mathcal{L}_{A_{-m}})\geq
\lambda_2(\mathcal{L}_{A})$, thus
\begin{equation}\label{eq:bound_l2_subset}
    \|\bbrtheta_{-m}^{(m)}-\bbrtheta_{-m}^{*}-a_{m} \mathbf{1}_{n-1}\|_2^2\leq \frac{e^{2\kappa_E}}{c^2\lambda_2(\mathcal{L}_{A})^2}\|\nabla\ell_n^{(-m)}(\bbrtheta_m^*)\|_2^2.
\end{equation}
Now by \Cref{lm:lem7}, it holds with probability at least $1 - O(n^{-10})$ that
\begin{equation*}
    \|\bbrtheta_{-m}^{(m)}-\bbrtheta_{-m}^{*}-a_{m} \mathbf{1}_{n-1}\|_2^2\leq C\frac{e^{2\kappa_E }n_{\max} n}{L \lambda_2(\mathcal{L}_{A})^2},
\end{equation*}
and the conclusion is guaranteed by a union bound.
\eprfof


}{}

\end{document}